\documentclass[11pt]{article}

\usepackage{amscd}
\usepackage{amstex}
\usepackage{theorem}
\usepackage{epsf}

\title{\bf Homeomorphisms between Limbs of the Mandelbrot Set}

\author{Bodil Branner\\  The Tech.~Univ.~of Denmark \\ Building 303 \\ DK-2800 Lyngby, Denmark\\ {\tt branner\@mat.dtu.dk} \and  N\'uria Fagella \\ Math.~Sci.~Research Inst.\\ 1000 Centennial Drive \\ Berkeley, CA 94720, USA \\ {\tt nuria\@math.bu.edu}}

\newtheorem{theorem}{Theorem}[section]
\newtheorem{lemma}[theorem]{Lemma}
\newtheorem{proposition}[theorem]{Proposition}
\newtheorem{definition/lemma}[theorem]{Definition/Lemma}
\newtheorem{corollary}[theorem]{Corollary}
{\theoremstyle{break} \theorembodyfont{\rmfamily} \newtheorem{remarks}[theorem]{Remarks}}
\newenvironment{remark}{\vspace{\baselineskip} \noindent{\bf Remark \addtocounter{theorem}{1}\thetheorem \ \nobreak}}{$\bullet$}
\newenvironment{definition}[1]{\noindent{\newline \bf Definition #1\ \nobreak}}{}
\newenvironment{theorem*}[1]{\par \vspace{\baselineskip}\noindent{ \bf Theorem #1 \nopagebreak}\par \em}{\vspace{\baselineskip}}
\newenvironment{proof}[1]{\noindent{\bf Proof #1\/: }}{\boxend}
\newcommand{\boxend}{\\ \hspace*{5.5in} {\bf q.e.d.}}
\newcommand{\fullref}[1]{\ref{#1} on page~\pageref{#1}}

\newcommand{\la}{\lambda}
\newcommand{\lla}{$\lambda$}
\newcommand{\pol}{P_{q,\lambda}}
\newcommand{\ppol}{$P_{q,\lambda}$}

\newcommand{\potil}{\widetilde{P}_\nu}
\newcommand{\potiltil}{\widetilde{\widetilde{P}}_\nu}
\newcommand{\potilq}{\widetilde{P}_{q,\nu}}
\newcommand{\form}{$P_{q,\lambda}(z)=\lambda z (1+z/q)^q$}
\newcommand{\ho}{\phi_{p/q}}
\newcommand{\bfihat}{\widehat{\Phi}}
\newcommand{\fihat}{\widehat{\phi}}
\newcommand{\hho}{$\phi_{p/q}$}
\newcommand{\hoi}{\xi_{p/q}}

\newcommand{\finho}{\Phi_{p p^\prime}^q}
\newcommand{\finte}{\Theta_{p p^\prime}^q}

\newcommand{\qu}{$Q_c(z)=z^2+c$}
\newcommand{\qc}{$Q_c $}

\newcommand{\Vtil}{\widetilde{V}}
\newcommand{\vtil}{\widetilde{v}}
\newcommand{\latil}{\widetilde{\lambda}}
\newcommand{\Latil}{\widetilde{\Lambda}}

\newcommand{\estil}{\widetilde{S}}

\newcommand{\rtil}{\widetilde{R}}
\newcommand{\ftil}{\widetilde{f}}
\newcommand{\pitil}{\widetilde{\Pi}}
\newcommand{\Ttil}{\widetilde{T}}

\newcommand{\Ltil}{\widetilde{L}}
\newcommand{\psitil}{\widetilde{\psi}}
\newcommand{\Ktil}{\widetilde{K}}
\newcommand{\ktil}{\widetilde{k}}
\newcommand{\otil}{\widetilde{\Omega}}
\newcommand{\pp}{p^\prime}
\newcommand{\thetatil}{\widetilde{\theta}}

\newcommand{\bigthat}{\widehat{\Theta}}
\newcommand{\bigthatp}{{\widehat{\Theta}}_{p/q}}
\newcommand{\bigthatpp}{{\widehat{\Theta}^q_{p {p^\prime}}}}

\newcommand{\linte}{\overline{\Theta}}
\newcommand{\goesto}{\stackrel{1-1}{\longrightarrow}}
\newcommand{\goestoo}{\stackrel{2-1}{\longrightarrow}}

\newcommand{\cala}{{\mathcal A}}

\newcommand{\bc}{{\Bbb  C}}
\newcommand{\bt}{{\Bbb   T}}

\newcommand{\bi}{{\Bbb   I}}
\newcommand{\bd}{{\Bbb   D}}

\newcommand{\br}{ {\Bbb   R}}
\newcommand{\bq}{ {\Bbb   Q}}
\newcommand{\bh}{ {\Bbb   H}}
\newcommand{\bn}{ {\Bbb   N}}
\newcommand{\bz}{ {\Bbb   Z}}
\newcommand{\ical}{{\mathcal I}}
\newcommand{\mcal}{{\mathcal M}}
\newcommand{\icalhat}{\widehat{\mathcal I}}
\newcommand{\icalp}{{\mathcal I}_{p/q}}
\newcommand{\one}{^{(1)}}
\newcommand{\two}{^{(2)}}
\newcommand{\three}{^{(3)}}
\newcommand{\p}{_{p/q}}

\newcommand{\qo}{_{q,0}}
\date{}

\begin{document}

\maketitle

\vspace{0.5cm}
\centerline{\bf Abstract}
\begin{quotation} {\small
Using a family of higher degree polynomials as a bridge, together with complex surgery techniques, we construct a homeomorphism between any two limbs of the Mandelbrot set of common denominator. Induced by these homeomorphisms  and complex conjugation we obtain an involution between each limb and itself, whose fixed points form a topological arc.  All these maps have counterparts at the combinatorial level relating corresponding external arguments. Assuming local connectivity of the Mandelbrot set we may conclude that the constructed homeomorphisms between limbs are compatible with the embeddings of the limbs in the plane. As usual we plough in the dynamical planes and harvest in the parameter space.}

\end{quotation}

\newpage

\tableofcontents
\newpage
\reversemarginpar
%*****************************************************************************
%*****************************************************************************
%*****************************************************************************

\section{Introduction } \label{intro}  

The Mandelbrot set  has been a main object of interest and study in the recent years. It is associated to  the quadratic family of complex polynomials \qu, for $c\in \bc$. Let $K_c$ denote the {\em filled Julia set} of $Q_c$, i.e. 
\[ 
	K= K_c=\{z \in \bc  \mid \{Q^n_c(z)\}_{n \in \bn} \mbox{ is bounded }\},
\] 
and let $J=J_c$ be the boundary of $K_c$, the {\em Julia set} of $Q_c$. 

As for all rational functions, the dynamical behavior of the critical point $\omega=0$ dominates the dynamical behavior of the polynomial. The filled  Julia set $K_c$ is connected when the orbit of $0$ is bounded and totally disconnected when it is unbounded.

The Mandelbrot set, $M$, is defined as the set of parameter values $c$, for which $K_c$ is connected, or equivalently as the set of parameters for which the orbit of $0$ is bounded.

\begin{figure}[htbp]
%\vspace{2cm}
\centerline{\epsfysize=9cm\epsffile{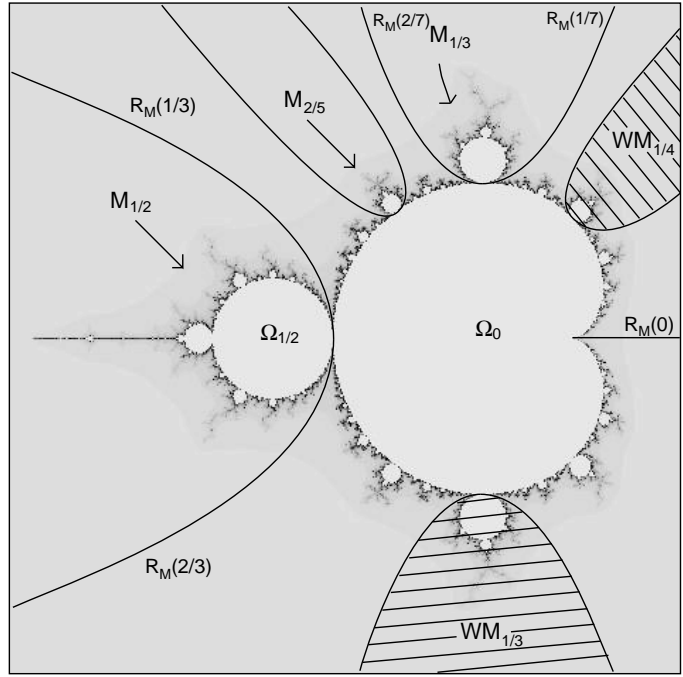}}
\caption{\small The Mandelbrot set. }
\label{mandelbrot}
\end{figure}

The works of  A. Douady, J.~H. Hubbard \cite{dh} \cite{dh2}, D. Sullivan \cite{sullivan} and J-C. Yoccoz in the last decade contributed enormously to the understanding of the Mandelbrot set, but there remain many interesting open questions. The main one is to prove that $M$ is
locally connected (MLC conjecture). The following are some well known results about $M$:
\begin{enumerate}
\item $M$ is full, compact and connected.

\item The interior of $M$ contains  connected components for which $Q_c$ has an attracting periodic orbit. These are called {\em hyperbolic components} and it is  conjectured that their  union equals the interior of $M$.

\item For each hyperbolic component $ \Omega$, there is a conformal isomorphism $\rho_\Omega: \bd \rightarrow \Omega$ such that $Q_c$ has a cycle of multiplier $t$ when $c=\rho_\Omega(t)$. The point $\rho_\Omega(0)$ is called the {\em center} of $\Omega$.  This map extends to a homeomorphism of $\overline{\bd}$ onto $\overline{\Omega}$. The function $\gamma(t)=\gamma_\Omega (t)=\lim_{r \rightarrow 1} \rho_\Omega (r e^{2\pi i t})$ defines a parametrization of the boundary of $\Omega$. For each $t\in \br / \bz$, the point  $\gamma(t)$ in $\partial \Omega$ is said to have {\em internal argument} $t$. The point $\gamma(0)$ is called {\em the root} of the hyperbolic component $\Omega$.

\item There is a unique hyperbolic component,  namely $\Omega_0$, that is bounded by the main cardioid,  for which $Q_c$ has an attracting fixed point . For any internal argument $p/q \in \bq / \bz$, there is a hyperbolic component $\Omega_{p/q}$ attached to the cardioid at the point $\gamma_0(p/q)$. This hyperbolic component contains  $c$-values for which $Q_c$ has an attracting periodic cycle of period $q$. The point $\gamma_0(p/q)$ is the root  of the hyperbolic component $\Omega_{p/q}$. 

\item There exists a conformal isomorphism $\varphi_M: \bc \setminus M \longrightarrow \bc \setminus \overline{\bd}$ such that $\frac{\varphi_M(c)}{c} \rightarrow 1$ when $c \rightarrow \infty$. An {\em external ray} of {\em external argument}  $\theta$ is defined by
\[ 
	R_M(\theta)=\varphi_M^{-1}(\{re^{2\pi i \theta}\}_{1<r<\infty}).  
\]  
(See Figure \ref{mandelbrot} and Sect.~\ref{rays}). Rays with rational arguments play a special role. They are known to have a limit on the boundary of $M$ as $r\rightarrow 1$. The set of {\em landing points} of these rays consists of all the roots of the hyperbolic components and all the {\em Misiurewicz points}, i.e. $c$ values for which $\omega=0$ is strictly preperiodic.
\end{enumerate}

Let $p$ and $q$ be positive integers, $p<q$, $q \geq 2$ and $\gcd(p,q)=1$.

We define the  $p/q${\em -limb} of $M$, $M_{p/q}$, to be the connected component of $M\setminus {\overline{\Omega}}_0$ attached to the main cardioid at the point with internal argument $p/q$, i.e. at $c=\gamma_0(p/q)$  (see Figure \ref{mandelbrot}).

We define the $p/q${\em -wake} of $M$, $WM_{p/q}$, to be the open subset of $\bc$ that contains the $p/q$-limb of $M$ and is bounded by the two rays landing at $c=\gamma_0(p/q)$ union this point \cite{atela}.

\begin{definition}{}
Given two compact sets $K$ and $K'$ in $\bc$, we say that a homeomorphism $\phi: K \rightarrow K'$ is {\em compatible } (resp.~{\em reversely compatible}) {\em with the embeddings} of $K$ and $K'$ in $\bc$ if there are neighborhoods $N$ and $N'$ of $K$ and $K'$ respectively such that $\phi$ extends to a homeomorphism $\widehat{\phi}:N \rightarrow N'$, preserving (resp.~reversing) orientation.
\end{definition}

Some years ago, J.~C. Yoccoz observed (unpublished) that the limbs $M_{1/5}$ and $M_{2/5}$ are homeomorphic and that corresponding polynomials have conjugated dynamics on the filled Julia sets. Recently E.~Lau and D.~Schleicher \cite{dierk} have announced the existence of homeomorphisms between any two limbs with common denominator compatible with the dynamics, but not compatible with the embeddings of the limbs in the plane.

In this paper, we construct a homeomorphism between any two limbs with common denominator. The homeomorphism is not compatible with the dynamics but compatible with the embeddings assuming local connectivity of  $M$. In fact, we believe that the hypothesis of local connectivity could be removed but we will discuss this fact in a forthcoming paper.

Our goal is to prove the following theorems:

\begin{theorem*}{A} \nobreak
 Given $p/q$ and $\pp/q$, there exists a homeomorphism 
\[ 	
	\Phi=\finho :M_{p/q} \longrightarrow M_{\pp/q} 
\] 
which is holomorphic on the interior of $M_{p/q}$.

{\bf Moreover}, for each $c \in M_{p/q}$, there exists a homeomorphism
\[ 
	\bfihat_c: K_c \longrightarrow K_{\Phi(c)} 
\] 
which is holomorphic in the interior of $K_c$ and compatible with the embeddings of $K_c$ and $K_{\Phi(c)}$ in $\bc$.  
\end{theorem*}

\begin{figure}[htbp]
%\vspace{2cm}
\centerline{\epsfysize=7cm\epsffile{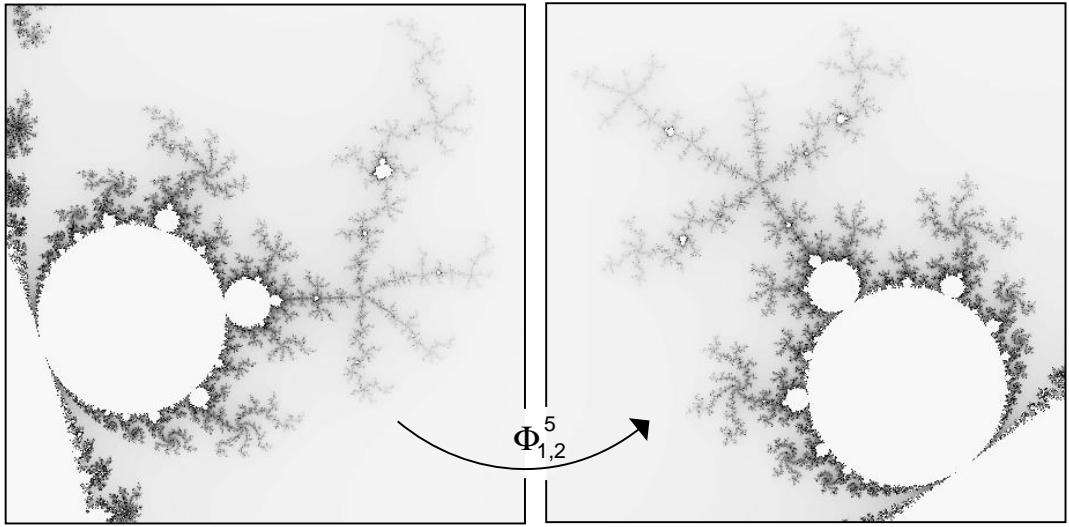}} 
\caption{\small The 1/5-limb and the 2/5-limb of the Mandelbrot set. }
\label{fifths}
\end{figure}

\begin{figure}[htbp]
%\vspace{2cm}
\centerline{\epsfysize=7cm\epsffile{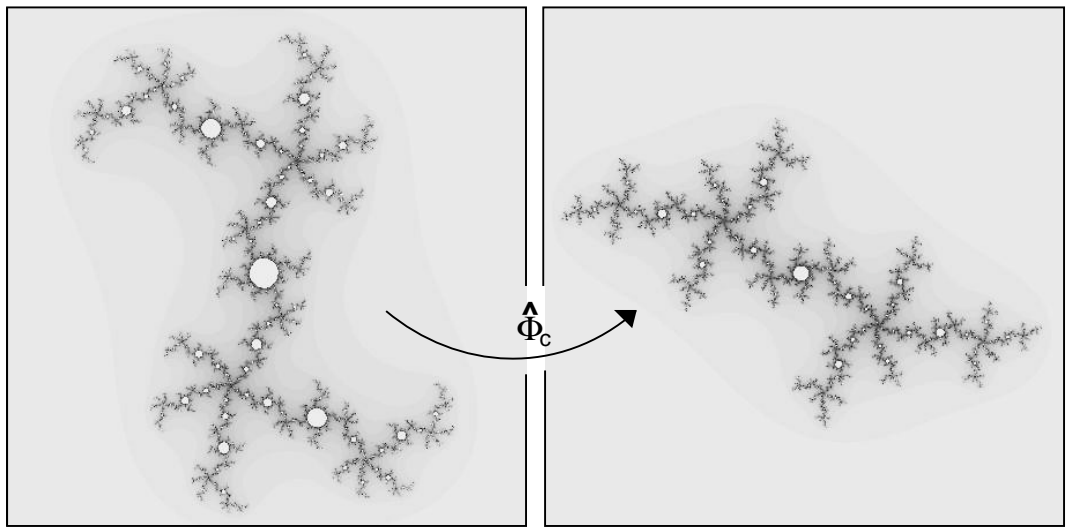}} 
\caption{\small Corresponding filled Julia sets under $\widehat{\phi}_c$. Left: $c  \in M_{1/5}$, center of a period 6 hyperbolic component. Right: $c \in M_{2/5}$, center of a period 8 hyperbolic component. }
\label{juliafifths}
\end{figure}

\begin{figure}[htbp]
%\vspace{2cm}
\centerline{\epsfysize=7cm\epsffile{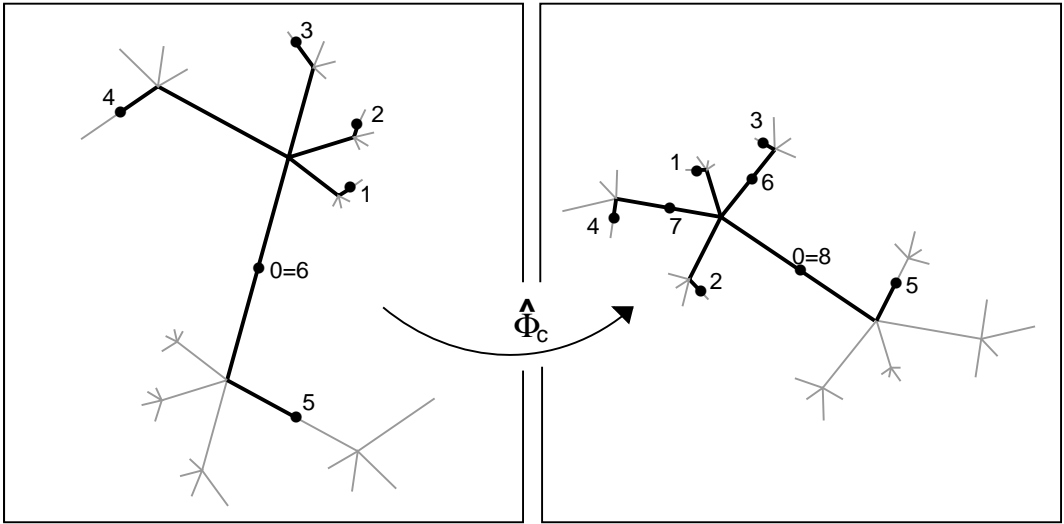}} 
\caption{\small Hubbard trees for the filled Julia sets in Fig.~\ref{juliafifths} }
\label{treesfifths}
\end{figure}

See Figs.~\ref{fifths}, \ref{juliafifths} and \ref{treesfifths}. In fact, a stronger statement is true: for all $c \in M\p$ the homeomorphism $\widehat{\Phi}_c$ matches with the combinatorially defined map described in Theorem C.

\begin{theorem*}{B} \nobreak 
Given $p/q$, there exists a homeomorphism,
\[ 
	\ical=\ical_{p/q} :M_{p/q} \longrightarrow M_{p/q}. 
\] 
which is an involution, i.e.~$\ical^2={\mbox Id}$, and is anti-holomorphic in the interior of $M_{p/q}$.  The set of points that are fixed by $\ical$ form a topological arc through the limb.

{\bf Moreover}, for each $c\in M_{p/q}$, there exists a homeomorphism
\[ 
	\widehat{\ical }=\widehat{\ical }_c: K_c \longrightarrow K_{\ical (c)} 
\] 
reversely compatible with the embeddings, which is anti-holomorphic in the interior of $K_c$. 
\end{theorem*}

\begin{figure}[htbp]
%\vspace{2cm}
\centerline{\epsfysize=7.5cm\epsffile{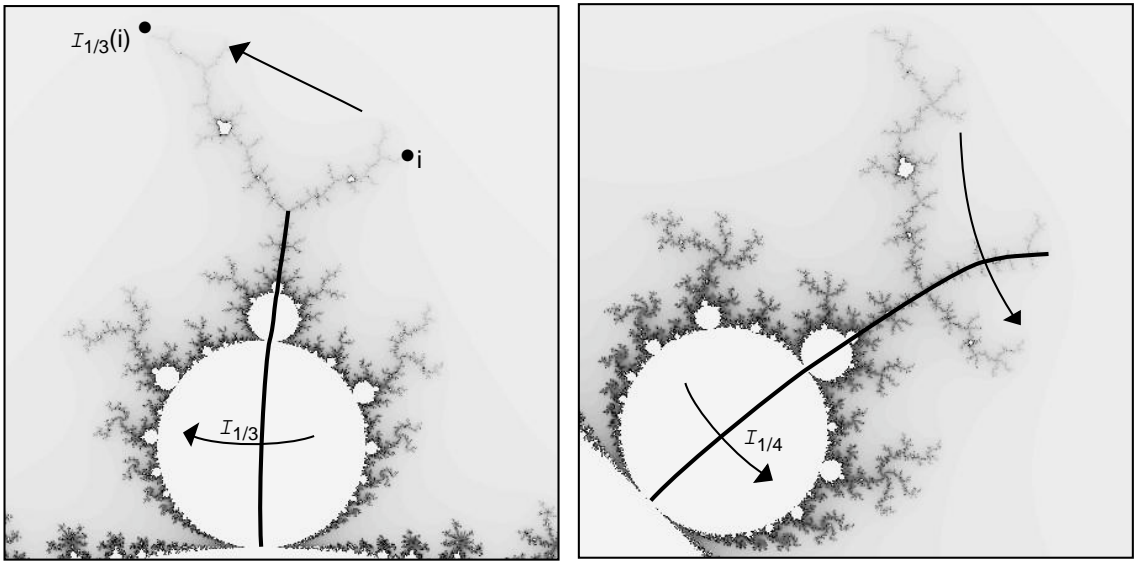}}
\caption{\small Arc of symmetry for the involutions $\ical _{1/3}$ and $\ical _{1/4}$.} 
\label{inv1}
\end{figure}

\begin{figure}[htbp]
%\vspace{2cm}
\centerline{\epsfysize=7.5cm\epsffile{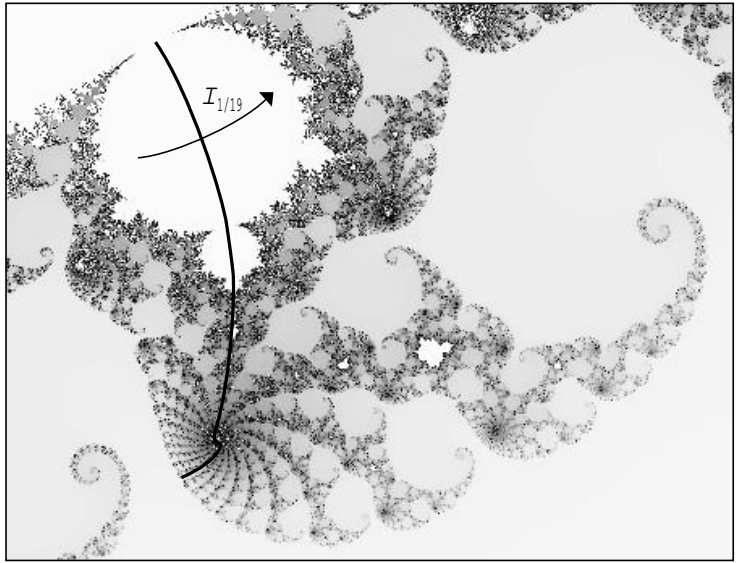}}
\caption{\small Arc of symmetry for the involution $\ical _{1/19}$.} 
\label{inv1-19}
\end{figure}

\begin{figure}[htbp]
%\vspace{2cm}
\centerline{\epsfysize=7cm\epsffile{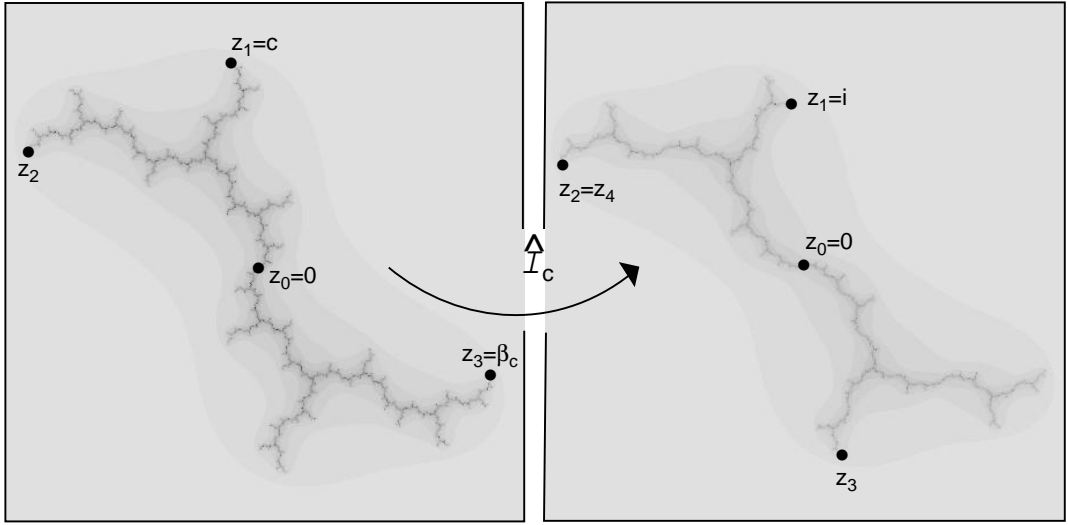}}
\caption{\small The Julia sets for the landing point of the ray $R_M(1/4)$ and $\ical _{1/3}(c)=i$.}
\label{julinv}
\end{figure}

See Figs.~\ref{inv1}, \ref{inv1-19} and \ref{julinv}.

\begin{theorem*}{C} Given $p/q$ and $\pp/q$, let $\theta^\pm_{p/q}$ and $\theta^\pm_{\pp/q}$ be the arguments of the external rays landing at the root points of the limbs $M_{p/q}$ and $M_{\pp/q}$ respectively. Then, there exists an orientation preserving homeomorphism
\[ 
	\finte: [\theta^-_{p/q},\theta^+_{p/q}] \longrightarrow [\theta^-_{\pp/q},\theta^+_{\pp/q}] 
\] 
such that
\begin{enumerate}
\item For $\theta \in [\theta^-_{p/q},\theta^+_{p/q}] \cap \bq$, the ray $R_M(\theta)$ lands at a point $c \in M_{p/q}$ if and only if $R_M(\finte(\theta))$ lands at $\finho(c) \in M_{\pp/q}$.

\item The map $\finte$ is induced by a homeomorphism
\[ 
	\bigthat=\bigthatpp: \bt \longrightarrow \bt 
\] 
such that for any $\theta \in \bt$ and any $c \in M_{p/q}$, the ray $R_c(\theta)$ lands at $z\in K_c$ if and only if the ray $R_{\finho(c)}(\bigthat(\theta))$ lands at $\widehat{\Phi}_c(z) \in K_{\finho(c)}$.
\end{enumerate}
\end{theorem*}

It is remarkable that the map $\finte$ in Theorem C does not depend on the point $c \in M_{p/q}$. A stronger statement than part one of Theorem C is in fact true (see remark \ref{stronger}).
 
\begin{theorem*}{D} 
Given $p/q$ and $\theta^\pm_{p/q}$ as in Theorem C, there exists an orientation reversing homeomorphism $\linte=\linte_{p/q} $ of $[\theta^-_{p/q},\theta^+_{p/q}]$ onto itself such that
\begin{enumerate}
\item $\linte^2=$ Id.

\item For $\theta \in [\theta^-_{p/q},\theta^+_{p/q}] \cap \bq $, the ray $R_M(\theta)$ lands at a point $c \in M_{p/q}$ if and only if $R_M(\linte(\theta))$ lands at $\ical _{p/q}(c) \in M_{p/q}$,

\item The argument $\theta^s$ fixed by $\linte$ is rational.

\item The map $\linte\p$ is induced by an orientation reversing homeomorphism
\[
	\widehat{\linte}=\widehat{\linte}\p: \bt \longrightarrow \bt
\]
such that given $c\in M_{p/q}$ and $\theta \in \bt$, the ray $R_c(\theta)$ lands at a point $z \in K_c$, if and only if the ray $R_{\ical (c)}(\widehat{\linte}(\theta))$ lands at $\widehat{\ical }(z) \in K_{\ical (c)}$.
\end{enumerate}  
\end{theorem*}

\begin{figure}[htbp]
%\vspace{2cm}
\centerline{\epsfysize=7.5cm\epsffile{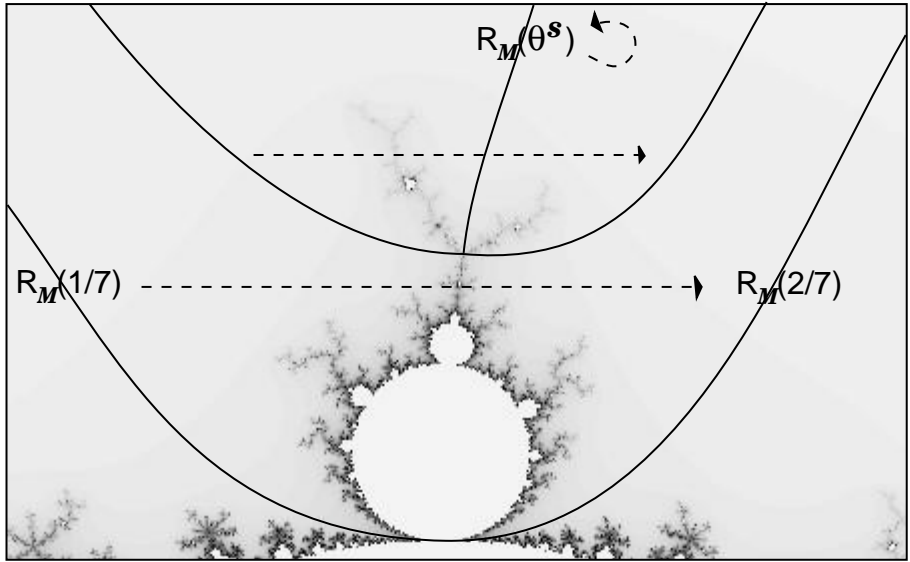}}
\caption{\small Some rays with corresponding arguments under the map $\linte_{1/3}$.}
\label{raysinv}
\end{figure}

\begin{theorem*}{E}
Assume $M$ is locally connected. Then  the map $\finho$ (resp.~$\ical\p$) is compatible (resp.~reversely compatible) with the embeddings of the limbs in $\bc$. 
\end{theorem*}

To prove Theorem A, we use complex surgery techniques which are analogous to those used in \cite{Branner}, where a homeomorphism was constructed between the $1/2$-limb of M and a limb in a cubic parameter space. The way we use this technique  is somewhat unexpected, since we use different families  of higher degree polynomials, $P_{q,\lambda}=\la z (1+z/q)^q$, as a bridge for the construction of the homeomorphisms. That is,  we actually show that each of the $p/q$-limbs  of $M$ is homeomorphic to the 0-limb $L_{q,0}$ to be defined in section \ref{polynomials}. These families of polynomials were introduced in \cite{me} as a way to approximate the family $z \mapsto \la z e^z$, inspired by \cite{dgh} where they considered polynomial families approximating $z \mapsto \la e^z$. 

The homeomorphisms could be defined directly without reference to the higher degree polynomials. But the similarity between corresponding quadratic polynomials becomes more transparent when passing to these special higher degree polynomials. Moreover the involutions together with the arc of symmetry on the limbs are immediate from the symmetry of the 0-limb $L_{q,0}$ with respect to the real axis. Although the involutions also follow from the symmetry of $M$ with respect to the real axis together with the homeomorphisms in Theorem A, it is not clear that the fixed points form a topological arc.

The statements given above for theorems A to D just mention existence of homeomorphisms. In this paper we actually construct such homeomorphisms. The homeomorphisms constructed are not compatible with the dynamics.

The paper is structured as follows:

Instead of proving Theorem A-D directly we prove analogous theorems F-H, forming the bridges to higher degree polynomials. In Sect.~\ref{prelims} we recall some general results which we use throughout the paper. In Sect.~\ref{polynomials} we introduce the families of higher degree polynomials that we use as bridges, and formulate Theorem F, the parameter space analogue of Theorem A. In Sect.~\ref{dynchar} we emphasize the common dynamical features to be applied. Sect.~\ref{surgery} is by far the most important part of the paper. Using surgery techniques we prove Theorem F and also Theorem G, the dynamical analogue of Theorem A. Theorem C follows immediately. Sect.~\ref{involutions} is dedicated to Theorem B. In the last section we state and prove Theorem H, the combinatorial analogue of Theorem C. Theorem C-E follow.

We thank Bob Devaney, Alfredo Poirier, Dan S{\o}rensen, Dierk Schleicher and Tan Lei for many helpful conversations and especially, we thank Adrien Douady not only for innumerable comments and suggestions but also for his constant encouragement. We also thank Christian Mannes for creating the program that produced all the pictures in this paper. We wish to thank the Mathematical Sciences Research Institute at Berkeley, the Mathematics Department at Boston University and the Mathematical Institute at the Technichal University of Denmark for their hospitality. Research at MSRI was supported in part by NSF grant DMS-9022140.

%*****************************************************************************
%*****************************************************************************
%*****************************************************************************

\section{Preliminaries} \label{prelims}
%*****************************************************************************

\subsection{Dynamics of Polynomials} \label{dynpol}

If $f$ is a  complex polynomial of degree $d>1$,  the point at infinity is always a superattracting fixed point. We define its filled Julia set as
\[	K_f=\{ z \in \bc \mid \{f^n(z) \}_{n\in \bn} \text{ is bounded} \}, \]
and its Julia set $J_f$ as the boundary of $K_f$.

\begin{proposition} \label{dichotomy}
Let $f$ be a polynomial of degree $d$. Then, $K_f$ is connected if and only if the orbits of all critical points of $f$ are bounded.
\end{proposition}
For a proof see \cite{Blanchard}

Assume $f$ is monic and $K_f$ connected. Then, there exists a unique analytic isomorphism tangent to the identity at infinity,  $\psi_f: \bc \setminus \overline{\bd}  \longrightarrow \bc \setminus K_f$, which conjugates $f_0(z)=z^d$ to $f$. The map $\psi_f$ is called the B\"{o}ttcher coordinate.

We can also lift $f_0$ to a map $\widetilde{f}_0$ on the universal covering space, the right half plane $\bh$. If we require that  $\widetilde{f}_0(\br_+)=\br_+$ we have that
\[ \widetilde{f}_0(z)=\mcal_d(z) :=d z\]
In summary, the following diagram commutes
\[
  \begin{CD}
	\bh   @>{\widetilde{f}_0}>> \bh \\
	 @V{\exp}VV @VV{\exp}V   \\
   \bc\setminus\overline{\bd} @>f_0>> \bc\setminus\overline{\bd}\\
	@V\psi_fVV @VV\psi_fV   \\
   \bc\setminus K_f @>f>> \bc\setminus K_f					 
	\end{CD}
\]
	
We remark that in the case of $K_f$ being locally connected, $\psi_f$ extends continuously to the boundary of $\bd$, so that $\psi_f \circ \exp$ is defined on $\overline{\bh}$. Even in the case when $K_f$ is not locally connected, there is a  set of points of full measure on $\partial \bd$ where $\psi_f$ is well defined. This set always includes the points with rational arguments.

Given $t\in \br$ we denote by $R(t)$ the horizontal line in $\bh$ with imaginary part equal to $2 \pi t$ i.e. 
\[ R(t):=\{\rho+2\pi i t \in \bh \mid \rho>0\}. \]
We define the {\em external ray of argument} $t$ to be 
\begin{eqnarray*} 
	R_f(t)&=&\psi_f(\exp(R(t)))\\
	         &=& \psi_f( \{(re^{2\pi i t} )\mid r\in (1,\infty)\}). 
\end{eqnarray*}   
If $R_f(t)$ has a limit when  $r \rightarrow1$ (or $\rho \rightarrow 0$), then it  tends to a point of the Julia set which we denote by $R_f^\ast(t)$. We say that $R_f(t)$ {\em lands} at $R_f^\ast(t)$ and we have:
\[ 
	f(R_f^\ast ( t))=R_f^\ast(d t). 
\] 
All external rays with rational arguments land and if $K_f$ is locally connected all external rays land.

\begin{figure}[htbp]
%\vspace{2cm}
\centerline{\epsfysize=4cm\epsffile{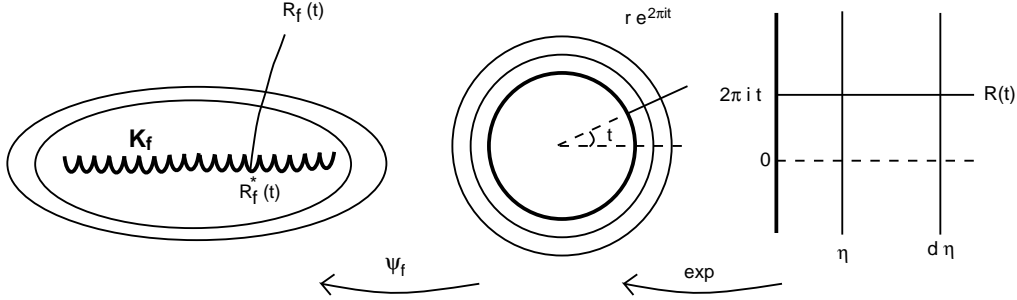}}
\caption{\small Sketch of an external ray and some level sets of the potential function.}
\label{julia}
\end{figure}

The {\em potential}, $G_f$,  (Green's function) of $K_f$ satisfies $G_f=\log(|\psi_f^{-1}|)$ on $\bc \setminus K_f$. We may extend it to $\bc$ defining $G_f=0$ on $K_f$. We have that 
\[ 
	G_f(f(z))=d \; G_f(z) 
\]
for all $z\in \bc $.  The positive level sets of the potential function are simple closed curves around $K_f$ called {\em equipotentials}, and the potential  measures the rate of escape of  points under iteration.  See Figure \ref{julia}.  A given equipotential $\{ z\in \bc\setminus K_f \mid G_f(z)=\eta\}$ corresponds in the exterior of the unit disc to a circle around the origin of radius $e^\eta$ and on the right half plane to a vertical line of real part $\eta$.

\subsubsection{The Quadratic family} \label{rays}

Let $Q_c(z)=z^2+c$ and  set $K_c=K_{Q_c}$, $\psi_c=\psi_{Q_c}$, etc. Let $M$ be the Mandelbrot set (see Sect.~\ref{intro}). The results in this section can be found in \cite{dh2}.

The map $\varphi_M : \bc \setminus M \longrightarrow \bc \setminus \overline{D}$ defined as $\varphi_M(c)=\psi_c^{-1}(c)$ is a conformal isomorphism. So we define an {\em external ray} of {\em external argument}  $\theta$ as
\[ 
	R_M(\theta)=\varphi_M^{-1}(\{re^{2\pi i \theta}\}_{1<r<\infty}).  
\]  
If $R_M(\theta)$ has a limit $c \in \partial M$ when $r \rightarrow 1$ we say that $R_M(\theta)$ {\em lands} at $c$.  We denote the landing point of the ray $R_M(\theta)$ by $R_M^\ast (\theta)$. 

The conjecture that the Mandelbrot set is locally connected is equivalent to the continuous landing of all external rays.

It is known that all external rays with rational arguments land. In what follows we relate the external rays of $M$ with the  external rays in the dynamical plane of the corresponding quadratic polynomials, as well as the relation between the landing points in the case of rational arguments.

\begin{enumerate}
\item Suppose $c \in \bc \setminus M$. Then 
\begin{equation} \label{reading}
	c \in R_M(\theta)  \Longleftrightarrow c \in R_c(\theta)
\end{equation}
\item Suppose $c$ is a Misiurewicz point in $M$. Then 
\[
R_c(\theta) \text{\ lands  at } c\in K_c \Longleftrightarrow R_M(\theta) \text{\ lands at } c \in M.
\]
An argument $\theta$ satisfying the above is rational with even denominator.

\item Suppose $c$ is the center of a hyperbolic component $\Omega$ of $M$ of period $k$. Let $c_0=\gamma_\Omega(0)$ denote the root point of the hyperbolic component. Let $U_c$ denote the connected component of Int$(K_c)$ containing $c$ and let $z_0$ denote the root point on $\partial U_c$, that is the boundary point  which is fixed by $Q_c^k$. Then
\[
R_c(\theta) \text{\ lands  at } z_0 \in K_c \text{\ adjacent to } U_c \Longleftrightarrow R_M(\theta) \text{\ lands at } c_0 \in M.
\]
An argument $\theta$ satisfying the above is rational with odd denominator.
\end{enumerate}

%************************************************************************************   
\subsection{Tools} \label{tools}

For the surgery in Sect.~\ref{surgery} we will use the theory of quasi-conformal mappings, the theory of integrability by Morrey-Ahlfors-Bers and the theory of Polynomial-like mappings of Douady and Hubbard (compare with \cite{Branner}, \cite{dh3}).

\subsubsection*{Quasi-conformal mappings, Beltrami forms and almost complex structures}

Let $X$ and $Y$ be two Riemann surfaces isomorphic either to $\bd$ (the unit disc) or $\bc$.

An {\em almost complex structure} $\sigma$ on $X$ is a measurable field of ellipses $(E_x)_{x \in X}$, equivalently defined by a measurable Beltrami form $\mu$ on $X$
\[ \mu = u \frac{d \overline{z}}{dz}. \]
The correspondence between Beltrami forms and complex structures is as follows: the argument of $u(x)$ is twice the argument of the major axis of $E_x$, and $|u(x)|=\frac{K-1}{K+1}$ where $K\geq 1$ is the ratio of the lengths of the axes. 

The {\em standard complex structure}  $\sigma_0$ is defined by circles.

Suppose $\varphi:X \longrightarrow Y$ is a quasi-conformal homeomorphism. Then $\varphi$ gives rise to an {\em almost complex structure} $\sigma$ on $X$. For almost every $x \in X$, $\varphi$ is differentiable and the $\br$-linear tangent map $T_x \varphi : T_x X \longrightarrow T_{\varphi(x)} Y$ defines, up to multiplication by a positive factor, an ellipse $E_x$ in $T_x X$:
 \[         E_x=(T_x \varphi)^{-1} (S^1).	\]
Moreover, there exists a constant $K>1$ such that the ratio of the axes of $E_x$ is bounded by $K$ for almost every $x \in X$. The smallest bound is called the {\em dilatation ratio} of $\varphi$.

Equivalently,  $\varphi $ defines a measurable Beltrami form on $X$
\[	\mu=\frac{\overline{\partial} \varphi}{\partial \varphi}= \frac{\frac{\partial \varphi}{\partial \overline{z}}}{\frac{\partial \varphi}{\partial z}} \frac{d \overline{z}}{dz}=u(z) \frac{d \overline{z}}{dz}. 
\]

An almost complex structure is quasi-conformally equivalent to the standard structure if it is defined by a measurable field of ellipses with bounded dilatation ratio.

Given $\varphi:X \longrightarrow Y$ a quasi-conformal homeomorphism, an almost complex structure $\sigma$ on $Y$ can be {\em pulled back} into an almost complex structure $\varphi^\ast \sigma$ on $X$. If $\sigma$ is defined by an infinitesimal field of ellipses $(E_y)_{y \in Y}$ then $\varphi^\ast \sigma$ is defined by $(E_x)_{x \in X}$ where $E_x= (T_x \varphi)^{-1} E_{\varphi(x)}$ whenever defined.

To {\em integrate} an almost complex structure $\sigma$ means to find a quasi-conformal homeomorphism $\varphi$ such that $ (T_x \varphi)^{-1} (S^1) = \rho(x) E_x$ for almost every $x \in X$. Informally, we will say that $\sigma$ is {\em transported} to $\sigma_0$ by $\varphi$. 

A {\em quasi-regular} mapping is of the form $\psi=h \circ \varphi$ where $\varphi$ is quasi-conformal and $h$ is holomorphic, but $h$ may have critical points.

The Theorem of Integrability can be stated as follows:
\begin{theorem*}{of Integrability (Measurable Riemann Mapping Theorem)}
Let $X$ be a Riemann surface isomorphic to $\bd$ or to $\bc$. Let $\sigma_\mu$ be any  almost complex structure on $X$ given by the Beltrami form 
\[ \mu= u \frac{d \overline{z}}{dz} \]
with bounded dilatation ratio, i.e.~
\[ \| \mu \|_\infty := sup |u(z)| <m<1. \]
Then $\sigma_\mu$ is integrable i.e.~there exists a quasi-conformal homeomorphism $\varphi$ such that 
\[	\mu=\frac{\overline{\partial} \varphi}{\partial \varphi}. \]
If $X$ is isomorphic to $\bd$ then $\varphi:X \longrightarrow \bd$ is unique up to composition with a M\"{o}bius transformation mapping $\bd$ to $\bd$.  If $X$ is isomorphic to $\bc$ then  $\varphi: X \longrightarrow \bc$ is unique up to composition with an affine map.
\end{theorem*}

In Sect.~\ref{continuity} we will need the dependence on parameters of the integrating maps, when dealing with Beltrami forms defined in $\bc$. 

\begin{theorem*}{of Integrability (Dependence on Parameters)}
Let  $\Lambda$ be an open set of $\bc$ and let  $(\mu_\la=u_\la \frac{d \overline{z}}{dz})_{\la \in \Lambda} $ be  a family of measurable Beltrami forms on $\bc$. Suppose that $\la \mapsto u_\la(z)$ is holomorphic for each fixed $z \in \bc$, and that there is $m <1$ such that  $\| \mu_\la \|_\infty<m$ for all $\la \in \Lambda$. Let $\varphi_\la :\bc \rightarrow \bc$  be the unique  quasi-conformal homeomorphisms that fixes two given points in $\bc$ and such that $\mu_\la=\frac{\overline{\partial} \varphi_\la}{\partial \varphi_\la}$. Then for each $z \in \bc$ the map $\la \mapsto \varphi_\la(z)$ is holomorphic.
\end{theorem*}

\subsubsection*{Polynomial-like mappings}

A {\em polynomial-like mapping}  is a proper holomorphic mapping 
\[ f: U' \longrightarrow U \]
where $U'$ and $U$ are open sets in a Riemann surface $X$ isomorphic to a disc, and $\overline{U'} \subset U$. The {\em degree} of $f$ is the number of preimages of any point $z \in U$ counted with multiplicity. 

For a polynomial-like mapping $f$ we define
\[ K_f= \{ z \in U' \mid f^n(z) \in U' \text{\ for all } n \in \bn \}. \]

Given two polynomial-like mappings $f:U' \rightarrow U$ and $g: V' \rightarrow V$ with $K_f$ and $K_g$ connected we say that $f$ and $g$ are {\em holomorphically equivalent} (resp.~{\em quasi-conformally equivalent}) if there exists a holomorphic (resp.~quasi-conformal) homeomorphism
\[ \varphi:U_1 \longrightarrow V_1 \]
where $U_1$ and $V_1$ are neighborhoods of $K_f$ and $K_g$ conjugating $f$ to $g$ on $U'_1=f^{-1}(U_1)$ i.e.~satisfying $g \circ \varphi = \varphi \circ f$ on $U'_1$.

A {\em hybrid equivalence} is a quasi-conformal equivalence which satisfies $\overline{\partial} \varphi=0$ almost everywhere on $K_f$  so which is holomorphic in the interior of $K_f$. 

If $f$ and $g$ are polynomials, a {\em co-hybrid equivalence } between $f$ and $g$ is a quasi-conformal equivalence $\varphi: \bc \rightarrow \bc$ which is holomorphic on $\bc \setminus K_f$.

\vspace{\baselineskip}

\noindent{\bf The Straightening Theorem}
{\em Let  $f:U^\prime \longrightarrow U$ be a polynomial-like mapping of degree $d$ with $K_f$ connected. Then, $f$ is hybrid equivalent to a polynomial $P$ of degree d. Moreover, this polynomial is unique up to an affine conjugacy.}

\vspace{\baselineskip}

We will also need the following result:

\begin{proposition} \label{cor2}
Let $f$ and $g$ be polynomials with $K_f$ and $K_g$ connected. If $f$ and $g$ are hybrid equivalent then they are conjugate by an affine map.
\end{proposition}

%*****************************************************************************
%*****************************************************************************
%*****************************************************************************

\section{Families of Higher Degree Polynomials} \label{polynomials} 

Consider the family of  complex polynomials  \form, where $\lambda \in \bc \setminus \{0\}$ and $q \geq 2$.  Each \ppol\  has degree $q+1$ and has exactly  two critical points: $-q$ with multiplicity $q-1$, and $\omega=\frac{-q}{q+1}$ with multiplicity one. Moreover,  $\pol (-q)=0$ and $0$ is fixed. We now show that polynomials with these features, can always be expressed in the form $P_{q,\la}(z)$.
 
\begin{definition}{} A polynomial $P(z)$ of degree $q+1 \geq 3$  is of type $ E_q$ if it satisfies the following properties:
\begin{enumerate}
\item $P$ has two critical points: $-q$ of multiplicity $q-1$ and $\omega \in \bc$ of multiplicity one.
\item $P$ has a fixed point at $z=0$.
\item $P(-q)=0$.
\end{enumerate}
\end{definition}

\begin{proposition} \label{esubq}
 Any polynomial $P(z)$ of degree $q+1$ which  is of type $E_q$ is of the form
\[
	 P_{q,\la}(z)=\la z (1+\frac{z}{q})^q .
\] 
Moreover, if $P_{q,\la}$ and $P_{q,\la^\prime}$ of type $E_q$ are affine conjugate with $q \geq 3$, then $\la=\la^\prime$.
\end{proposition}

From here and throughout the paper, we will omit the subindex $q$ whenever it causes no confusion.
 
\begin{proof}{} 
The polynomial $P$ has a zero of order $q$ at the point $-q$ and a simple zero at $0$ and, since $P$ is of degree $q+1$, it is of the form $\la z (1+z/q)^q$.

For $q \geq 3$, if an affine map $A$ conjugates $P_{q,\la}$ to $P_{q,\la'}$ it must fix $-q$ (because it is the only multiple critical point for both polynomials) and $0$ (because it is the image of $-q$), so $A$ is the identity.
\end{proof} 
 
Since the critical point of high multiplicity is always prefixed, the dynamics of \ppol\  are dominated by the behavior of the free  critical point, $\omega$. We have the following lemma, which follows from Prop.~\ref{dichotomy}:

\begin{lemma} 
Let $K_{q,\la}$ be the filled Julia set of \ppol. Then, $K_{q,\la}$ is connected if and only if the orbit of $\omega$ is bounded.
\end{lemma}

\begin{proposition} \label{choices}
Let $P_{q,\la}$ and $P_{q,\la'}$ be  polynomials of type $E_q$ with $K_\la$ and $K_{\la'}$ connected. If they are hybrid equivalent, then $\la=\la'$.
\end{proposition}

\begin{proof}{}
This is a direct corollary of Prop.~\ref{cor2}. By  this proposition $P_{q,\la}$ and $P_{q,\la'}$ are affine conjugate. By the Prop.~\ref{esubq} $\la=\la'$.
\end{proof}

\begin{figure}[htbp]
%\vspace{2cm}
\centerline{\epsfysize=7cm\epsffile{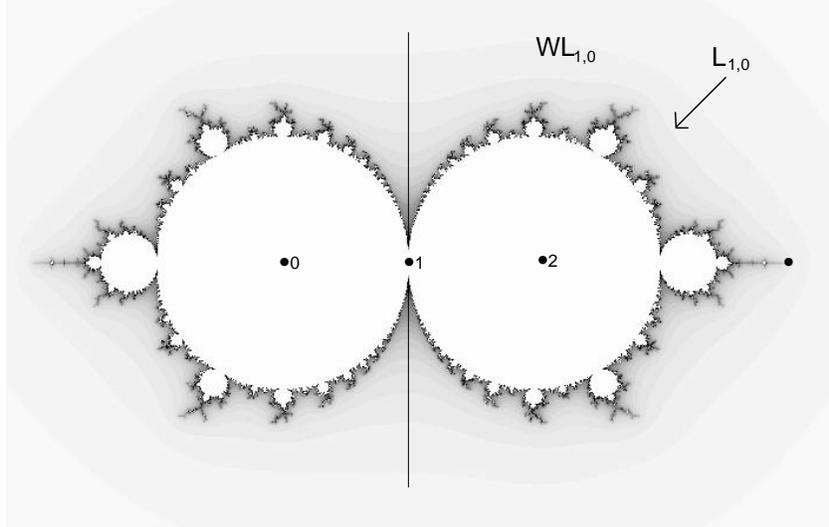}}
\caption{\small Parameter plane $\Lambda_1$: the Connectedness Locus $L_1$ and the 0-wake. Range: $[-2.5,9.14]\times[-2.1,2.1]$. }
\label{pol1}
\end{figure}

\begin{figure}[htbp]
%\vspace{2cm}
\centerline{\epsfysize=8cm\epsffile{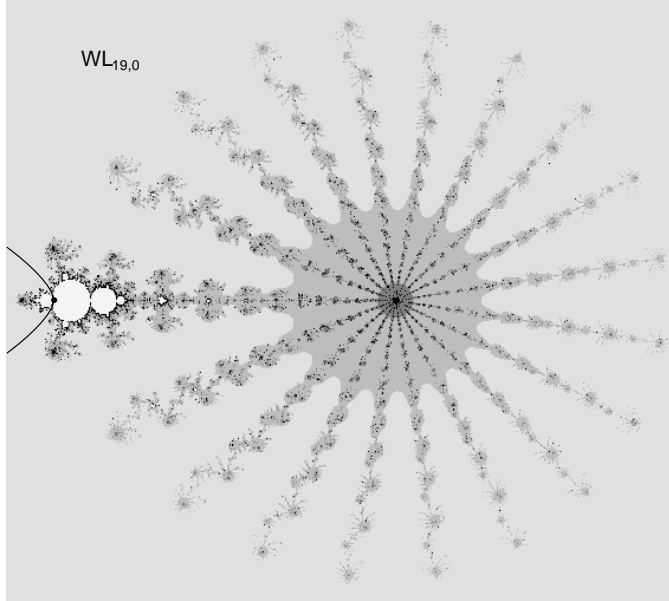}}
\caption{\small Parameter plane $\Lambda_{19}$: the Connectedness Locus $L_{19}$ and the $0-wake$. Range: $[-6,95]\times [-45.6,45.6]$.}
\label{pol19}
\end{figure}

For each $q\in \bn$ let $\Lambda_q$ be $\bc^\ast$ viewed as the parameter space of \ppol. Let $L_q$ denote the {\em connectedness locus}, i.e. the set of $\la \in \Lambda_q$  for which $K_{q,\la}$ is connected.  Observe that the set $L_q$ associated with \ppol\  is the analogue of the Mandelbrot set for \qc.  Figures \ref{pol1}, through \ref{pol19} show the sets $L_1$, $L_2$, $L_3$, $L_4$, $L_5$ and  $L_{19}$, respectively.  Following this analogy, we will call a connected component $\Omega$ of the interior of $L_q$ a {\em hyperbolic component} if the free critical point $\omega$ is attracted to an attracting cycle (despite the fact that the polynomial is not hyperbolic, since the orbit $\{-q,0\}$ is contained in the Julia set, see Sect.~\ref{monicfamily}). As in the Mandelbrot set, a hyperbolic component $\Omega$ is parametrized by the multiplier map $\rho_\Omega: \bd \longrightarrow \Omega$, a biholomorphic map that extends continuously to the boundary. For all $q$, the sets $L_q$ share the following features:

\begin{enumerate}
\item If  $0<|\la|<1$, $z=0$ is an attracting fixed point of \ppol, with $\omega$ in its basin of attraction. Hence the punctured closed unit disk is contained in $L_q$.

\item Let $r/s \in \bq / \bz$ such that $r \geq 0$, $s>0$ and $\gcd(r,s)=1$. At $\la_{\frac{r}{s}}=e^{2\pi i \frac{r}{s}}$, the polynomial $\pol$ experiences a period $s$-tupling bifurcation. That is, when we exit the unit disk through $\la_{\frac{r}{s}}$, $z=0$ becomes repelling and a cycle of period $s$ becomes attracting. Hence, there is a hyperbolic component of period $s$ attached to the unit disk at $\la_{\frac{r}{s}}$.
\end{enumerate}

\begin{figure}[htbp]
%\vspace{2cm}
\centerline{\epsfysize=7cm\epsffile{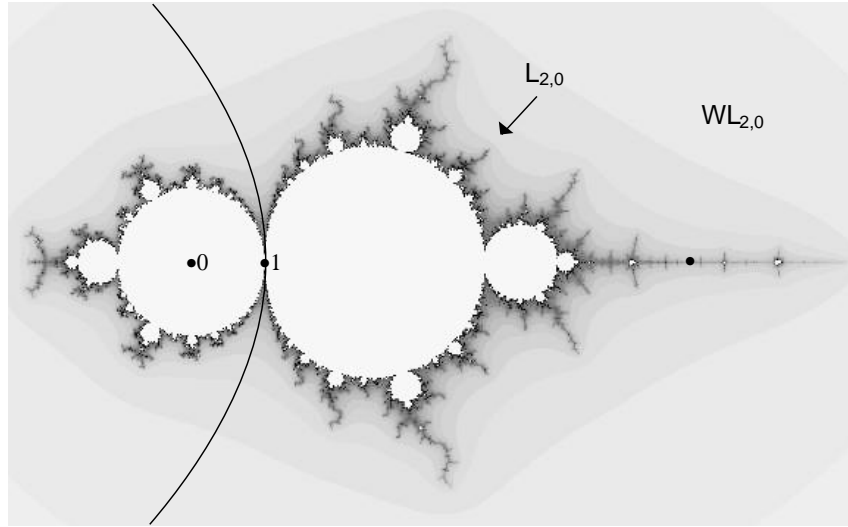}}
\caption{\small Parameter plane $\Lambda_2$: the Connectedness Locus $L_2$ and the o-wake. Range: $[-2.2,4.4]\times[-3.5,3.5]$. }
\label{pol2}
\end{figure}

\begin{figure}[htbp]
%\vspace{2cm}
\centerline{\epsfysize=7cm\epsffile{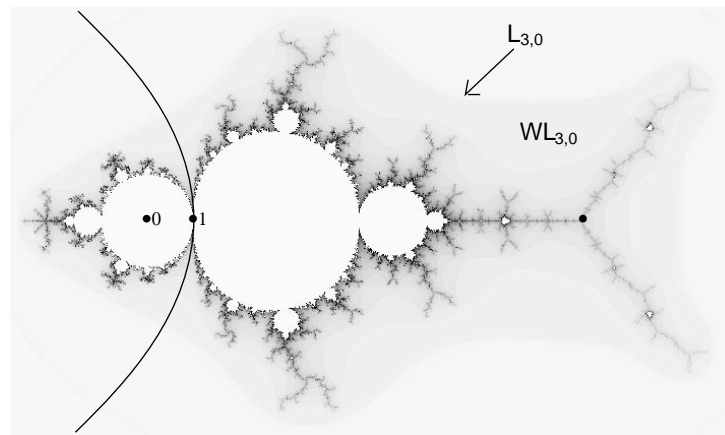}}
\caption{\small Parameter plane $\Lambda_3$: the Connectedness Locus $L_3$ and the 0-wake. Range: $[-3,13]\times [-4.7,4.7]$. }
\label{pol3}
\end{figure}

\begin{figure}[htbp]
%\vspace{2cm}
\centerline{\epsfysize=8cm\epsffile{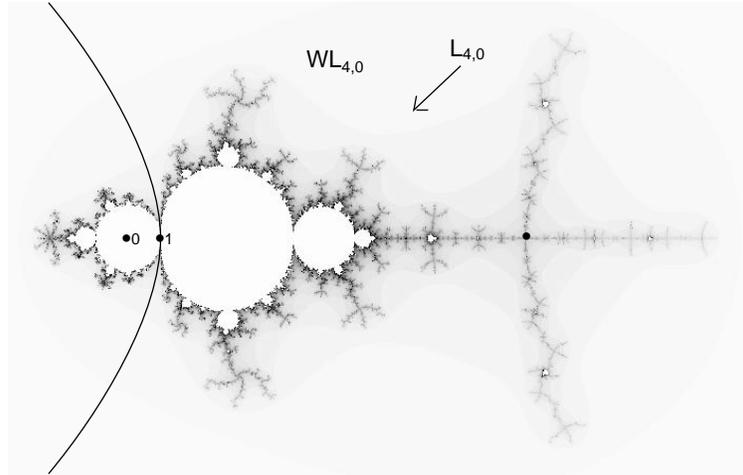}}
\caption{\small Parameter plane $\Lambda_4$: the Connectedness Locus $L_4$ and the 0-wake. Range: $[-3.5,23.4]\times [-9.5,9.5]$.}
\label{pol4}
\end{figure}

\begin{figure}[htbp]
%\vspace{2cm}
\centerline{\epsfysize=8cm\epsffile{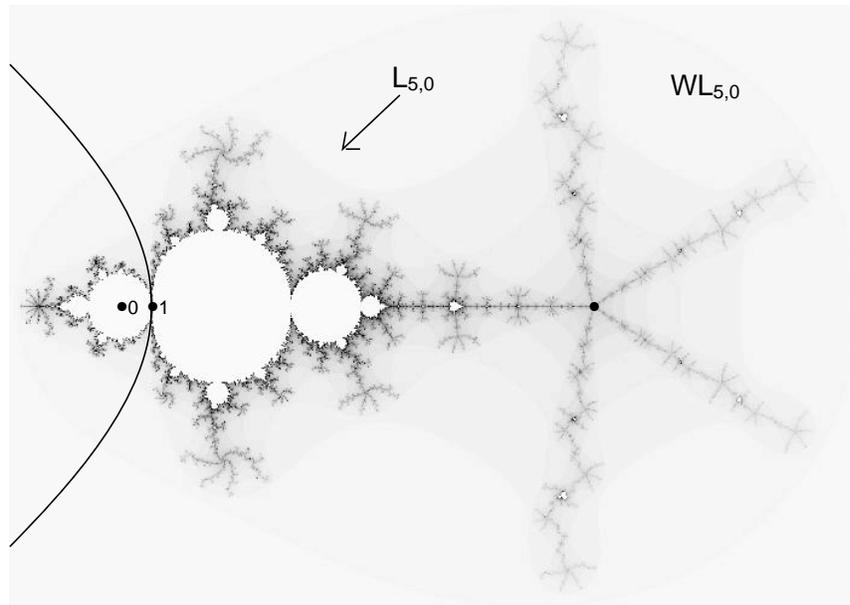}}
\caption{\small Parameter plane $\Lambda_5$: the Connectedness Locus $L_5 and the 0-wake$.}
\label{pol5}
\end{figure}

For $L_q$, the unit circle is the analogue to the main cardioid of $M$. Here, we do not have a cusp and a hyperbolic component of period $1$ is attached to $\la_{\frac{0}{1}}=1$. 

As we did for $M$, we define the $\frac{r}{s}$-{\em limb of } $L_q$, $L_{q,\frac{r}{s}}$, to be the connected component of $L_q \setminus \overline{\bd}$ attached to the unit circle at the point $\la_{\frac{r}{s}}$. The $0$-limb, $L_{q,0}$,  will be of special importance in our work.

One can define external rays in $\bc \setminus L_q$ similarly to $\bc \setminus M$, to be explained in section \ref{extrays}. There are exactly two rational rays landing at $\la=1$, the root of 0-limb. Hence we may define the 0-{\em wake} $WL_{q,0}$ of  $L_q$, as the open subset of $\bc$ that contains the $0$-limb of $L_q$ and is bounded  by the two rays landing at its root. In section \ref{polyn}, the 0-wake is discussed in detail.

Theorem A is now a corollary of the following:

\begin{theorem*} {F}  Let $p$ and $q$ be positive integers such that $p<q$, $q\geq 3$ and $\gcd(p,q)=1$. Then,  there exists a homeomorphism
\[		
	\phi=\ho: M_{p/q} \longrightarrow L_{q,0}	,
\]
which is holomorphic in the interior of $M_{p/q}$. \newline 
\end{theorem*}

Note that we have restricted to $q \geq 3$. The reason is that the connectedness locus $L_1$ is a ramified double covering of $M$ and not related to limbs of $M$; the 0-limb $L_{2,0}$ of $L_2$ is homeomorphic to $M_{1/2}$ \cite{Branner} and gives no new information with respect to Theorems A to E.

%*****************************************************************************
%*****************************************************************************

\subsection{Monic Families of Higher Degree Polynomials and External Rays} \label{monicfamily}

Our goal in this section is to make precise definitions of external rays in the dynamical plane of $\pol$ as well as in the parameter plane. Also, we characterize the polynomials $\pol$ with $\la$ in the 0-wake. 

As explained in Sect.~\ref{dynpol}, when $P$ is a monic polynomial, the ray $R_P(\theta)$ is define using the B\"{o}ttcher coordinates  $\psi_P$ which is tangent to the identity at infinity. Here we are dealing with $\pol$ which is not monic, hence there is no canonical choice of a map conjugating $\pol$ to $z \mapsto z^{q+1}$ on the complement of the filled Julia set. This will cause some minor complication.

Set 
\[	\potil(z)=\potilq(z)=z (z+\nu)^q.\]
The polynomial $\potilq$ is affine conjugate to $\pol$ with $\la=\nu^q$ through the map $g_\nu(z)=\frac{q}{\nu}z$. Under this map, the special fixed point $0$ is still at $0$, the multiple critical point is at $-\nu$, the free critical point is at $\widetilde{\omega}_\nu=\frac{-\nu}{q+1}$ and the free critical value is at 
\begin{equation} \label{critvalue}
\vtil_\nu=\potil(\frac{-\nu}{q+1})= -\frac{q^q}{(q+1)^{q+1}} \; \nu^{q+1}.
\end{equation}
Let $\Ktil_\nu=K_{\potil}$ and $\psitil_\nu:\bc \setminus \overline{\bd} \longrightarrow \bc \setminus \Ktil_\nu$ denote the B\"{o}ttcher coordinates for $\potil$, i.e.~conjugating $\potil$ to $z \mapsto z^{q+1}$ and tangent to the identity at infinity.

For $\la \in L_q \setminus \br_-$ let $\psi_\la :\bc \setminus \overline{\bd} \longrightarrow \bc \setminus K_\la$ denote the unique B\"{o}ttcher coordinate conjugating $P_\la$ to $z \mapsto z^{q+1}$ and being asymptotic to $z \mapsto \frac{q}{\la^{1/q}} z$ at $\infty$, where $()^{1/q}$ denotes the principal branch of the $q$-th root. In other words $\psi_\la=g_\nu \circ \psitil_\nu$ with $\nu=\la^{1/q}$.

We may now define the external ray of argument $t$ to be 
\begin{eqnarray*}
	R_\la(t)&=& \psi_\la(\exp(R(t))) \\
		    &=& \psi_\la (\{ r e^{2\pi i t} \mid r \in (1,\infty) \}) 
\end{eqnarray*}
and denote $G_\la= \log |\psi_\la^{-1}|$  the potential function.

\begin{proposition} \label{0wake}
A parameter  $\la$ belongs to the  0-wake $WL\qo$ if and only if $z=0$ is a repelling fixed point of $\pol$ and the landing point of a fixed ray.
\end{proposition}

In order to show this proposition we will need to define the external rays in parameter space. We will work again with the monic polynomials.

\begin{figure}[htbp]
%\vspace{2cm}
\centerline{\epsfysize=7cm\epsffile{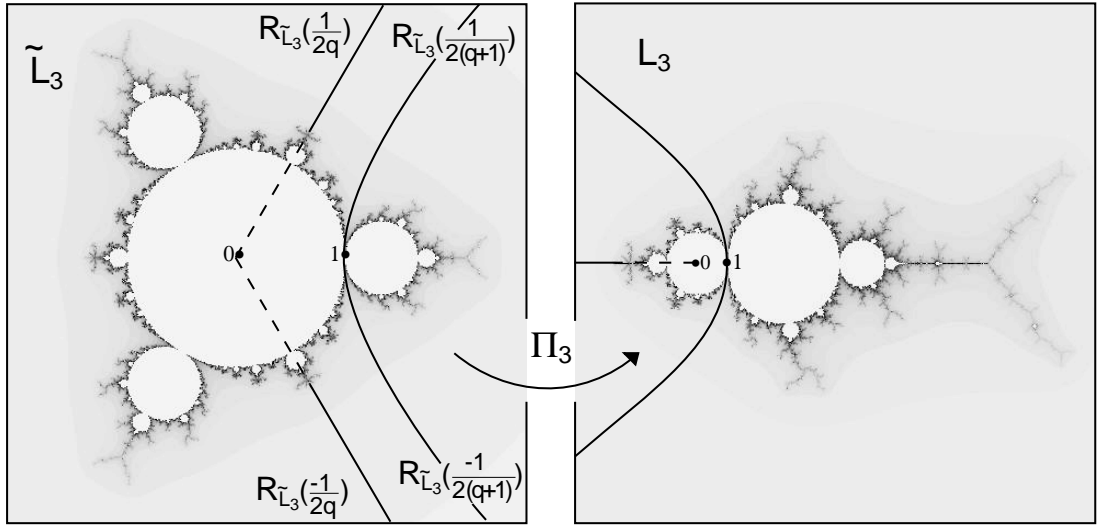}}
\caption{\small The connectedness loci $\Ltil_3$ and $L_3$. }
\label{monic}
\end{figure}

 Let $\Ltil_q$ denote the connectedness locus of the monic family (see Fig. \ref{monic}), i.e.
\[
	\Ltil_q=\{ \nu \in \bc \mid \Ktil_\nu \text{\ connected} \}.
\]
Set $\widetilde{\Lambda}_q=\bc^\ast$ and let $\Pi=\Pi_q:\widetilde{\Lambda}_q \longrightarrow \Lambda_q$ denote the $q$-fold covering map defined by $\Pi(\nu)=\nu^q$. Note that $\Pi$ maps the wedge $\Sigma_q=\{\nu \in \widetilde{\Lambda}_q \mid |\text{Arg}(\nu)| < \frac{\pi}{2q}\}$ conformally onto $\Lambda_q \setminus \br_-$  and in particular, the 0-limb $\Ltil\qo$ onto the 0-limb $L\qo$.

Let $\psi_{\Ltil_q}:\bc \setminus \Ltil_q \longrightarrow \bc \setminus \overline{\bd}$ denote the unique Riemann mapping asymptotic to $\nu \mapsto \ktil \nu$ at $\infty$, where $\ktil$ is some positive real number. The relation between the dynamical maps and the parameter maps is contained in the next proposition.
\begin{proposition}
Let $\nu \in \bc \setminus \Ltil_q$. Then
\begin{equation} \label{noname}
	 \psitil_\nu^{-1} ( \vtil_\nu)=-(\psi_{\Ltil_q}(\nu))^{q+1} .
\end{equation}
\end{proposition}

\begin{proof}{}
Note that the map 
\[ 	\begin{array}{rcl}
	\bc \setminus \Ltil_q  & \longrightarrow  &\bc \setminus \overline{\bd} \\
	\nu & \longmapsto & \psitil_\nu^{-1} ( \vtil_\nu)
	\end{array}
\]
is a $(q+1)$-covering map and that the asymptotic behavior as $\nu \rightarrow \infty$ is given by eq. (\ref{critvalue}). Moreover the map 
\[ 	\begin{array}{rcl}
	\bc \setminus \Ltil_q  & \longrightarrow  &\bc \setminus \overline{\bd} \\
	\nu & \longmapsto & (\psi_{\Ltil_q}(\nu))^{q+1}
	\end{array}
\]
is a $(q+1)$-covering map with the asymptotic behavior as $\nu \rightarrow \infty$ given by $\nu \mapsto \ktil^{q+1} \nu^{q+1}$. The proposition follows.
\end{proof}

\begin{remark}
From the proof we also obtain $\ktil=\frac{q^{\frac{q}{q+1}}}{q+1}$. Note that the radius of capacity of $\Ltil_q$ is $\frac{1}{\ktil}$ so it depends on $q$.
\end{remark}

Instead of proving prop.~\ref{0wake} directly, we prove the analogous proposition in the closed wedge $\overline{\Sigma_q} \subset \widetilde{\Lambda}_q$.

\begin{proposition}

\begin{enumerate}
\item The external rays $R_{\Ltil}(\pm \frac{1}{2(q+1)})$ land at $\nu=1$, the root point of the 0-limb $\Ltil\qo$.

\item A parameter $\nu \in \overline{\Sigma_q}$ belongs to the  0-wake $W\Ltil\qo$ if and only if $z=0$ is a repelling fixed point for $\potilq$ and the landing point of $R_\nu(0)$.
\end{enumerate}
\end{proposition}

\begin{proof}{}
Note that if the ray $R_\nu(0)$ does not branch then it must land at a fixed point. The situation is stable in $\nu$ if the fixed point is repelling. It is unstable either if the fixed point is indifferent (in which case it is parabolic of multiplier 1, due to the Snail lemma, see \cite{milnor}) or if the ray $R_\nu(0)$ branches at the free critical point $\widetilde{\omega}_\nu$. The stable set is open in $\widetilde{\Lambda}_q$ and the unstable set is closed in the same space.

The ray $R_\nu(0)$ passes through the free critical value $\vtil_\nu$ if and only if 
\[
	\text{Arg } (\psitil_\nu^{-1} (\vtil_\nu)) =0
\]
From eq.~(\ref{critvalue}) it follows that this is equivalent to 
\[
	\nu \in R_{\Ltil}(\frac{2n-1}{2(q+1)}), \quad n\in \bz
\]
i.e., restricted to the wedge $\overline{\Sigma_q}$ we have $\nu \in R_{\Ltil}(\pm \frac{1}{2(q+1)})$.

For $\nu=\frac{q+1}{q}$ the polynomial $\potil$ has a superattracting fixed point and $\br_+ \subset \bc \setminus \Ktil_\nu$. By symmetry, $R_\nu(0)=\br_+$ and $R_\nu(0)$ lands at $z=0$. 

There exists $\rho>1$ so that the polynomials $\potil$ for $\nu=\pm \rho \exp(\pi i /q) \in \partial \Sigma_q$ have a superattracting cycle of period two with $\widetilde{\omega}_\nu < 0 < \vtil_\nu$. By symmetry, $R_\nu (0) \subset \br_+$ but $R_\nu(0)$ does not land at $z=0$.

For $\nu$ in the wedge the situation can change only if $\nu=1$ for which 0 is a parabolic fixed point of multiplier 1, or if $\nu \in R_{\Ltil}(\pm \frac{1}{2(q+1)})$ for which the ray $R_\nu(0)$ branches at $\widetilde{\omega}_\nu$.
\end{proof}

\begin{remark}
For completeness we add that the connectedness locus $\Ltil_q$ can be thought of as the connectedness locus of a different family of polynomials of degree $(q+1)$, namely
\[
	\potiltil(z)=z(\nu +z^q).
\]
Note that $\potiltil$ and $\potil$ are semiconjugate by the ramified covering map $\Pi_q(z)=z^q$, i.e.~$\Pi_q \circ \potiltil = \potil \circ \Pi_q$. A polynomial $\potiltil$ is characterized by having a fixed point at $0$ and $q$ ordinary critical points arranged symmetrically around 0. Symmetrical points have symmetrical orbits, that is $\potiltil (e^{\frac{2\pi}{q}i}z)=e^{\frac{2\pi}{q}i} \potiltil(z)$. The symmetry is collapsed by $\Pi_q$ when passing to $\potil$. Note that a hyperbolic component $\Omega$ of $\Ltil_q$ deserves the name hyperbolic when viewed in this other family, i.e.~if $\nu \in \Omega$ then all critical points of $\potiltil$ are attracted to attracting cycles. All our results could be obtained from the family $\potiltil$ with a different kind of surgery.
\end{remark} 

%*******************************************************************************
\subsection{External Rays in the 0-wake} \label{extrays}

To end this section and for later use, we are going to rename the external rays in the 0-wake so that their arguments may be read directly from the dynamical plane as it is the case for the quadratic family (compare with eq.~\fullref{reading}).

Observe that eq.~(\ref{noname}) gives the following result:

If $\nu \in W\Ltil\qo \setminus \Ltil\qo$ then
\begin{equation} \label{noname2}
	\text{Arg}_{\Ktil_\nu}(\vtil_\nu)= [ \frac{1}{2} + (q+1) \text{ Arg}_{\Ltil_q}(\nu)] \pmod 1
\end{equation}
where $\text{Arg}_{\Ktil_\nu}$ (resp.~$\text{Arg}_{\Ltil_q}$) means the external argument with respect to $\Ktil_\nu$ (resp.~$\Ltil_q$).

Inspired by eq.~(\ref{noname2}) we define the affine map
\[	
	\begin{array}{cccc}
		A_q: & (0,1) & \longrightarrow &(-\frac{1}{2(q+1)},\frac{1}{2(q+1)}) \\
			 &   \theta & \longmapsto  & \frac{1}{2} + (q+1) \theta.
	\end{array}
\]
Recall that $\Pi_q : \Latil_q \longrightarrow \Lambda_q$ denotes the covering map $\nu \mapsto \lambda=\nu^q$. External rays of $\Ltil_q$ are mapped by $\Pi_q$ onto external rays of $L_q$. For our purpose, the important external rays of $L_q$ are these in the  0-wake. We parametrize them by arguments in $(0,1)$ as follows:
\[	R_{L\qo}(\theta) = \Pi_q (R_{\Ltil_q} (A_q (\theta))). \]
We say that $\theta$ is the external argument of the ray {\em relative} to the 0-limb.

The following proposition explains why we choose this parametrization:
\begin{proposition} \label{misus}
Let $\la \in L\qo$ be a Misiurewicz point. Then, $R_\la(\theta)$ lands at the critical value $v_\la=\la(-1)^q (\frac{q}{q+1})^{q+1}$ if and only if $R_{L\qo}(\theta)$ lands at $\la$.
\end{proposition}

The proof can be copied from  \cite{dh2} Part I, p.~74-75.

%*****************************************************************************
%*****************************************************************************
%*****************************************************************************

\section{Dynamical Characterizations} \label{dynchar}

In this section, we describe the common dynamical behavior of all polynomials $Q_c$ in a rational wake of $M$ as well as that of all polynomials $\pol$  in the 0-wake of $L_q$. We will use the Julia sets to construct the maps in Theorem F.

%*****************************************************************************

\subsection{$Q_c$ with $c$ in a rational wake} \label{julia of Q} 

Let $K_c=K_{Q_c}$, $\psi_c=\psi_{Q_c}$, etc.
Since $Q_c(-z)=Q_c(z)$, the filled Julia set $K_c$ is always symmetric with respect to $\omega=0$, the critical point.  Let us denote by $\beta_c$ the repelling fixed point where the ray $R_c(0)$ lands. The other preimage of $\beta_c$ under $Q_c$, $-\beta_c$, is the landing point  for $R_c(1/2)$. For $c$ in $\Omega_0$, the other fixed point, $\alpha_c$, is attracting. When $c$ leaves the main cardioid to enter a $p/q$-limb, $\alpha_c$ becomes repelling and a ``pinching'' point in the Julia set. 

Fix $p/q$. For each $c \in M_{p/q}$, in fact for all $c$ in the  wake $WM\p$, there are $q$ rays landing at $\alpha_c$. Their arguments are fixed throughout the wake and form a period $q$ cycle under doubling, with combinatorial rotation number $p/q$. As a consequence they are rational numbers which can be written with denominators $2^q-1$ and numerators depending on $p$ (see \cite{gm}).

\begin{figure}[htbp]
%\vspace{2cm}
\centerline{\epsfysize=7.5cm\epsffile{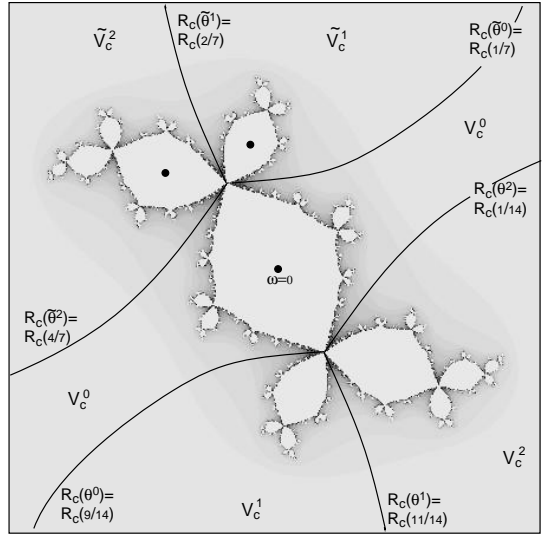}}
\caption{\small The filled Julia set for the center of the main hyperbolic component in $M_{1/3}$ and the five  subsets of the plane. The dots correspond to the period three orbit of the critical point $\omega=0$.}
\label{center13}
\end{figure}

\begin{figure}[htbp]
%\vspace{2cm}
\centerline{\epsfysize=7.5cm\epsffile{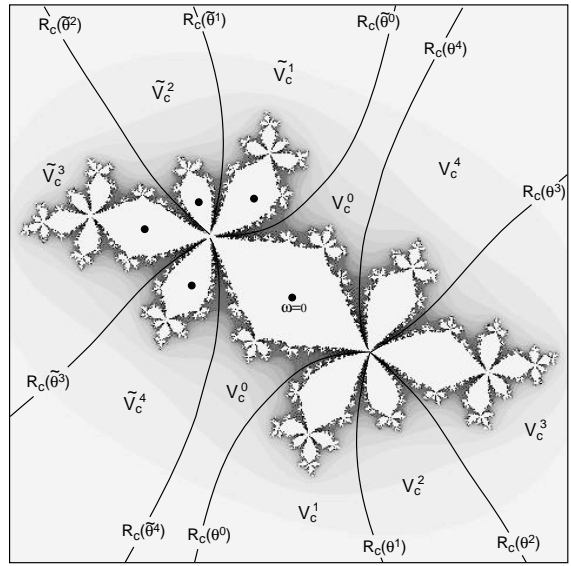}}
\caption{\small The filled Julia set for the center of the main hyperbolic component in $M_{2/5}$ and the nine  subsets of the plane. The dots correspond to the period five orbit of the critical point $\omega=0$.}
\label{center25}
\end{figure}

The other preimage of $\alpha_c$ under $Q_c$ is the point $\alpha^\prime_c=-\alpha_c$. There are $q$ additional rays landing at $\alpha_c^\prime$, and their arguments are  preimages under doubling of the arguments of the rays landing at $\alpha_c$.  Figs.~\ref{center13} and \ref{center25} show examples of Julia sets in the $1/3$ and $2/5$ limbs respectively, together with the rays described above.

The rays landing at $\alpha_c$ and $\alpha^\prime_c$ partition the dynamical plane into $2q-1$ closed subsets. We denote the subset containing the critical point by $V_c^0$, and the others by $V_c^i$ or $\widetilde{V}_c^i=-V_c^i$ for $i=1,2,\ldots q-1$ as shown in Figs.~\ref{center13} and \ref{center25}.  Note that these subsets have their counterparts in the right half plane for which we will not use the subscript $c$. For $0\leq i \leq q-1$ we let $\theta^i \in (0,1)$ be the argument of the ray on the common boundary of $V_c^i$ and $V_c^{[i+1 \pmod q]}$. In the same fashion, $\thetatil^i$ denotes the argument of the ray $R_c(\thetatil^i)=-R_c(\theta^i)$. Note that $R_c(\theta^i)=R_c(\thetatil^i+1/2)$.

Let  $V_c^{0\ldots q-1}=\bigcup_{i=0}^{q-1} V_c^i$.  Then, \qc\  acts on these sets as follows:
\begin{equation} \label{action} 
\begin{array}{llll}   
		V_c^0          	&  \goestoo  & \Vtil^p_c      &    \\
		V_c^i,\Vtil_c^i & \goesto & \Vtil_c^{[i+p \pmod q]} & \mbox{ for $0<i<q-1, i \neq q-p$}\\
      	V^{q-p}_c, \Vtil^{q-p}_c &  \goesto  &  V_c^{0\ldots q-1}& 
\end{array}
\end{equation}    

We establish the following conventions: in the dynamical plane and in expressions with integer indices like $[i+p \! \! \pmod q]$ we will omit $\pmod q$, while in expressions with arguments, we will omit $\pmod 1$. In both cases, it should be  understood that expressions should  be taken $\pmod q$ and $\pmod 1$ respectively.
%*****************************************************************************
\subsubsection{Sectors} \label{the sectors} 

In this part assume $c \in M_{p/q}$. For later purposes, we need to define some subsets which we call {\em sectors}. They should be viewed as neighborhoods of rays $R_c(\theta)$ that land.

\begin{figure}[htbp]
%\vspace{2cm}
\centerline{\epsfysize=4cm\epsffile{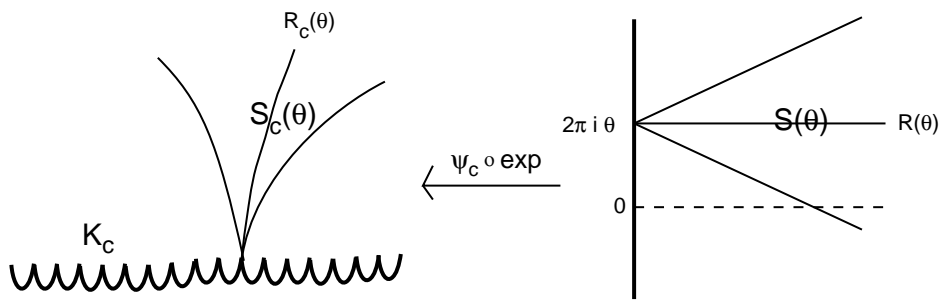}}
\caption{\small Sketch of a sector in the right half plane and its correspondent in the dynamical plane. }
\label{sectors}
\end{figure}

Instead of viewing the sectors in the dynamical plane, it is better to think about them in the exterior of the unit disk or, even better, in the right half plane (see Fig. \ref{sectors}).

\begin{definition}{}
 For a fixed  $s>0$,  we define the sector centered at $R(\theta)$ with slope $s$ as :
\[ 
	S(\theta) = S^s(\theta)= \{\rho + 2\pi i t\in \overline{\bh} \mid |t-\theta|\leq s \rho\}.
\] 
\end{definition}

The boundaries of the sector $S(\theta)$ are straight lines of slope $\pm 2 \pi s$ which cross exactly at $2\pi i\theta $ (see Figure \ref{sectors}).  The doubling map $\mcal_2$ maps each sector $S(\theta)$ homeomorphically onto $S(2\theta)$.

On the dynamical plane of $Q_c$ we define 
\[ S_c(\theta)=S_c^s(\theta)=\psi_c(\exp(S(\theta)). \]

We will only need to consider the sectors $S(\thetatil^i)$ and their preimages under $\mcal_2$ which, in the dynamical plane, correspond to sectors around the rays landing at $\alpha_c$ and their preimages. Since these rays land, $S_c(\theta)$ includes the landing point $R_c^\ast(\theta)$. 

Note that, as they are defined, any two sectors in $\bh$ overlap and hence they  also overlap in $\bc\setminus K_c$. To avoid this we restrict ourselves to a vertical strip as follows.

Choose $\eta >0$ and let $W_n$ denote the vertical strip in $\bh$ with real part less than or equal to  $\frac{\eta}{2^n}$. That is
\[ 
	W_n=\{\rho +2\pi i t \in \bh \mid  \rho \leq \frac{\eta}{2^n}\},
\] 
where $n \in \bn \cup \{0\}$. In dynamical plane, define the {\em filled level set} for the potential function as 
\[ W_{c,n}= \psi_c(\exp(W_n)) \cup K_c.\]

Note that $\mcal_2$ maps $W_{n+1}$ onto $W_n$ homeomorphically while, in the dynamical plane, the polynomial \qc\  maps every $W_{c,n+1}$ onto $W_{c,n}$ with degree $2$. 

\begin{figure}[htbp]
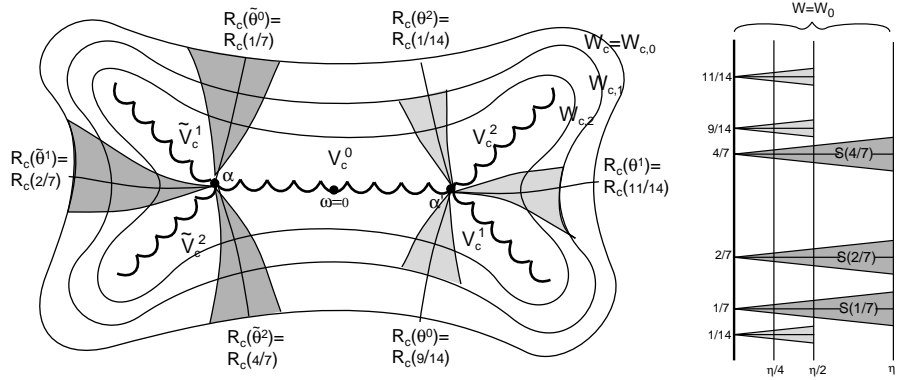

%\vspace{2cm}
\centerline{\epsfysize=5cm\epsffile{pictures/sectorsexamples.ai}
	         \epsfysize=5cm\epsffile{pictures/sectorsexamples2.ai}}
\caption{\small Examples of relevant sectors in the right half plane for the $1/3$-limb and their correspondents in the dynamical plane. }
\label{sectorsexamples}
\end{figure}

We remark that this construction depends on the parameter $c$ only through the function $\psi_c$. Hence, sectors are defined in the same way for all parameter values.

\begin{figure}[htbp]
%\vspace{2cm}
\centerline{\epsfysize=5cm\epsffile{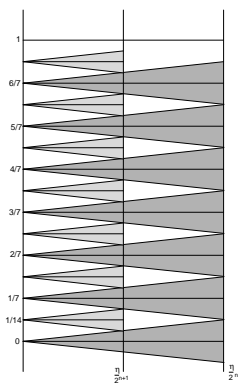}}
\caption{\small $n$-sectors and $(n+1)$ sectors drawn with slope $s=\frac{1}{2\eta(2^q-1)}$. }
\label{limitsectors}
\end{figure}

Set $S_n^s(\theta)=S^s(\theta)\cap \overline{W}_n$. The following proposition assures that the sectors landing at $\alpha_c$ and their preimages do not overlap, if the slope $s$ is chosen sufficiently small (see Fig. \ref{limitsectors}).

\begin{proposition} 
Fix $q\in \bn$, $\eta>0$ and  $0 < s < \frac{1}{2\eta (2^q-1)}$. 
\begin{enumerate}
\item For each $n\in \bn \cup \{0\}$ the $n$-sectors
\[ S_n^s(\frac{m}{2^n(2^q-1)}), \;\;\; m\in \bz \]
are disjoint.
\item The sectors 
\[ S_n^s(\frac{m}{2^n (2^q -1)}) \quad  \text{$m \in \bz$, $m$ odd, $n \in \bn \cap \{0\}$} 
\]
are disjoint.
\end{enumerate}
\end{proposition}

\begin{proof}{} 
The second statement follows from the first. The first statement is proved by induction on $n$.  It is true for $n=0$ since the lines 
\[ 
	t= s \rho + \frac{l}{2^n (2^q-1)}   \quad \text{and} \quad   t= -s \rho + \frac{l+1}{2^n (2^q-1)}  
\] 
do not intersect in $W_0$.

Suppose it is true for $n\geq 0$. Since
\[ \mcal_2 \left [ \bigcup_m S_{n+1}^s (\frac{m}{2^{n+1}(2^q-1)}) \right ] =\bigcup_m S_n^s(\frac{m}{2^n(2^q-1)})
\]
it is true for $n+1$.  
\end{proof}

In figure \ref{limitsectors} the $n$-sectors and the $n+1$ sectors are drawn with slope $s=\frac{1}{2\eta(2^q-1)}$. Two adjacent $n$-sectors intersect at the vertical line $\rho=\frac{\eta}{2^n}$.

%*****************************************************************************
\subsection{\ppol\  with $\la$   in the $0$-wake} \label{polyn} 

In this section we characterize the dynamical behavior of polynomials $\pol$ with $\la$ in the $0$-limb $L_{q,0}$ or, more generally, polynomials in the $0$-wake $WL_{q,0}$.

Recall from sect.~\ref{monicfamily} that  polynomials $\pol$ with $\la$ in $WL_{q,0}$ are characterized by the fixed ray $R_\la(0)$ landing at $z=0$. Since $-q$ is the only preimage of $0$ other than itself, the preimages of $R_\la(0)$ 
\[ 
	R_\la (\frac{1}{q+1}), R_\la (\frac{2}{q+1}), \ldots , R_\la (\frac{q}{q+1})
\] 
land at the critical point $-q$, independently of $\la$. 

\begin{figure}[htbp]
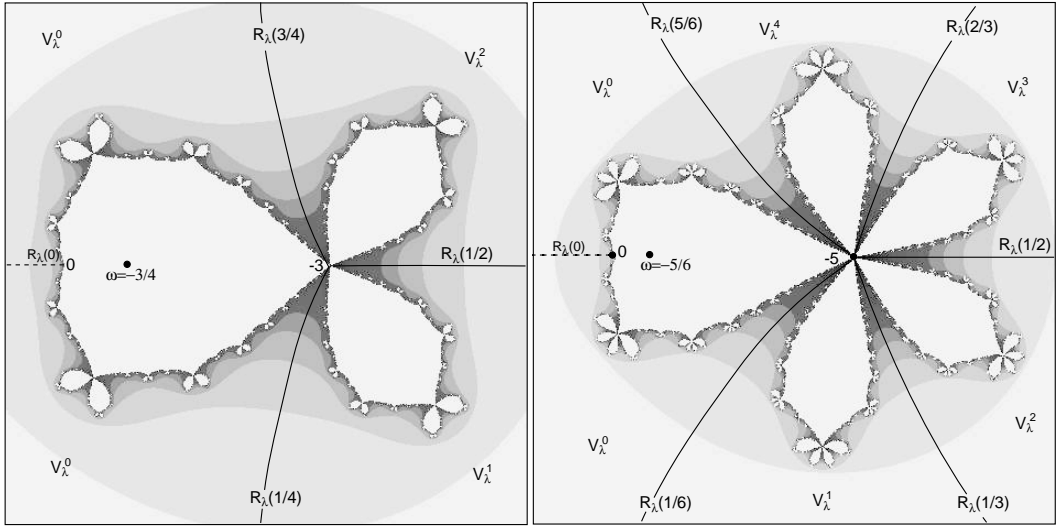

%\vspace{2cm}
\centerline{\epsfysize=7cm\epsffile{pictures/center3.ai}		 \epsfysize=7cm\epsffile{pictures/center5.ai}}
\caption{\small The filled Julia set for the center of the main hyperbolic component in $L_{3,0}$ (left) and $L_{5,0}$ (right) and the three (resp.~five) subsets of the plane. In both cases $\omega$ is fixed. The planes have been rotated 180 degrees. }
\label{center3/center5}
\end{figure}

These rays divide $\bc$ into $q$ closed subsets. Let us denote these sets by $V_\la^0,V_\la^1,\ldots , V_\la^{q-1}$ as shown in Figure~\ref{center3/center5}. We will always draw the dynamical plane for $P_{q,\la}$ rotated by $180^\circ$. The reason will become clear in sect.~\ref{polike}. The polynomial \ppol\  maps each of these subsets as follows:
\begin{equation} \label{action2} 
\begin{array}{lll}   
		\smash[t]{\overset{\circ}{V}}_\la^0   &  \goestoo  & \bc \setminus (R_\la(0) \cup \{0\})          \\
		\smash[t]{\overset{\circ}{V}}_\la^i   &  \goesto  & \bc \setminus (R_\la(0) \cup \{0\})  \quad \text{for $1 \leq i\leq q-1$}
\end{array}
\end{equation}    

As in section \ref{the sectors}, for $\theta \in \br$ let $S(\theta)=S^s(\theta)$ denote a sector in $\overline{\bh}$ with slope $s$. The map $\mcal_{q+1}$ (multiplication by $(q+1)$)  maps $S(\theta)$ homeomorphically onto $S((q+1)\theta)$. When $\la \in L\qo$ we may transport these sectors to the dynamical plane by defining
\[ 
	S_\la(\theta)=\psi_\la(\exp(S(\theta))).
\] 

To avoid the overlapping, fix $\eta>0$ and let $W_n$ be as above, substituting $\frac{\eta}{2^n}$ by $\frac{\eta}{(q+1)^n}$. Let also $W_{q,\la,n}$ be the {\em filled level set} for the potential function, i.e. 
\[
	W_{q,\la,n}=\psi_\la(\exp(W_n)) \cup K_\la
\]  

\begin{figure}[htbp]
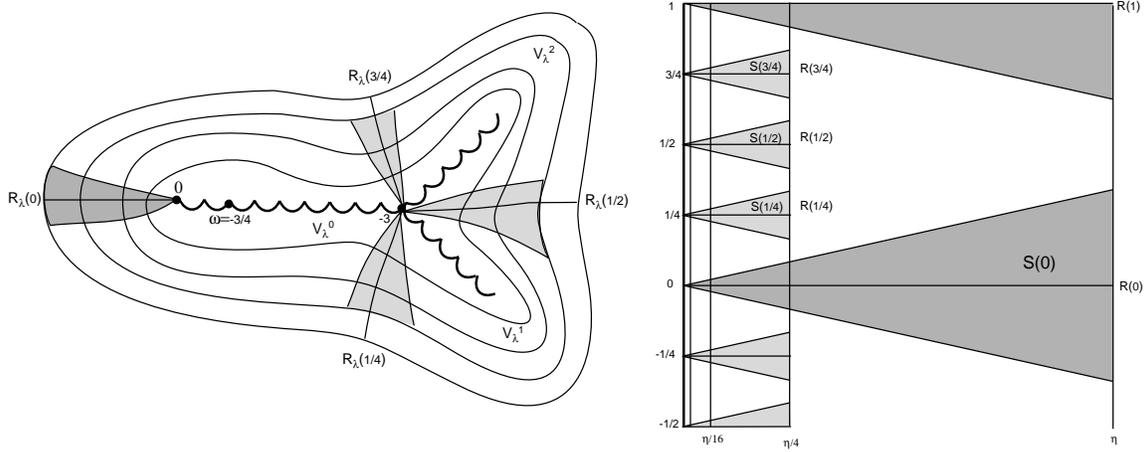

%\vspace{2cm}
\centerline{\epsfysize=6cm\epsffile{pictures/Psectexamples.ai} \epsfysize=6cm\epsffile{pictures/Psectexamples2.ai}}
\caption{\small Sketch of relevant sectors in the right half plane for any $\la \in L_{3,0}$ and their correspondents in dynamical plane. }
\label{Psectexamples}
\end{figure}

We will consider only the sector centered at the $0$-ray and its preimages (see Fig. \ref{Psectexamples}). As in section \ref{the sectors}, one can show that if $0<s< \frac{1}{2 \eta (q+1)}$, these sectors do not overlap for $\rho< \eta$.

%*****************************************************************************
%*****************************************************************************
%*****************************************************************************

\section{Surgery: Proof of Theorem F} \label{surgery}

%*****************************************************************************

\subsection{Idea of the Proof} \label{ideaproof}

We prove Theorem F using {\em surgery } following the methods
in \cite{Branner}. In this section we sketch  the steps of the proof .
 
Let $p$ and $q$ be positive numbers, $p<q$, $q \geq 3$ and $\gcd(p,q)=1$. Then, for any value
$c\in M_{p/q}$ we will assign a value $\la \in L_{q,0}$.          

We start  in the dynamical plane of $Q_c$. Through cutting and sewing, we construct a new map $f^{(1)}_c$, the {\em first return map}. This  map has several lines of discontinuity.

Next we smooth $f\one_c $ on sectors and obtain a new map $f\two_c$ which is quasi-regular and has two critical points: one with multiplicity $q-1$ and the other with multiplicity one.

To obtain an analytic map, we define an $f\two_c$-invariant almost complex structure, $\sigma$. By the Measurable Riemann Mapping Theorem,  this structure can be integrated by a quasi-conformal homeomorphism $\varphi_c$, inducing the analytic map $f^{(3)}_c=\varphi_c \circ f^{(2)}_c \circ \varphi_c^{-1}$.

We then show  that $f\three_c$  is a polynomial-like map (see sect.~\ref{tools} or \cite{dh3}) of degree $q+1$. By the Straightening Theorem, it is hybrid equivalent to a polynomial $P$ of degree $q+1$. Finally, we show that $P$ may be chosen to belong to the family $L_{q,0}$, and that it is independent of the choices made in the above steps.

This process defines the map $\ho: M_{p/q} \longrightarrow L_{q,0}$. In Sec. \ref{polike}, we explain the construction in the dynamical plane. In Sect. \ref{bijectivity}, \ref{analyticity} and \ref{continuity} we prove that \hho\  is  bijective, analytic in the interior of $M_{p/q}$ and continuous at points in the boundary of $M\p$, respectively.

%*****************************************************************************

\subsection{Construction of a polynomial of degree $q+1$. Definition of $\phi_{p/q}$.}\label{polike}

We start in the dynamical plane of a polynomial $Q_c$ with $c\in
M_{p/q}$, for a fixed $p/q$, and with $W_c$, $V_c^i$, $S_c(\theta)$, etc. as in section \ref{julia of Q}. Let $V_c^{0 \ldots q-1}=\bigcup_{i=0}^{q-1} V_c^i$. Our goal is to obtain a value $\la=\la(c) \in L_{q,0}$.

 We construct a truncated space $\bc_c^T$ obtained from $V_c^{0 \ldots q-1}$ by identifying 
\[
	R_c(\thetatil^{0}) \mbox{\ \ \ \ with \ \ \ \ } R_c(\thetatil^{q-1})
\] 
equipotentially (see Figs.~\ref{domainf13} and \ref{domainf25}), that is, we identify a point $z\in R_c(\thetatil^{0})$ with the unique point $w\in R_c(\thetatil^{q-1})$ such that $G_c(z)=G_c(w)$. Throughout the rest of the section we let $\thetatil$ denote this identified argument.

\begin{figure}[htbp]
%\vspace{2cm}
\centerline{\epsfysize=6cm\epsffile{pictures/domainf13.ai} \epsfysize=6cm\epsffile{pictures/domainf13rhp.ai}}
\caption{\small Left:The space $\bc_c^T$ where a random equipotential has been drawn and its preimage under $f_c\one$, for $c \in M_{1/3}$. Right: The space $\bh^T$ with a vertical line and its preimage under $f\one$, when working in the $1/3$-limb.}
\label{domainf13}
\end{figure}

\begin{figure}[htbp]
%\vspace{2cm}
\centerline{\epsfysize=6cm\epsffile{pictures/domainf25.ai} \epsfysize=6cm\epsffile{pictures/domainf25rhp.ai}}
\caption{\small Left:The space $\bc_c^T$ where  a random equipotential has been drawn and its preimage under $f_c\one$, for $c \in M_{2/5}$. Right: The space $\bh^T$ with a vertical line and its preimage under $f\one$, when working in the $2/5$-limb.}
\label{domainf25}
\end{figure}

Note that the space $\bc_c^T$ can be viewed as the quotient of $V_c^{0\ldots q-1}$  union some neighborhoods of the rays $R_c(\thetatil^{q-1})$ and $R_c(\thetatil^{0})$ by the equivalence relation identifying the points in the two neighborhoods 
\[ 
	\psi_c(\exp(\rho+2\pi i (\thetatil^{0}+t))) \sim \psi_c(\exp(\rho+2\pi i (\thetatil^{q-1}+t))). 
\] 
Thus, $\bc_c^T$ is a Riemann surface isomorphic to $\bc$. 

Define the {\em truncated filled Julia set} as 
\[
	K_c^T=K_c \cap \bc_c^T=K_c \cap V_c^{0 \ldots q-1}.
\]
Note that no identification takes place in $K_c$, so $K_c^T$ can be viewed as a subset of $K_c$.

We define now a map $f_c^{(1)}:\bc_c^T \longrightarrow \bc_c^T$, the {\em first return map}, so that $f_c\one$ acts on the sets $\smash[t]{\overset{\circ}{V}}_c^i$ similarly to the polynomial $P_{q,\la}$ with $\la \in L_{q,0}$ (compare with (\ref{action2}) on page \pageref{action2}). That is we define
\[
	f_c^{(1)}=
		\begin{cases}
			Q_c^q	  & \text{if $z\in \smash[t]{\overset{\circ}{V}}_c^0$}, \\
			Q_c^{q-1} & \text{if $z\in \smash[t]{\overset{\circ}{V}}_c^p$} \\
			Q_c^{q-2} & \text{if $z\in \smash[t]{\overset{\circ}{V}}_c^{2p}$}	\\
			\ldots	  & 					\\
			Q_c	  	  & \text{if $z\in \smash[t]{\overset{\circ}{V}}_c^{(q-1)p} =\smash[t]{\overset{\circ}{V}}_c^{q-p}$}	\\
			\alpha_c	  & \text{if $z=\alpha_c^\prime$}
		\end{cases}
\]
\begin{figure}[htbp]
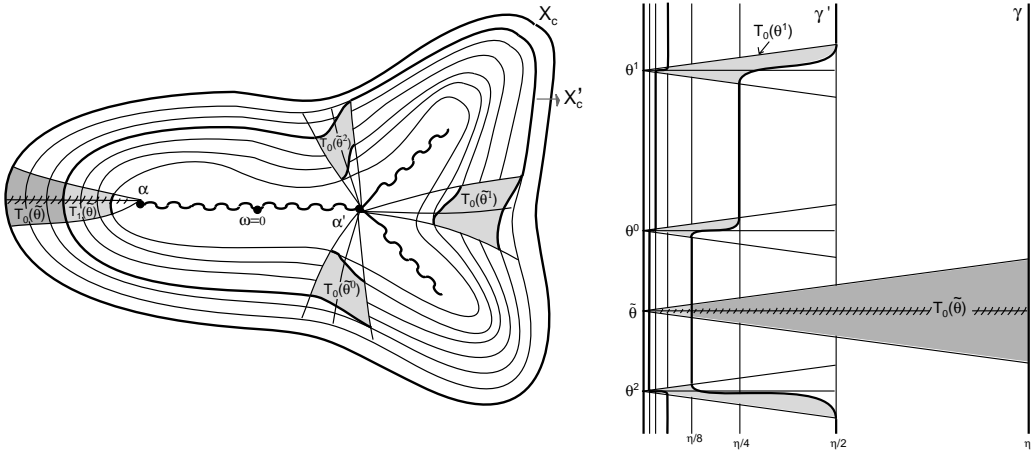

%\vspace{2cm}
\centerline{\epsfysize=6cm\epsffile{pictures/domain13dplane.ai} \epsfysize=6cm\epsffile{pictures/domain13.ai}}
\caption{\small Left:The space $\bc_c^T$ and the domain and the image of $f_c\two$, for $c\in   M_{1/3}$. Right: The space $\bh^T$ and the domain and the image of $f\two$ when working in the $1/3$-limb.}
\label{domain13}
\end{figure}

Note that $f_c\one$ is well defined on the identified ray $R_c(\thetatil)$ and holomorphic everywhere except at $\alpha_c^\prime$ and at 
\[
	R'_c=R_c(\theta^0)\cup R_c(\theta^1)\cup \cdots \cup R_c(\theta^{q-1})
\]
where it is discontinuous. 

Informally, just to get some intuition, note as well that the first return
map already has some of the characteristics of a map of the family \ppol. For example, each of the subsets $\smash[t]{\overset{\circ}{V}}_c^i$ is mapped under $f_c^{(1)}$ to the space $\bc_c^T\setminus (R_c(\thetatil) \cup \{ \alpha_c \} )$. The map has degree two in $\smash[t]{\overset{\circ}{V}}_c^0$, and degree one on the other subsets $\smash[t]{\overset{\circ}{V}}_c^i$, which justifies the hope that some variation of the first return map will be a map of degree $q+1$. Also, we can observe the creation of a new critical point, at $\alpha_c^\prime$. This is still a prefixed point, but a neighborhood around it now wraps $q$ times around $\alpha_c$ under $f_c^{(1)}$ , making $\alpha^\prime_c$ a critical point of multiplicity $q-1$. 

This ends the first step of the construction. We will now modify the map $f_c\one$  to construct a new map $f_c\two$  which will be quasi-regular. The
modification is done only on neighborhoods of the rays (the sectors defined in section \ref{the sectors}) for which $f_c\one$  is discontinuous. We denote these sectors by $S'_c$, i.e.,
\[
	S'_c=S_c(\theta^0)\cup S_c(\theta^1)\cup \cdots \cup S_c(\theta^{q-1}).
\]
We will have to restrict the space so that these sectors do not overlap. We will do this construction once and for all in the right half plane.

We define $\bh^T$ from $\bh$ by cutting along $R(n+\thetatil^{0})$ and $R(n+\thetatil^{q-1})$, $n\in \bz$ and identifying the two rays as we did for $\bc_c^T$. Then, the universal covering map $\psi_c \circ \exp :\bh \longrightarrow \bc \setminus K_c$ restricts to a universal covering map 
\[
	(\psi_c \circ \exp)^T :\bh^T \longrightarrow \bc_c^T \setminus K_c^T
\]
which is holomorphic in both variables.

Let $R', S', S$ and  $V^i$ denote the sets in $\bh^T$ that project to $R'_c, S'_c, S_c(\thetatil)$ and $V^i_c$ respectively.  

Choose a lift $f\one$ of $f_c\one$. Then, the following diagram commutes
\[
  \begin{CD}
	\bh^T   @>f\one>>  \bh^T \\
	@V{(\psi_c \circ \exp )^T}VV @VV{(\psi_c \circ \exp )^T}V   \\
    \bc_c^T\setminus K_c^T @>{f_c\one}>> \bc_c^T \setminus K_c^T 
  \end{CD}
\]
Note that $f\one$ can be chosen independently of $c$ and that one can actually write down such a map explicitly. The map $f\one$ is holomorphic everywhere except at $R(\theta^i+n)$, $0 \leq i \leq q-1$, $n \in \bz$ where it is discontinuous. 

We proceed now to restrict the domain. We will define two sets $X'$ and $X$ which will be the domain, respectively the image, of the new map $f\two$. Set $X=W \cap \bh^T$ and let $\gamma$ denote the curve in $\bh^T$ which bounds $X$. Let $\gamma'$ denote a ${\mathcal C}^\infty$ curve which projects to a Jordan curve $\gamma'_c=(\psi_c \circ \exp)^T(\gamma')$ in $\bc_c^T$ and which equals
\[ 
   \mbox{Re}(\gamma')= 
		\begin{cases} 
			\frac{\eta}{2^q} & \text{if $z\in V^0\setminus S'$} \\
			\frac{\eta}{2^{q-1}} & \text{if $z\in V^p \setminus S'$} \\
				\ldots		  & 					\\
			\frac{\eta}{2}  & \text{if $z\in V^{(q-1)p}\setminus S' = V^{q-p} \setminus S'$} 		
		\end{cases}
\]
Let $X'$ be the vertical strip bounded by $\gamma'$ (see Fig.~\ref{domain13}). In dynamical plane, set $X_c=W_c \cap \bc_c^T$ and let $X'_c$ denote the filled set in $\bc_c^T$ bounded by $\gamma'_c$.

We shall modify the map $f\one$ on the sectors around the rays of discontinuity as shown in Fig. \ref{domain13} so that it induces a quasi-regular map
\[
	f_c\two : X'_c \longrightarrow X_c
\]
First we divide  the sectors in $S'$ into quadrilaterals as described below.

For $S(\theta) \subset (S' \cup S)$ we have a homeomorphism 
\[
	\Pi_\theta : S(\theta) \longrightarrow S(\theta)
\]
which is multiplication by $2^q$ followed by a vertical translation that depends on $\theta$. 

For $S(\theta) \subset S$ we define
\begin{align}
	T_0(\theta)&=\overline{(S(\theta)  \cap X) \setminus \Pi_\theta^{-1}(S(\theta) \cap X)} \\
\intertext{and for $k>0$}
	T_k(\theta)&=\Pi_\theta^{-1}(T_{k-1}(\theta)).
\end{align}
For $S(\theta) \subset S'$, define
\[
	T_0(\theta)=\overline{(S(\theta)\cap X') \setminus \Pi_\theta^{-1}((S(\theta) \cap X'))}
\] and $T_k(\theta)$ as above.

Note that the map $\Pi_\theta$ restricted to $T_k(\theta)$ is a conformal isomorphism from $T_k(\theta)$ to $T_{k-1}(\theta)$.

We define $f\two:X' \longrightarrow X$ to equal $f\one$ outside of $S'$. On these sectors, $f\two$ is defined by induction as follows. Let $n_{i}$ denote the integer determined by
\[
	f\one: \partial(S(\theta^{i})) \longrightarrow \partial(S(n_{i}+\thetatil)).
\]
Then choose for $0 \leq i \leq q-1$ a diffeomorphism
\[
	f_{0,i}\two :T_0(\theta^{i}) \longrightarrow T_0(n_{i}+\thetatil)
\]
so that $f\one$ and $f\two_{0,i}$ determine the same tangent map on the boundary of the sectors. For $k>0$ define
\[
	f_{k,i}\two: T_k(\theta^{i}) \longrightarrow T_k(n_{i} +\thetatil)
\]
by induction as the map satisfying
\[
	\Pi_{n_{i}+\thetatil} \circ f\two_{k,i}=f\two_{k-1,i} \circ \Pi_{\theta^{i}}.
\]

Finally extend $f\two:X' \longrightarrow X$ as a covering transformation so it is compatible with the projection. The map $f_c\two:X'_c \longrightarrow X_c$ is defined on $\bc_c^T \setminus K_c^T$ by
\[
  \begin{CD}
  X' @>f\two>> X\\
  @V(\psi_c \circ \exp )^TVV @VV(\psi_c \circ \exp )^TV\\
  X'_c\setminus K_c^T @>f_c\two>> X_c\setminus K_c^T
  \end{CD}
\]
and $f_c\two|_{K_c^T}=f_c\one|_{K_c^T}$.

The map $f\two: X' \longrightarrow X$ is a diffeomorphism. The map $f\two_c: X'_c \longrightarrow X_c $ is a degree $(q+1)$-ramified covering.

\begin{lemma} \label{quasiregular} The  map $f\two$   is quasi-conformal. Hence $f_c\two$ is quasi-regular.
\end{lemma}

\begin{proof}{} We must show that the field of ellipses $E_x=(T_x f^{(2)})^{-1}(S^1)$ for $x\in X'$, has bounded dilatation ratio.

If $x \notin S^\prime$, then $E_x$ is a circle, since the map $f\two$ outside $S^\prime$ equals $f\one$  which is analytic.

If $x \in T_0(\theta^{i})$ for some $i$,  then the ratio of the axes of $E_x$ is bounded by some constant $K$, since $f^{(2)}_{0,i}$ is a diffeomorphism and $T_0(\theta^{i})$ is compact.

If $x \in T_k(\theta^{i})$ for some $k>0$ then the ratio of the axes of $E_x$ is bounded by the same constant $K$, since $f^{(2)}_{k,i}$ is
$f^{(2)}_{0,i}$ composed with a finite number of  analytic maps.

The same can be said for points that belong to the integer translations of this sectors since $f\two$ is a covering transformation.

Hence, the dilatation ratio of the field of ellipses $(E_x)_{x \in X'}$ is bounded by the same constant for all $x\in X'$, and $f\two$  is quasi-conformal. Since $(\psi_c \circ \exp)^T$ is conformal, $f_c\two$ is quasi-regular.  
\end{proof}

\begin{remarks}{} \label{shishi}
\begin{enumerate}
\item It is  important to note that orbits enter $S^\prime$ at most once. This allows us to bound the dilatation ratio by the same constant everywhere, since the diffeomorphism $f^{(2)}_{0,i}$ has to be applied only once. Moreover it is essential for the next step in the construction when we change the complex structure. This is the Shishikura principle for surgery of this type.

\item The map $f_c\two$  is a degree $(q+1)$ ramified covering. It has two critical points: one at $\omega=0$ with multiplicity one and another one at $\alpha^\prime_c$ with multiplicity $q-1$, mapped to the fixed point $\alpha_c$. These are the topological characteristics of maps in the family $P_{q,\lambda}$, but the map $f_c\two$  is not holomorphic.
\end{enumerate}
\end{remarks}

The next step (holomorphic smoothing) is to construct a map $f_c\three$  which has the same properties as $f_c\two$  but so that it also is a polynomial-like map of degree $q+1$. In order to do so, we change the complex structure into a new almost complex structure with the property that it is preserved  under the map $f_c\two$ .

\begin{lemma} \label{newstructure} There  exists an
almost complex structure $\sigma_ c$ on $X_c$  quasi-conformally equivalent to the standard complex structure $\sigma_0$, such that $(f_c^{(2)})^\ast \sigma_c=\sigma_c$ and $\sigma_c=\sigma_0$ on the set $K_c^T$.
\end{lemma}

This lemma could also be phrased as follows: There exists a measurable field of ellipses in $X_c$ with bounded  ratio of the axes. These ellipses are invariant under  $f_c\two$ and the ellipses are circles on  $K_c^T$. As before, we will construct this almost complex structure in the right half plane.

\begin{proof}{} 
For $x\in S^\prime$, we define the complex structure $\sigma$ to be the same as in lemma \ref{quasiregular}, i.e. 
\[ 
	E_x= (T_x f^{(2)})^{-1} (S^1),
\] 
and then, the ratio of the axes is bounded by some constant $K$. 

If $x \notin S^\prime$, then either there exists $n>0$ such that $(f^{(2)})^n(x) \in S^\prime$ or the orbit of $x$ never enters $S^\prime$.

If there exists such $n$, then it is unique since a point in $S'$ is mapped to a point in $S$ and does not leave $S$ as long as the map is defined. In this case define 
\[ 
	E_x= (T_x (f^{(2)})^n)^{-1} (E_{(f^{(2)})^n(x)}). 
\] 
Since $f\two$  is holomorphic outside of the sectors in $S^\prime$, we have that $(f^{(2)})^n$ is holomorphic in a neighborhood of $x$. Hence, the dilatation ratio of $E_x$ is the same as the dilatation ratio of $E_{(f^{(2)})^n(x)}$, and therefore bounded by $K$. 

Finally, if the orbit of $x$ never enters $S^\prime$, define $E_x=S^1$.

By construction, $\sigma$ is invariant under $f\two$. Define $\sigma_c$ on $X_c\setminus K_c^T$ as the pull-back of $\sigma$ by the map $((\psi_c \circ \exp )^T)^{-1}$ on a fundamental domain, and set $\sigma_c=\sigma_0$ on $K_c^T$. By construction, $\sigma_c$ is the required almost complex structure.  
\end{proof}

By doing this construction once and for all in the right half plane and pulling back by $(\psi_c \circ \exp)^T$ we have that $\sigma_c$ and $f_c\two$ vary holomorphically with respect to $c$.

Let $\sigma_c$ be the almost complex structure given by Lemma \ref{newstructure}. We now apply the Measurable Riemann Mapping Theorem (see Sect.~\ref{tools}) on $\stackrel{\circ}{X}_c$,  a Riemann surface isomorphic to $\bd$, to obtain a unique quasi-conformal homeomorphism $\varphi_c: \stackrel{\circ}{X}_c \longrightarrow \bd$ integrating $\sigma_c$ and satysfying $\varphi_c(\omega_c)=0$ and $\varphi_c(\alpha_c) \in \br_+$. We have
\[ 
	(T_x \varphi_c)^{-1}(S^1)=\rho_c(x) E_x  
\]
where $\rho_c(x)\in \br_+$.

Now set $ D_c^\prime= \varphi_c(\smash[t]{\overset{\circ}{X}}^\prime_c) \subset \bd $ and define
$ f_c^{(3)}$ by
\[
  \begin{CD}
	  \smash[t]{\overset{\circ}{X}}^\prime_c  @>{f^{(2)}_c}>> \stackrel{\circ}{X}_c \\
	  @V{\varphi_c }VV @VV{\varphi_c}V   \\
      D'_c @>{f_c\three}>> \bd 
  \end{CD}
\]

\begin{remark}\label{hybrid} In  this particular case, the quasi-conformal homeomorphism $\varphi_c$ is actually a hybrid equivalence, since $\sigma_c=\sigma_0$ on $K_c^T$.
\end{remark}
 
We show that $f_c\three$  is holomorphic as follows: $\varphi_c^{-1}$ takes the standard structure $\sigma_0$ of $U_c^\prime$ to the almost complex structure $\sigma_c$. Then, the map $f_c\two$  preserves $\sigma_c$ which is taken back to $\sigma_0$ by $\varphi_c$. We conclude that $f_c\three$  takes $\sigma_0$ to $\sigma_0$ and therefore  is holomorphic.
 
Since $\varphi_c$ is a homeomorphism, $f_c\three$  is still a ramified covering of degree $q+1$ with two critical points (see remarks \ref{shishi}). Since $D_c^\prime$ is relatively compact in $\bd$, we can conclude that $f_c\three$  is a polynomial-like map of degree $q+1$.

We now apply the Straightening Theorem to obtain a polynomial $P_{q,\la(c)}$ of degree $q+1$ in $L_{q,0}$, and a hybrid equivalence $\chi_c$ that conjugates $f_c\three$  to $P_{q,\la(c)}$ on neighborhoods of $K_c^T$ and $K_{\la(c)}$.

To justify this claim notice that the polynomial $\chi_c \circ f_c\three \circ \chi_c^{-1}$ satisfies the following properties:

\begin{itemize}
\item the point $z_1=\chi_c(\varphi_c(\alpha_c))$ is a repelling fixed point,
\item the point $z_2=\chi_c(\varphi_c(0))$ is a critical point of multiplicity one,
\item the point $z_3=\chi_c(\varphi_c(\alpha_c^\prime))$ is a critical point of multiplicity $q-1$, 
\item $z_2$ is mapped to $z_1$, and
\item $z_1$ is the landing point of one fixed ray.
\end{itemize} 
Since we may conjugate by an affine transformation sending $z_1$ to $0$ and $z_3$ to $-q$, we may assume that the polynomial is of the form $P_{q,\la(c)}$ . 

This concludes the construction on the dynamical plane. By this construction, given a parameter value $c\in M_{p/q}$ we have obtained a unique $\lambda(c)\in L_{q,0}$. We define the  map
\[ 
\ho:   M_{p/q}   \longrightarrow   L_{q,0}          
\] by $\ho(c)=\lambda(c)$.

\begin{remarks} \label{welldefined}
\begin{enumerate}
	\item The value of $\lambda$ does not depend on the choices made throughout the construction, i.e. the choice of the boundaries of $X$ and $X'$, the slope $s$, the diffeomorphisms $f_{0,i}^{(2)}$ and the integrating map $\varphi_c$. Indeed, suppose that by other choices $\widetilde{\eta}$, $\widetilde{\gamma'}$, $\widetilde{s}$, $\widetilde{f}^{(2)}_{0,i}$ and $\widetilde{\varphi}_c$ we obtain a map $\widetilde{f}^{(2)}$ and a polynomial $P_{q,\widetilde{\la}}$. Then, $f_c^{(2)}$ and $\widetilde{f}_c^{(2)}$ are hybrid equivalent and also \ppol\  and
$P_{q,\widetilde{\la}}$ would be hybrid equivalent (because of remark \ref{hybrid}). Hence  by Prop.~\ref{choices} we have $\la=\widetilde{\la}$.

	\item Following the construction of the Straightening map in \cite{dh3}, one can check that   $\chi_c$ can be defined in all of $\varphi_c(X_c)$, since the boundary of $X_c$ is smooth and $\varphi_c$ is quasi-conformal.
\end{enumerate}
\end{remarks} 

The following proposition addresses the question of uniqueness of the conjugating maps.
\begin{proposition}
Suppose $c \in M\p$ and $\la=\ho(c) \in L\qo$. A hybrid equivalence between $f_c\three$ and $P_\la$ is uniquely determined on $\varphi_c(K_c^T)$, hence $\chi_c \circ \varphi_c$ is uniqueley determined on $K_c^T$.
\end{proposition}
\begin{proof}{}
Suppose $\chi_c$ and $\chi'_c$ are hybrid equivalences between $f_c\three$ and $P_\la$. Then $\varphi=\chi'_c \circ \chi_c^{-1}$ is a hybrid equivalence between $P_\la$ and itself. By copying the proof of Prop.~6 in \cite{dh3}, p.~302 defining
\[
	\Phi=
		\begin{cases}
			\varphi & \text{on } K_\la \\
			\text{Id} & \text{on } \bc \setminus K_\la
		\end{cases}
\]
we obtain that $\Phi$ is holomorphic and hence the identity on all of $\bc$.
\end{proof}

The following theorem is the dynamical counterpart of Theorem F.
 
\begin{theorem*}{G}
Suppose $c \in M\p$ and $\la=\ho(c) \in L\qo$. Recall that $X_c$ is a neighborhood of $K_c^T$ bounded by the truncated equipotential $\gamma_c$ and that $W_\la$ is a neighborhood of $K_\la$ namely the filled level set of a chosen equipotential of $G_\la$. There exists a homeomorphism
\[ H_c: X_c \longrightarrow W_\la \]
which satisfies $\overline{\partial} H_c=0$ on $K_c^T$ and conjugates $f_c\two$ to $P_\la$ and such that $H_c$ maps relevant sectors and rays for $Q_c$ to their counterparts for $P_\la$, i.e. $H_c$ maps $S_c(\thetatil)$, $R_c(\thetatil)$, and $R_c(\theta^i)$ to $S_\la(0)$, $R_\la(0)$, and $R_\la(\frac{i+1}{q+1})$ respectively.

Moreover, $H_c|_{K_c^T}$ is uniquely determined and it conjugates the first return map $f_c\one$ to $P_\la$. \newline
\end{theorem*}

\begin{proof}{}
The map
\[
	  \chi_c \circ \varphi_c :X_c \rightarrow \chi_c ( \varphi_c(X_c))
\]
is a quasi-conformal homeomorphism from $X_c$ to a neighborhood of $K_\la$ conjugating by construction the map $f_c\two$ to the polynomial  $P_\la$. In particular, $H_c$ conjugates the first return map $f_c\one$ on $K_c^T$ to the polynomial $P_\la$ on $K_\la$. By remark \ref{hybrid}, it satisfies $\overline{\partial} H_c=0$ on $K_c^T$. 

The map $\chi_c \circ \varphi_c$ sends $S_c(\thetatil)$, $R_c(\thetatil)$, $R_c(\theta^i)$ and $\gamma_c$ to some quasi-conformal image of these objects. We claim that there exists a hybrid equivalence $h_\la: \chi_c ( \varphi_c(X_c)) \longrightarrow W_\la$ conjugating $P_\la$ to itself and such that $H_c:=h_\la \circ \chi_c \circ \varphi_c$ satisfies the properties in the theorem. In the remainder of this proof we give the idea of the construction of $h_\la$.

Let $\frak{A}$ and $A$ denote the annuli 
\begin{align*}
	\frak{A} &=\overline{\chi_c (\varphi_c (X_c)) \setminus P_\la^{-1}(\chi_c (\varphi_c (X_c)))} \\
	A &= \overline{W_\la \setminus W_{\la,1}}.
\end{align*}
First define $h_\la :\frak{A} \longrightarrow A$ as a diffeomorphism conjugating $P_\la$ to itself on the inner boundaries and sending  $\chi_c (\varphi_c (S_c(\thetatil))) \cap \frak{A}$, $\chi_c (\varphi_c (R_c(\thetatil)))\cap \frak{A}$ and $\chi_c (\varphi_c (R_c(\theta^i))) \cap \frak{A}$ to $S_\la(0) \cap A$, $R_\la(0) \cap A$ and  $R_\la(\frac{i+1}{q+1}) \cap A$ respectively (see Fig.~\ref{annuli}). This can easily be done constructing first such a map in the right half plane.

\begin{figure}[htbp]
%\vspace{2cm}
\centerline{\epsfysize=5cm\epsffile{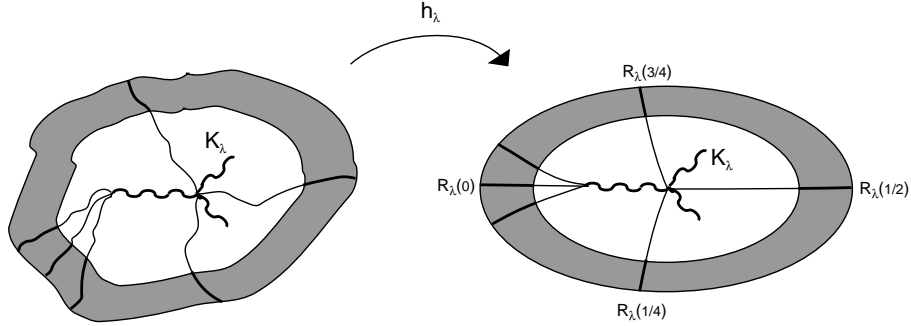}}
\caption{\small The annuli $\frak{A}$ and $A$ and the function $h_\la$. }
\label{annuli}
\end{figure}

Once $h_\la$ is correctly defined on the annulus we use a pull-back argument to define it everywhere else outside $K_\la$, that is if we denote $\frak{A}^n=P_\la^{-n}(\frak{A})$ and $A^n=P_\la^{-n}(A)$,  we define $h_\la=P_\la^{-n} \circ h_\la \circ P_\la^{n}$ on $\frak{A}^n$. Finally, define $h_\la$ to be the identity on $K_\la$. By an argument analogous to Lemma 1 in page 302 of \cite{dh3}, followed by Rickman's lemma in page 303, one can show that $h_\la$ is a hybrid equivalence between $P_\la$ and itself with the required properties.

That $H_c|_{K_c^T}$ is uniquely defined follows from the proposition above.
\end{proof}

\begin{remarks} \label{conju}
\begin{enumerate}
	\item If $c$ is Misiurewicz, i.e. if $\omega_c=0$ is strictly
preperiodic for $Q_c$, then $\omega_\la=\frac{-q}{q+1}$ is strictly preperiodic for $P_\la$. In particular if $\omega_c$ is eventually mapped to $\alpha_c$ then $\omega_\la$ is eventually mapped to the fixed point 0.

	\item If $c$ is hyperbolic, i.e. if $\omega_c$ is attracted to an attracting cycle by $Q_c$, then $\omega_\la$ is attracted to an attracting cycle by $P_\la$. Moreover, the multiplier of the cycle by $Q_c$ equals the multiplier of the cycle by $P_\la$, since $H_c$ is holomorphic in the interior of $K_c^T$.

	\item The main branch point (Misiurewicz) in the $p/q$-limb of $M$ corresponds to $Q_c^q(\omega_c)=\alpha'_c$. It follows that the main branch point  in $L\qo$ corresponds to the polynomial $\pol(\omega_\la)=-q$ or $\la=(q+1)(1+\frac{1}{q})^q$. Hence this branch point tends to $\infty$ as $q$ tends to $\infty$.
\end{enumerate}
\end{remarks}

\begin{figure}[htbp]
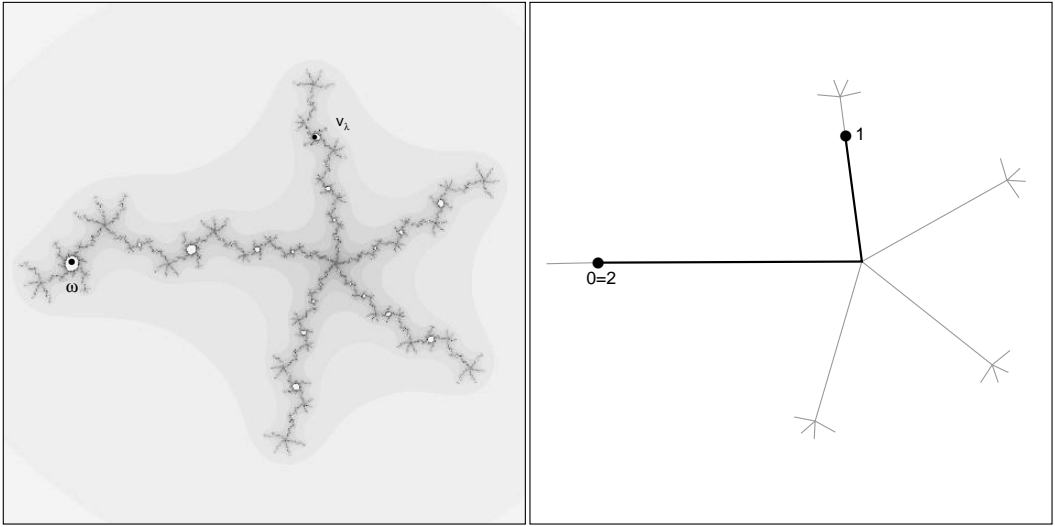

%\vspace{2cm}
\centerline{\epsfysize=7cm\epsffile{pictures/julpol5.ai} \epsfysize=7cm\epsffile{pictures/treejulpol5.ai}}
\caption{\small Left: The filled Julia set for $P_{5,\la}$ where $\la$ is the center of a hyperbolic component of period two in $L_{5,0}$; $\la$ is the image of both $c$-values in Fig.~\ref{juliafifths} under $\phi_{1/5}$ and $\phi_{2/5}$ respectively. Right: The Hubbard tree (compare with Fig.~\ref{treesfifths}) }
\label{julpol5}
\end{figure}

%*****************************************************************************
\subsection{Bijectivity of \hho} \label{bijectivity}

In this section, we construct a map
\[ 
	\hoi : L_{q,0} \longrightarrow M_{p/q}.
\] 
and show that it is the inverse of $\phi_{p/q}$.

We start with a fixed $\la \in L_{q,0}$ and our goal is to obtain a unique quadratic polynomial
$Q_c$ such that $c \in M_{p/q}$. Note that this means that,  from each $L_{q,0}$, we have to construct several maps, one for each $p$.

Let $V_\la^i$, etc. be defined as in Sect. \ref{polyn}. Since this construction is supposed to be the opposite of the construction of Sect. \ref{polike}, we first want to ``add'' the piece of dynamical plane that we removed before.

Hence, we cut $V_\la^0$ along $R_\la(0)$ in such a way that now $V_\la^0$ contains two copies of $R_\la(0)$, which we denote by $R_\la(0)$ and $R_\la(1)$. For $1 \leq i \leq q-1$, let $\Vtil_\la^i$ be an identical copy of $V_\la^i$ and let
\[ 
	\tau_\la^i : V_\la^i \longrightarrow \Vtil_\la^i 
\] 
be the identity map.  The copies $\Vtil_\la^i$ are sewn together in order to resemble Figs.~\ref{center13} and \ref{center25}, as shown in figures \ref{Pdomain13} or \ref{Pdomain25}. We glue the boundaries (the rays) together in order to construct a Riemann surface $\bc_\la^E$ isomorphic to $\bc$. However, we will not identify the rays equipotentially as we did before. (If we did, the map that we are going to  define  would not be continuous on these rays). In what follows, we describe these identifications. Figs. \ref{Pdomain13} and \ref{Pdomain25} show two examples.

\begin{figure}[htbp]
%\vspace{2cm}
\centerline{\epsfysize=7cm\epsffile{pictures/Pdomain13.ai}}
\caption{\small The space $\bc_\la^E$ and the domain of $f_\la\two$ for $\la \in L_{3,0}$ and $p=1$. }
\label{Pdomain13}
\end{figure}

\begin{figure}[htbp]
%\vspace{2cm}
\centerline{\epsfysize=7cm\epsffile{pictures/Pdomain25.ai}}
\caption{\small The space $\bc_\la^E$ and the domain of $f_\la\two$ for $\la \in L_{5,0}$ and $p=2$. }
\label{Pdomain25}
\end{figure}

Let $N$ be the bijective function that assigns to each integer $1 \leq i \leq q-1$ the unique integer $1 \leq N(i) \leq q-1$ such that
\[ 
	N(i) \; p \equiv i \pmod{q}. 
\] 
Intuitively, if this were the dynamical plane of a quadratic polynomial $Q_c$ (with $c\in M_{p/q}$), $N(i)$ would be  the number of iterates that it takes for $V_\la^0$ to be mapped to $\Vtil_\la^i$ for the first time.  

Pick $1 < i < q-1$ and suppose $N(i-1) < N(i)$. Then, for each point $x \in R_\la(i/(q+1))$ of potential $\rho$ (i.e. $G_\la(x)=\rho$), we identify
\[ 
	\tau_\la^{i-1} (x) \sim \tau_\la^i (y), 
\] 
where $y \in R_\la(i/(q+1))$ and has potential $G_\la (y)= \rho \, / \, (q+1)^{\frac{N(i)-N(i-1)}{q}}$.

Otherwise, if $N(i-1) > N(i)$, then we identify
\[ 
	\tau_\la^{i} (x) \sim \tau_\la^{i-1} (y), 
\] 
where $x,y \in R_\la(i/(q+1))$, $G_\la(x)=\rho$, and $ G_\la (y)= \rho  \, / \, (q+1)^{\frac{N(i-1)-N(i)}{q}}$. 

If  $i=1$, for each $x \in R_\la(1)$ with potential $\rho$ we identify
\[ 
	x \sim \tau_\la^1(y), 
\]  
where $y \in R_\la (1/(q+1))$ and has potential $G_\la(y)= \rho \, / \, (q+1)^{\frac{N(1)}{q}}$.

If $i=q-1$ then for each $x\in R_\la(0)$, we identify
\[ 
	x \sim \tau_\la^{q-1} (y), 
\] 
where $y \in R_\la(q/(q+1))$ and has potential $G_\la(y) = \rho (q+1)^{\frac{N(q-1)}{q}}$. Finally, we identify the landing points of the rays bounding the $\Vtil_\la^i$ with $0$. We shall denote these $q$ rays by $\rtil_\la(\frac{i}{q+1})$, for $1\leq i \leq q$.

We define 
\[ 
	\bc_\la^E=\bigcup_{i=0}^{q-1} V_\la^i \quad \cup \quad \bigcup_{i=1}^{q-1} \Vtil_\la^i / \sim
\]
with the identifications above.

Define the extension of $K_\la$ to be $K_\la^E=K_\la \cup \bigcup_{i=1}^{q-1}  \tau_\la^i(K_\la) / \sim$. Note that $K_\la$ can be viewed as a subset of $K_\la^E$.

We define
\[
	f_\la\one :\bc_\la^E \longrightarrow \bc_\la^E
\]
so that $f_\la\one$ acts on the sets $V_\la^0, V_\la^i, \Vtil_\la^i$ similarly to the quadratic polynomial $Q_c$ with $c \in M_{p/q}$ (compare with (\ref{action}) on page \pageref{action}). It should be viewed as the ``inverse'' of the first return map. 

Recall that $P_\la: \smash[t]{\overset{\circ}{V}}^i_\la \longrightarrow \bc \setminus (R_\la(0) \cup \{ 0\} )$ is a homeomorphism and define then 

\[
 f_\la^{(1)}= 
	\begin{cases} 
        \tau_\la^{(i+1)p} \circ P_\la|_{V_\la^{(i+1)p}}^{-1} \circ P_\la	  &	\text{on $V_\la^{ip}$, \ \ for $0\leq i \leq q-2$} 	\\
      P_\la						 &	\text{on $V_\la^{q-p}$}	\\
      f_\la^{(1)} \circ (\tau_\la^{ip})^{-1} & 	\text{on $\Vtil^{ip}_\la$,
\ \ for $1\leq i \leq q-1$}
   \end{cases}
\]

Note that $f_\la\one$  is well defined on the identified rays 
\[
	\rtil_\la= \rtil_\la(\frac{1}{q+1}) \cup  \rtil_\la(\frac{2}{q+1}) \cup \cdots \cup \rtil_\la(\frac{q}{q+1}) 
\]

and holomorphic everywhere except at $-q$ and at  
\[ 
	R_\la= R_\la(\frac{1}{q+1}) \cup  R_\la(\frac{2}{q+1}) \cup \cdots \cup R_\la(\frac{q}{q+1}) 
\]
where it is discontinuous. As before, we shall smoothen the map on sectors and restrict the space so that the sectors do not overlap. We will do this - once and for all - in the right half plane.

We define $\bh^E$ from $\bh$ by cutting along $R(n)$, for $n \in \bz$, and gluing in copies of $V^i$ with similar identifications along the rays $\rtil(n+\frac{i}{q+1}), 1\leq i \leq q, n\in \bz$, as in $\bc_\la^E$. The universal covering map
\[
	\psi_\la \circ \exp  : \bh \longrightarrow \bc \setminus K_\la
\]
extends in an obvious way to a universal covering map
\[
    (\psi_\la \circ \exp )^E : \bh^E \longrightarrow \bc_\la^E \setminus K_\la^E.
\]
The map $(\psi_\la \circ \exp )^E$ is holomorphic in both variables. Let $S$, $\widetilde{S}$, $\smash[t]{\Vtil}^i$ denote the sets in $\bh^E$ that project to $\underset{i}{\bigcup} S_\la(\frac{i}{q+1})$, $\underset{i}{\bigcup} \smash[t]{\widetilde{S}}_\la(\frac{i}{q+1})$ and $\smash[t]{\Vtil}^i$ respectively.

Choose a lift $f\one$ of $f_\la\one$. Then the following diagram commutes:
\[
  \begin{CD}
	\bh^E   @>{f\one}>> \bh^E \\
	@V{(\psi_\la \circ \exp )^E}VV   @VV{(\psi_\la \circ \exp )^E}V   \\
     \bc_\la^E\setminus K_\la^E @>{f_\la\one}>> \bc_\la^E\setminus K_\la^E 
  \end{CD}
\]
Note that we can choose $f\one$ independently of $\la$. The map $f\one$ is holomorphic, except at $R(n+\frac{i}{q+1}), 1\leq i \leq q-1, n\in \bz$ where it is discontinuous.

Let $X \subset \bh^E$ denote the vertical strip in $\bh^E$ bounded by the curve $\gamma$ which equals
\[ 
	\mbox{Re}(\gamma)= 
		\begin{cases}
			\eta 	& \text{in $V^{0 \ldots q-1}$}		\\
			\frac{\eta}{(q+1)^{i/q}} & \text{in $\Vtil^{ip}$}
		\end{cases}
\]
Due to the identification of the rays, $\gamma$ is a vertical line. Set $\gamma_\la=(\psi_\la \circ \exp )^E(\gamma)$. Let $\gamma'$ denote a ${\mathcal C}^\infty$ curve which projects to a Jordan curve
\[
	\gamma'_\la=(\psi_\la \circ \exp )^E(\gamma')
\]
in $\bc_\la^E$, and which equals
\[ 
	\mbox{Re}(\gamma')= 
		\begin{cases}
			\frac{\eta}{(q+1)^{(i+1)/q}} 	& \text{in $V^{ip}\setminus S', 0\leq i \leq q-1$}		\\
			\frac{\eta}{(q+1)^{(i+1)/q}} & \text{in $\Vtil^{ip}, 1\leq i \leq q-1$}
		\end{cases}
\]
Let $X'$ be the vertical strip in $\bh^E$ bounded by $\gamma'$ and let $X_\la$, respectively $X'_\la$ denote the filled sets in $\bc_\la^E$ of $\gamma_\la$, respectively $\gamma'_\la$.

We shall modify the map $f\one$ on the sectors around the rays of discontinuity so it induces a quasi-regular map
\[
	f_\la\two : X'_\la \longrightarrow X_\la.
\]

For $\theta=n +\frac{i}{q+1}$, $n\in \bz$, $0\leq i \leq q$ we have a homeomorphism
\[ 
	\Pi_\theta : S(\theta) \longrightarrow S(\theta)
\]
which is multiplication by $q+1$ composed with a vertical translation.

For the same values of $\theta$, the maps $\Pi_\theta$ induce homeomorphisms 
\[
	\pitil_\theta: \estil(\theta) \longrightarrow \estil(\theta).
\]
Define
\[
	\Ttil_0(\theta)=\overline{(\estil(\theta) \cap X) \setminus \pitil_\theta^{-1}(\estil(\theta) \cap X)},
\]
and for $k>0$
\[
	\Ttil_k(\theta)=\pitil_\theta^{-1} (\Ttil_{k-1}(\theta)).
\]
For $\theta=n +\frac{i}{q+1}$, $n\in \bz$, $1\leq i \leq q$, define  
\[
	T_0(\theta)=\overline{(S(\theta) \cap X') \setminus \Pi_\theta^{-1}(S(\theta) \cap X')},
\]
and $T_k(\theta)$ as above.

We define $f\two:X' \longrightarrow X$ to equal $f\one$ outside of $S$.  On these sectors, $f\two$ is defined by induction. Let $n_i$ denote the integer determined by
\[
	f\one: \partial (S(\frac{i}{q+1})) \longrightarrow \partial(\estil(n_i+\frac{i+p}{q+1})).
\] 
Choose for $1 \leq i \leq q$  a diffeomorphism
\[
	f\two_{0,i}: T_0(\frac{i}{q+1}) \longrightarrow \Ttil_0(n_i+\frac{i+p}{q+1})
\]
so that $f\one$ and $f\two_{0,i}$ determine the same tangent map on the boundary of the sectors.

Define $f\two_{k,i}:T_k(\frac{i}{q+1}) \longrightarrow \Ttil_k(n_i+\frac{i+p}{q+1})$ by induction as the map satisfying 
\[
	\pitil_{n_i+\frac{i+p}{q+1}} \circ f\two_{k,i}=f\two_{k-1,i} \circ \Pi_{\frac{i}{q+1}}.
\]
Finally, extend $f\two:X' \longrightarrow X$ as a covering transformation so that it is compatible with the projection. The map $f_\la\two:X'_\la \longrightarrow X_\la$ is defined on $X'_\la\setminus K_\la^E$ by
\[
  \begin{CD}
	  X'   @>{f\two}>> X \\
	  @V{(\psi_\la \circ \exp )^E}VV @VV{(\psi_\la \circ \exp )^E}V   \\
      X'_\la \setminus K_\la^E @>{f_\la\two}>> X_\la \setminus K_\la^E 
  \end{CD}
\]
and $f_\la\two|_{K_\la^E}=f_\la\one|_{K_\la^E}$.

The map $f\two:X' \longrightarrow X$ is a diffeomorphism. The map $f_\la\two:X'_\la \longrightarrow X_\la$ is a degree 2-ramified covering.	
\begin{lemma} 
	The map $f\two$  is quasi-conformal. Hence $f_\la\two$ is quasi-regular.
\end{lemma} 
\begin{proof}{}
As in Lemma \ref{quasiregular}, the map $f\two:X' \longrightarrow X$ is quasi-conformal, since $f\two_{0,i}$ are diffeomorphisms on compact sets and $f_{k,i}\two$ is obtained from $f_{0,i}\two$ by composition with holomorphic maps.  
\end{proof}

We change the complex structure (as we did in Sect. \ref{polike}) into an almost complex structure $\sigma$, which is  quasi-conformally equivalent to the standard structure, $\sigma_0$, and such that it is invariant under the map $f\two$. Also, $\sigma$ coincides with $\sigma_0$ outside of $S$ and its preimages. Here, as in the quadratic case, orbits pass through these sectors at most once, since after one iteration,  a point in $S$ is mapped to a point in $\estil$ and do not leave $\estil$ as long as the map is defined.

By construction $\sigma$ is invariant under $f\two$. Define $\sigma_\la$ on $X_\la \setminus K_\la^E$ as the pullback of $\sigma$ by the map $((\psi_\la \circ \exp)^E)^{-1}$ and set $\sigma_\la=\sigma_0$ on $K_\la^E$.  (Compare with Lemma \fullref{newstructure}).

We apply the Measurable Riemann Mapping Theorem  to obtain the unique  quasi-conformal homeomorphism (in this case also a hybrid equivalence (see remark \ref{hybrid}))
\[ 
	\varphi_\la: \stackrel{\circ}{X}_\la \longrightarrow \bd, 
\] 
integrating $\sigma_\la$ and satisfying $\varphi_\la(\omega_\la)=0$ and $\varphi_\la(0) \in \br_+$. Then, let $D_\la^\prime= \varphi_\la (\smash[t]{\overset{\circ}{X}}^\prime_\la)$ and define
\[ f_\la^{(3)}= \varphi_\la \circ f_\la^{(2)} \circ \varphi_\la^{-1} :D_\la^\prime
\longrightarrow \bd. \]

The map $f_\la\three$  is a polynomial-like map of degree two, which has  one critical point at $\varphi_\la(-q/(q+1))$. We now apply the Straightening Theorem to obtain a uniquely determined polynomial $Q_{c(\la)}$ and a hybrid equivalence $\chi_\la$ conjugating $f_\la\three$ to $Q_{c(\la)}$. Hence we define
\[ 
	\hoi (\la)=c(\la) .
\]
\begin{remarks}
By the same arguments as in  remarks \fullref{welldefined}  the value of $c$ does not depend on the choices made throughout the construction, i.e.~the choice of the boundaries $\gamma$ and $\gamma'$, the slope $s$ and the diffeomorphisms $f_{0,i}\two$. We may also assume $\chi_\la$ to be defined on all $\varphi_\la(X_\la)$.
\end{remarks}
The following proposition is the analogue of Theorem G for the inverse function. 

\begin{proposition} \label{ginverse}
Suppose $\la \in L_{q,0}$ and $c=\hoi(\la) \in M\p$. Recall that $X_\la$ is a neighborhood of $K_\la^E$ bounded by the extended equipotential $\gamma_\la$ and that $W_c$ is a neighborhood of $K_c$, the filled level set of a chosen equipotential of $G_c$. There exists a homeomorphism
\[ H_\la: X_\la \longrightarrow W_c \]
which satisfies $\overline{\partial} H_\la=0$ on $K_\la^E$ and conjugates $f_\la\two$ to $Q_c$ and such that $H_\la$ maps relevant sectors and rays for $P_\la$ to their counterparts for $Q_c$, i.e. $H_\la$ maps $\estil_\la(\frac{i}{q+1})$, $R_\la(\frac{i}{q+1})$, and $\widetilde{R}_\la(\frac{i}{q+1})$ to $S_c(\thetatil^{i-1})$, $R_c(\theta^{i-1})$, and $R_c(\thetatil^{i-1})$ respectively.

Moreover, $H_\la|_{K_\la^E}$ is uniquely determined and it conjugates $f_\la\one$ to $Q_c$.
\end{proposition}
The proof is analogous to the proof of Theorem G.

Substituting $c$ by $\la$, remarks \ref{conju} hold for $\hoi$.

\begin{proposition}\label{injective} 
\ \\
\vspace{-.7cm}
\begin{enumerate}
\item $\ho \circ \hoi = \text{Id}_{L_{q,0}}$
\item $\hoi \circ \ho = \text{Id}_{M_{p/q}}$
\end{enumerate}
Hence $\ho$ is bijective.
\end{proposition}

\begin{proof}{}
To prove the first statement start with $\la \in L_{q,0}$ and, following the construction above, obtain a polynomial $Q_c$ and the map $H_\la$ given by Prop.~\ref{ginverse}, where $c=\hoi(\la)$. We now apply the construction in section \fullref{polike} to $Q_c$ choosing
\begin{align*}
	\gamma_c	&=	H_\la(\gamma_\la) \cap \bc_c^T\\
	\gamma'_c&=	H_\la(P_\la^{-1}(\gamma_\la)) \cap \bc_c^T\\
	f_{c,0,i}\two &= H_\la \circ f_{\la,0,i}\two  \circ H_\la^{-1}
\end{align*}
Let $\la'=\ho(c)$ which, by the remarks \ref{welldefined} on page \pageref{welldefined},  does not depend on these choices. Let $H_c:X_c \longrightarrow W_{\la'}$ be given by Theorem G.  It is easy to check that $H_c \circ H_\la |_{W_\la^T}$ is a hybrid equivalence between $P_\la$ and $P_{\la'}$. Hence $\la=\la'=\ho(\hoi(\la))$ by Proposition \ref{choices}.

To prove the second statement choose $c\in M_{p/q}$ and set $\la=\ho(c)$ and $c'=\xi_{p/q}(\la)$. We must show that $Q_c \sim_{hb} Q_{c'}$ and therefore $c=c'$.

Let $H_c$ be given by Theorem G. Besides the properties of $H_c$ stated in the theorem and since $H_c$ conjugates $f_c\two$ to $P_\la$ we have that
\begin{itemize}
\item the curves $(f_c\two)^{-n}(\gamma'_c)$ are mapped to equipotentials of $P_\la$ of potential $\eta' (q+1)^{(n+1)/q}$,
\item the rays landing at $\alpha'_c$ are mapped accordingly to the rays landing at $-q$. Moreover, if potential $\rho$ on $R_c(\thetatil)$ corresponds to potential $\rho'$ on $R_\la(0)$, then for all $i \in \br_+$, potential $\rho 2^{-i}$ in $R_c(\thetatil)$ corresponds to potential $\rho' (q+1)^{-i/q}$ on $R_\la(0)$.
\end{itemize}
We provide the idea of the construction, leaving the details to the reader. 

To define the domain of $f_\la\two:X'_\la \longrightarrow X_\la$ choose
\begin{align*}
	\gamma_\la \cap V_\la^{0\ldots q-1} &= H_c(\gamma_c) \\
	\gamma'_\la \cap V_\la^{0\ldots q-1} &= H_c(Q_c^{-1}(\gamma_c)),
\end{align*}
and complete the curves in $\Vtil_\la^i$ as usual. 

We now extend $H_c:X_c \rightarrow W_\la$ to $H_c^E:W_c \rightarrow X_\la$ as 
\[
	H_c^E=
		\begin{cases}
			\tau_\la^{ip} \circ P_\la|_{V_\la^{ip}}^{-1} \circ H_c \circ Q_c^{q-i} & \text{if $z \in \Vtil_c^{ip}$} \\
			H_c	& \text{otherwise}
		\end{cases}
\]
One can check that $H_c^E$ is well defined along the rays landing at $0$. It satisfies $\overline{\partial}H_c^E=0$ on $K_c$.

Finally, we define $f_\la\two$ on the sectors where the first return map is modified so that 
\[
	f_\la\two=H_c^E \circ Q_c \circ H_c^{-1}
\]
and $f_\la\one$ elsewhere. One can check that this map is a valid choice and that the following diagram
\[
	\begin{CD}
		W_{c,1} @>{Q_c}>> W_c \\
		@V{H_c^E}VV	@VV{H_c^E}V \\
		X'_\la @>{f_\la\two}>> X_\la
	\end{CD}
\]
commutes.

From the above, $H_c^E$ is a hybrid equivalence between $Q_c$ and $f_\la\two$. Since $c'$ does not depend on the choices, we have that $f_\la\two$ is hybrid equivalent to $Q_{c'}$. Hence $Q_c$ and $Q_{c'}$ are hybrid equivalent.
\end{proof}

%*****************************************************************************

\subsection{Analyticity of $\ho$ in the interior of $M_{p/q}$} \label{analyticity} 

Let $\Omega_M$ be a hyperbolic component of $M_{p/q}$. Then $\phi_{p/q}$ maps $\Omega_M$ to $\Omega_L$, a hyperbolic component of $L_{q,0}$ (as defined in section \ref{polynomials}). By remark \fullref{conju}, the multiplier of the corresponding attracting cycles is preserved. Let $\rho_{\Omega_M}:\bd \rightarrow \Omega_M$ (resp. $\rho_{\Omega_L}: \bd \rightarrow \Omega_L$) be the multiplier function that parametrizes $\Omega_M$ (resp. $\Omega_L$). Then, by the observation above
\begin{equation} \label{multiplier}
	\ho \mid_{\Omega_M}=\rho_{\Omega_L} \circ \rho_{\Omega_M}^{-1} \quad \text{and} \quad \hoi \mid_{\Omega_L}=\rho_{\Omega_M} \circ \rho_{\Omega_L}^{-1} 
\end{equation}
and hence $\ho$ is holomorphic on all hyperbolic components.

If  the interior of $M$ were known to equal the union of its hyperbolic components we would be done. However, this being still a conjecture implies that there might be {\em non-hyperbolic components} in the interior of $M$ (also called {\em queer components}). Such components might exist similarly for $L_q$. In both cases, the filled Julia sets of polynomials that belong to non-hyperbolic components have empty interior due to the classification of Fatou components. 

Let $\otil_M$ be a non-hyperbolic component of $M\p$. In order to prove analyticity of $\ho$ in $\otil_M$ we will use a parametrization similar to the parametrization of a hyperbolic component given by the multiplier map. In what follows we describe how to obtain such a parametrization.

\begin{definition}{}
Let $X$ be a Riemann surface. A {\em line field} supported on $E \subset X$ is an $L^\infty$  Beltrami form $\mu$ supported on $E$ with $\| \mu \|_\infty=1$. 
\end{definition}

\begin{proposition} \label{queer}
A parameter  $c$ belongs to a non-hyperbolic component of $M$ if and only if $J_c$ has positive measure and carries an invariant line field.
\end{proposition}
(For the proof see \cite{sullivan} and in particular p.61 in \cite{mcmullen})

\begin{remark}
The same proposition is true for a polynomial $P_{q,\la}$, with $\la$ in a non-hyperbolic component of the interior of $L_q$.
\end{remark}

\begin{proposition} \label{resp}
Let  $\otil_M$ (resp.~$\otil_L$) be a non-hyperbolic component of $M$ (resp.~$L_q$), and let $c,c' \in \otil_M$ (resp.~$\la,\la' \in \otil_L$). Then, there exists a unique co-hybrid equivalence (see Sect.\ref{tools}) $\psi_{c {c'}} :\bc \rightarrow \bc$ (resp. $\psi_{\la \la'} :\bc \rightarrow \bc$) between $Q_c$ and $Q_{c'}$ (resp.~$P_\la$ and $P_{\la'}$) such that $\frac{\psi_{c {c'}}(z)}{z} \rightarrow 1$ (resp.~$\frac{\psi_{\la {\la'}}(z)}{z} \rightarrow (\frac{\la}{\la'})^{\frac{1}{q}})$ as $z \rightarrow \infty$.
Moreover, if $c_i \in \otil_M$ (resp.~$\la_i \in \otil_L$), $i=1,2$, and $\mu_i=\frac{\overline{\partial}\psi_{c c_i}}{\partial \psi_{c c_i}}$ (resp.~$\mu_i=\frac{\overline{\partial}\psi_{\la \la_i}}{\partial \psi_{\la \la_i}}$), $i=1,2$, then there exists $t \in \bc^\ast$ such that $\mu_1=t \mu_2$.
\end{proposition}
\begin{proof}{}
We give the idea of the proof in the case of $\otil_L$ (compare with \cite{mcmullen}). For the first part, set $\psi_{\la \la'}=\psi_{\la'} \circ \psi_\la^{-1} : \bc \setminus K_\la \longrightarrow \bc \setminus K_{\la'}$. By the $\la$-lemma, $\psi_{\la \la'}$ extends to a quasi-conformal map from $\bc$ to $\bc$ conjugating $P_\la$ to $P_{\la'}$. The uniqueness follows from the uniqueness of $\psi_\la$ and $\psi_{\la'}$ as defined in Sect.~\ref{polyn}. Finally, observe that 
\[
	\frac{\psi_{\la \la'}}{z} =\frac{g_{\nu'} g_{\nu}^{-1}(z)}{z} \rightarrow \frac{\nu}{\nu'}=(\frac{\la}{\la'})^{\frac{1}{q}} \quad \text{as $z \rightarrow \infty$}.
\]

To show the uniqueness of the Beltrami forms up to multiplication by a constant observe that if that were not the case we could construct from the two parameter family $t_1 \mu_1+t_2 \mu_2$ of Beltrami forms, an injective map $(t_1, t_2) \mapsto \la(t_1,t_2)$ giving rise to a two complex parameter family of polynomials of the form $P_{q,\la}$, which would be a contradiction.
\end{proof} 

\begin{corollary}
For $c \in \otil_M$ (resp.~$\la \in \otil_L$), an invariant line field on $K_c$ (resp.~$K_\la$) is unique up to rotation.
\end{corollary}

\begin{proposition} \label{rho}
Let $c \in \otil_M$ and let $\mu_c$ denote an invariant line field of $Q_c$. The Beltrami forms $(t \mu_c)_{t \in \bd}$ induce a conformal equivalence
\[
	\rho_{c,\mu_c}:\bd \longrightarrow \otil_M .
\]
Similarly, if $\la \in \otil_L$ and $\mu_\la$ denotes an invariant line field of $\pol$, the Beltrami forms $(t \mu_\la)_{t \in \bd}$ induce a conformal equivalence
\[
	\rho_{\la,\mu_\la}:\bd \longrightarrow \otil_L .
\]

\end{proposition}
\begin{proof}{}
We give the proof for the case of $\otil_L$ (compare with \cite{mcmullen}). For the  proof in the case of  $\otil_M$ just replace $\la$ and $P$ by $c$ and $Q$ respectively. 

For $t\in \bd$ consider the Beltrami form $t\mu_\la$. Then $\| t \mu_\la \|_\infty=|t|<1$ and $t \mu_\la$ is invariant under $P_\la$. Apply the Measurable Riemann Mapping Theorem (dependent on parameters) to obtain a quasi-conformal homeomorphism $\varphi_\la^t$ that integrates $t\mu_\la$ and such that it fixes $-q$, $0$ and $\infty$. The map $\varphi_\la^t\circ P_\la \circ (\varphi_\la^t)^{-1}$ is holomorphic and a degree $q+1$ polynomial of the form $P_{\la(t)}$. Define $\rho_{\la,\mu_\la}(t)=\la(t)$.

To show that $\la(t)$ is holomorphic note that the free critical point (resp.~free critical value) of $P_\la$ must be mapped under $\varphi_\la^t$ to the free critical point (resp.~free critical value) of $P_{\la(t)}$. Hence 
\[
	\la(t)=-\varphi_\la^t (-\la (\frac{q}{q+1})^{q+1}) (\frac{q+1}{q})^{q+1}.
\]
Since $t \mu_\la$ varies holomorphically with $t$, $\varphi_\la^t$ varies holomorphically with $t$ and therefore so does $\la(t)$.

To show that $\la(t)$ is injective, suppose that $\la(t_1)=\la(t_2)$. Then 
\[
	\varphi_\la^{t_1} \circ P_\la \circ (\varphi_\la^{t_1})^{-1}=\varphi_\la^{t_2} \circ P_\la \circ  (\varphi_\la^{t_2})^{-1}
\]
i.e. $P_\la=(\varphi_\la^{t_2})^{-1} \circ \varphi_\la^{t_1} \circ P_\la \circ  (\varphi_\la^{t_1})^{-1} \circ \varphi_\la^{t_2}$.

Hence $\varphi=(\varphi_\la^{t_2})^{-1} \circ \varphi_\la^{t_1}$ is a co-hybrid equivalence between $P_\la$ and itself. Let $\psi_\la$ be the B\"{o}ttcher coordinate on $\bc \setminus K_\la$. Then $\psi_\la^{-1} \circ \varphi \circ  \psi_\la : \bc \setminus \overline{\bd} \longrightarrow \bc \setminus \overline{\bd}$ is holomorphic and conjugates $z \mapsto z^{q+1}$ to itself. Moreover, $\varphi$ maps the cycle of rays landing at the repelling fixed point $0$ onto themselves. Hence $\varphi$ must be the identity. Since $\psi_\la$ is bijective, this implies $\varphi_\la^{t_1}=\varphi_\la^{t_2}$ on $\bc\setminus K_\la$.

Since $K_\la=\partial(\bc \setminus K_\la)$ and $\varphi_\la^{t_1}$ and $\varphi_\la^{t_2}$ are continuous, it follows that $\varphi_\la^{t_1}=\varphi_\la^{t_2}$ on $K_\la$ and therefore on all of $\bc$.

Finally, $\frac{\overline{\partial}\varphi_\la^{t_i}}{\partial \varphi_\la^{t_i}}=t_i \mu_\la$ for $i=1,2$, so $t_1 \mu_\la=t_2 \mu_\la$. Since $\mu_\la \neq 0$ on a set of positive measure we have that $t_1=t_2$. 

At this point we have obtained a holomorphic injective map from the unit disc to a neighborhood of $\la$ in $\otil_L$. To show that this neighborhood covers the whole component, pick a point $\la' \in \otil_L$. By Prop.~\ref{resp} there exists a co-hybrid equivalence $\psi_{\la \la'}$ between $P_\la$ and $P_{\la'}$ with $\mu' := \frac{\overline{\partial}\psi_{\la \la'}}{\partial \psi_{\la \la'}}$. Moreover, since $\varphi_\la^t$ is a co-hybrid equivalence between $P_\la$ and $P_{\la(t)}$ for any $t\in \bd$, we have that  $\mu'= t^\ast t \mu_\la$ for some $t^\ast \in \bc$. Hence $\psi_{\la \la'} \circ \psi_\la^{ t^\ast t}$ is a conformal equivalence between $P_{\la'}$ and $P_{\la(t^\ast t)}$ and therefore $\la'=\la(t^\ast t)$.  
\end{proof}

\begin{definition/lemma}
Let $c \in \otil_M$, $\la=\ho(c)$ and $H_c$ be a quasi-conformal homeomorphism as in Theorem G, uniquely defined on $K_c^T$. Given $\mu_c$ an invariant line field on $K_c$, let $\mu_c^T$ be the restriction of $\mu_c$ to $\bc_c^T$. Then the Beltrami form
\[ \mu_\la= 
		\begin{cases}
			(H_c^{-1})^\ast \mu_c^T 	& \text{on $ K_\la$}		\\
			0 							& \text{on $ \bc \setminus K_\la$.}
		\end{cases}
\]
defines an invariant line field in $K_\la$ which we call the {\bf induced line field on $K_\la$}. 
\end{definition/lemma}

\begin{proof}{}
Note that $\mu_\la$ is well defined since $H_c|_{K_c^T}$ is unique. The fact that $\|\mu_\la\|_\infty=1$ follows from  $\overline{\partial}H_c=0$ on $K_c^T$.

Let $E \subset K_c^T$ be the set of positive measure where $\mu_c^T\neq 0$. Since quasi-conformal mappings are absolutely continuous with respect to the Lebesgue measure, we have that $H_c(E) \subset K_\la$ has positive measure. Hence $\mu_\la \neq 0$ on a set of positive measure.

The invariance of $\mu_\la$ by $P_\la$ follows from the fact that $H_c$ conjugates $f_c\one$ to $P_\la$. By Prop.~\ref{queer}, $\la$ belongs to a non-hyperbolic component of $L_{q,0}$.  
\end{proof}

\begin{theorem} \label{hol}
Let  $\otil_M$ be a non-hyperbolic component of $M\p$. Then, $\otil_L := \phi_{p/q} (\otil_M)$ is a non-hyperbolic component of $L_{q,0}$. 

Moreover, if we fix $c \in \otil_M$, $\la=\ho(c)$ and $\mu_c$  an invariant line field on $K_c$, then
\[
	\ho|_{\otil_M}=\rho_{\la, \mu_\la} \circ \rho_{c,\mu_c}^{-1} 
\]
where $\mu_\la$ denotes the induced line field on $K_\la$;  hence $\ho$ is holomorphic on all non-hyperbolic components of $M$.
\end{theorem}

The rest of the section is dedicated to prove the theorem above.

Pick $c \in \otil_M$ and let $\mu_c$, $\rho_c$ and $\varphi_c^t$ be as in proposition \ref{rho}. To prove that $\ho$ is holomorphic in $\otil_M$ we will use notation from the construction of this map in section \fullref{polike}.

Let $(\varphi_c^t)^T=\varphi_c^t|_{\bc_c^T}$.

\begin{lemma}\label{mess}
The  diagram
\[
  \begin{CD}
	  X'_c   @>{f_c\two}>> X_c \\
	  @V{(\varphi_c^t)^T}VV  @VV{(\varphi_c^t)^T}V   \\
     X'_{c(t)} @>{f_{c(t)}\two}>> X_{c(t)} 
  \end{CD}
\]
commutes and $\frac{\overline{\partial}(\varphi_c^t)^T}{\partial (\varphi_c^t)^T}=t \mu_c^T$.
\end{lemma}
\begin{proof}{}
Let $\psi_c$ and $\psi_{c(t)}$ be as above. We have the following diagram
\[
  \begin{CD}
	  \bc \setminus \overline{\bd}  @>{z\mapsto z^2}>> \bc \setminus \overline{\bd} \\
	  @V{\psi_c}VV  @VV{\psi_c}V   \\
      \bc\setminus K_c @>{Q_c}>> \bc\setminus K_c\\ 
      @V{\varphi_c^t}VV   @VV{\varphi_c^t}V   \\
     \bc\setminus K_{c(t)} @>{Q_{c(t)}}>> \bc\setminus K_{c(t)}\\ 
     @A{\psi_{c(t)}}AA   @AA{\psi_{c(t)}}A   \\
     \bc \setminus \overline{\bd}   @>{z\mapsto z^2}>> \bc \setminus \overline{\bd} 
  \end{CD}
\]
By the same argument as before $\psi_{c(t)}^{-1} \circ \varphi_c^t \circ \psi_c=\mbox{Id}|_{\bc\setminus \overline{\bd}}$. It follows that $\varphi_c^t$ maps rays to rays and equipotentials to equipotentials. Hence $(\varphi_c^t)^T(\bc_c^T)=\bc_{c(t)}^T$.

Let 
\[
	(\bc\setminus \overline{\bd})^T=\{ \exp(\rho+2\pi i t) \in \bc\setminus \overline{\bd} \mid \rho+2\pi i t \in \bh^T\}
\]
and $\psi_c^T=\psi_c|_{(\bc\setminus \overline{\bd})^T}$. Then, $f\two :X' \rightarrow X$ ($X'\subset X \subset \bh^T$) induces a map $\ftil\two:\widetilde{X'} \rightarrow \widetilde{X}$ ($\widetilde{X'}\subset \widetilde{X} \subset (\bc\setminus \overline{\bd})^T$) such that the following diagram commutes
\[
  \begin{CD}
	  \widetilde{X'}   @>{\ftil\two}>> \widetilde{X} \\
	  @V{\psi_c^T}VV  @VV{\psi_c^T}V   \\
     X'_c @>{f_{c}\two}>> X_c 
  \end{CD}
\]
and the same is true replacing $c$ by $c(t)$. Then,
\[
  \begin{aligned}
	(\varphi_c^t)^T(f_c\two(\psi_c^T(z))) &= (\varphi_c^t)^T(\psi_c^T (\ftil\two(z)))\\
		&= \psi_{c(t)}^T (\ftil\two(z))\\
		&= f_{c(t)}\two(\psi_{c(t)}^T(z))\\
		&= f_{c(t)}\two((\varphi_c^t)^T	(\psi_{c(t)}^T(z)))\\
  \end{aligned}
\]
which proves that the diagram in the lemma commutes.

Since $(\varphi_c^t)^T$ is a restriction of $\varphi_c^t$  we have that $\frac{\overline{\partial}(\varphi_c^t)^T}{\partial (\varphi_c^t)^T}=t \mu_c^T$.  
\end{proof}

\begin{lemma}
The map $H:=H_{c(t)}\circ (\varphi_c^t)^T \circ H_c^{-1}$ conjugates $P_\la$ to $P_{\ho(c(t))}$. Moreover, $\frac{\overline{\partial}H}{\partial H}=t \mu_\la$ on $K_\la$.
\end{lemma}
\begin{proof}{}
The conjugacy follows from the construction of $\ho$ together with lemma \ref{mess}. Now consider the Beltrami form $t\mu_\la$ on $K_\la$. By $H_c$, $t \mu_\la$ is transported to $t \mu_c^T$ since $\overline{\partial}H_c=0$ on $K_c$. The map $(\varphi_c^t)^T$ integrates $t \mu_c^T$ which means that $t \mu_c^T$ is transported to the standard complex structure on $K_{c(t)}^T$. Now, $\overline{\partial} H_c=0$ on $K_{c(t)}$ so in summary, $t\mu_\la$ is transported to the standard complex structure by $H$.  
\end{proof}

We have proven that the following diagram commutes
\[
  \begin{CD}
	  N(K_{\ho(c(t))})   @>{P_{\ho(c(t))}}>> N(K_{\ho(c(t))}) \\
	  @A{H}AA @AA{H}A   \\
      N(K_\la) @>{P_\la}>>  N(K_\la)\\
	  @V{\varphi_\la^t}VV  @VV{\varphi_\la^t}V   \\
      N(K_{\la(t)}) @>{P_{\la(t)}}>> N(K_{\la(t)})
   \end{CD}
\]
Moreover, $\overline{\partial}(\varphi_\la^t \circ H^{-1})=0$ on $K_{\ho(c(t))}$ since both maps $H$ and $\varphi_\la^t$ integrate $t \mu_\la$. Hence, $\varphi_\la^t \circ H^{-1}$ is a hybrid equivalence between $P_{\la(t)}$ and $P_{\ho(c(t))}$ and therefore $\la(t)=\ho(c(t))$ by Cor.~\ref{choices}. 
\begin{flushright} {\bf q.e.d.~Thm.~\ref{hol}} \end{flushright}

\begin{remark}
One could state Theorem \ref{hol} in terms of  $\xi_{p/q}$ concluding that this map is also holomorphic. One should mimic the proof above considering the ``extended'' versions of $\varphi^t_\la$ and $\mu_\la$ as opposed to the ``truncated'' versions of $\varphi_c^t$ and $\mu_c$.
\end{remark}

%*****************************************************************************

\subsection{Continuity of $\ho$ at points on the boundary of $M_{p/q}$} \label{continuity} 

Assume that $c$ is in the boundary of the $p/q$-limb. The following
will be a key lemma:

\begin{lemma}\label{key} Suppose  $\la \in \partial L_q$,
$\la^\prime \in \bc $ and $P_\la \sim_{qc} P_{\la^\prime}$. Then
$\la=\la^\prime$.
\end{lemma}

\begin{proof}{} Let $\varphi$ be the quasi-conformal homeomorphism that
conjugates $P_\la$ to $P_{\la^\prime}$, i.e. 
\[  
	P_{\la'}=\varphi \circ P_\la \circ \varphi^{-1}
\]

Set $\mu=\frac{\overline{\partial} \varphi}{\partial \varphi}$. Since $\varphi$ is quasi-conformal, we have that $\| \mu \|_\infty <1$.

Define a new Beltrami form
$\widetilde{\mu}$ by
\[\widetilde{ \mu}  = 
		\begin{cases}
			\mu	&	\text{on $K_\la$} 		\\
			0	&	\text{on $\bc \setminus K_\la$}
		\end{cases}
\] 
Then $ k=\|\: \widetilde{\mu}\: \|_\infty \leq \|\: \mu \: \|_\infty <1.$

Consider the family of Beltrami forms $\mu_t=t \widetilde{\mu}$, for $|t|<1/k$. Then,
\[ 
	\| \: \mu_t\:  \|_\infty = |t|\:  \|\: \widetilde{\mu}\: \|_\infty = |t|\: k <1. 
\] 
We then apply the Measurable Riemann Mapping Theorem (dependent on parameters) to obtain the unique integrating maps $ \varphi_t: \bc \longrightarrow \bc $ such that $\frac{\overline{\partial} \varphi_t}{\partial \varphi_t}=  \mu_t$, and such that they fix $0$, $\infty$ and $-q$. Since $\mu_t$ depends analytically on $t$, the maps  $\varphi_t$ also depend analytically on $t$.

For each $t$, let $P_t=\varphi_t \circ P_\la \circ \varphi_t^{-1}$ which is a polynomial. Because of the choices made for $\varphi_t$, it is easy to check that $P_t$ must be of type $E$. Hence, $ P_t=P_{\la(t)}=P_{q,\la(t)} $ and $\varphi_t(\omega)=\omega$, where $\omega=-q/(q+1)$ is the simple critical point. Also, the orbit of $\omega$ by $P_{\la(t)}$ must be bounded, so $\la(t) \in L_q$.

By the same argument as in the proof of Prop. \ref{rho}, $\la(t)$ is an analytic function of $t$.

Now note that $\mu_0=0$, so $\varphi_0=\mbox{Id}$ and hence $\la(0)=\la$. Since
$\la(t)$ is analytic then, by the open mapping principle, it is either
open or constant. If it is open, a neighborhood of $t=0$ must be mapped
to a neighborhood of \lla, which is a contradiction since $\la \in
\partial L_q$ and $\la(t)$ belongs to $L_q$ for all $t$. Therefore
$\la(t)$ is constant and since $\la(0)=\la$ it follows that $\la(t) \equiv
\la$ for all $t$.

In particular $\la(1)=\la$. But for $t=1$ we have
$\mu_1=\widetilde{\mu}$. Hence, we have
\[  
\begin{CD}
	\bc   @>P_\la>> \bc 	\\
   @A{\varphi_1}AA @AA{\varphi_1}A   \\  
	\bc   @>P_\la>> \bc&		\\
   @V{\varphi}VV @VV{\varphi}V  \\
	\bc   @>P_{\la^\prime}>> \bc
\end{CD}	
\] 
and $\varphi \circ \varphi_1^{-1}$ is another quasi-conformal
homeomorphism that conjugates $P_\la$   to $P_{\la^\prime}$. In fact, we
claim that it is a hybrid equivalence since it pulls back the standard
complex structure on $K_{\la^\prime}$ to the standard complex structure
on $K_\la$. It follows that $P_\la \sim_{hb} P_{\la^\prime}$ and
therefore  by Prop.~\ref{choices}, $\la=\la^\prime$.  
\end{proof}

Given any sequence $c_n \rightarrow c$, let $\la_n=\ho(c_n)$ and
$\la=\ho(c)$. To prove continuity at $c$ we must show that
$\lim_{n\rightarrow \infty} \la_n =\la$.

We claim that it suffices to pick an arbitrary convergent subsequence
$\la_{n_k} \rightarrow \latil$, and show that $\latil=\la$. Indeed, since
$\la_{n_k}$ is arbitrary, it follows that all convergent subsequences
tend to $\la=\ho(c)$. Since $L_{q,0} \cup \{1\}$ is compact, it follows that
$\lim_{n\rightarrow \infty} \la_n =\la$.

So, from now on we assume that we have picked  a sequence $c_n
\rightarrow c \in \partial M_{p/q}$. Let  $\varphi_n$, $f^{(3)}_n$,
$\chi_n$, etc. be as in section \ref{polike}. To make the notation easier
we denote the convergent subsequence again by $\la_n \rightarrow \latil$.
We want to show that $\latil=\la$.

\begin{lemma}\label{dos} The  polynomials $P_\la$ and
$P_{\latil}$ are quasi-conformally conjugate.
\end{lemma}

\begin{proof}{} The mappings $\varphi_n$ are all quasi-conformal homeomorphisms and they have a uniform bounded dilatation ratio, since a new complex structure was constructed once and for all in the right half plane and then pulled back by $(\psi_\la \circ \exp)^T$. Let $K$ be this uniform bound. Then, $\{\varphi_n\}$ is a $K$-quasi-conformal family. Since the space of $K$-quasi-conformal mappings on a compact set is compact with respect to uniform convergence, we have that the $\{\varphi_n\}$ form an equi-continuous family. Therefore, by Arzela-Ascoli's theorem, we can find a convergent subsequence $\varphi_{n_k} \rightarrow
\varphi^\ast$. The limit map is a $K$-quasi-conformal homeomorphism. The map 
\[ 
	f^{(3)}_\ast = \varphi^\ast \circ f_c^{(2)} \circ (\varphi^\ast)^{-1}. 
\]
is analytic and 
\[
	f^{(3)}_\ast \overset{\varphi^\ast}{\sim}_{hb} f\two \overset{\varphi_c}{\sim}_{hb} f_c\three.	
\]  
Hence $f_c^{(3)} \sim_{qc} f^{(3)}_\ast $.

We now apply the lemma in p.~313 in \cite{dh3} to $f^{(3)}_\ast $ to obtain a subsequence
of hybrid equivalences $\{\chi_n \}$ (abusing notation again) that
converges to a quasi-conformal homeomorphism $\chi^\ast$ and a polynomial
$P^\ast$ such that $P_{\la_n} \rightarrow P^\ast$ and
\[ P_\ast= \chi^\ast \circ f^{(3)}_\ast \circ (\chi^\ast)^{-1} .\]
Then
\[  
	P^\ast \overset{\chi^\ast}{\sim}_{hb} f\three \overset{\chi_c}{\sim}_{hb} P_\la.	
\] 
so it follows that $P_\la \sim_{qc} P^\ast$.

Finally, note that $P^\ast$ must be of the form $P_{\la^\ast}$ (maybe after composing with an affine conjugacy).  Hence,  $P_{\la_n} \rightarrow P^\ast$ implies
that $\la_n \rightarrow \la^\ast$. It follows that $\la^\ast=\latil$ and
therefore we have proved $P_\la \sim_{qc} P_{\latil}$.  
\end{proof}

We now give the conclusion of the proof. 

Suppose the sequence $c_n \rightarrow c$ was such that all the $c_n$'s are Misiurewicz points. Then it is clear the $\la_n$'s are also Misiurewicz points and hence $\la_n \in \partial L_q$. Since the boundary is compact we have that any accumulation point of $\{\la_n\}$ must be on the boundary. So  $\latil \in \partial L_q$.

By  Lemma \ref{dos}, $P_\la \sim_{qc} P_{\latil}$. This fact, together with $\latil \in \partial L_q$  allows as to apply Lemma \ref{key} to obtain $\latil=\la$.

Now suppose the sequence $c_n \rightarrow c$ was  an arbitrary sequence, so that there is no assumption on the $c_n$'s. Choose also another sequence $c_n^\prime \rightarrow c$  formed by  Misiurewicz points, and let $\la^\prime$ be an accumulation point of $\{\la_n^\prime\}$. Then, by the above we have that $\la^\prime= \la$ and hence $\la \in \partial L_q$. Still, by lemma \ref{dos} we have $P_\la \sim_{qc} P_{\latil}$. This fact, together with $\la \in \partial L_q$, allows us to use again Lemma \ref{key} to conclude $\la=\latil$.

This concludes the proof of continuity at points on the boundary.

This concludes the proof of Theorem F.

\begin{remark} \label{roots}
It follows that (\ref{multiplier}) is satisfied on the closure $\overline{\Omega_M}$, so if $c$ has internal argument $t$ then $\la=\ho(c)$ has also internal argument $t$. In particular, roots of hyperbolic components are mapped to roots of hyperbolic components.
\end{remark}
%*****************************************************************************
\subsection{Proof of theorem A} \label{conclusion} 

At this point, for every $p/q$, we have constructed a homeomorphism
\[ \ho:M_{p/q} \longrightarrow L_{q,0} \] which is analytic in the interior of $M_{p/q}$.

Now, for any $p/q$ and $p^\prime/q$, define the map $\finho$ to be
\[ 
\finho: M_{p/q} @>{\ho}>> L_{q,0} @>{\phi_{p^\prime/q}^{-1}}>> M_{p^\prime/q} 
\].
 
Thus, $\finho$ is a homeomorphism between $M_{p/q}$ and $M_{p^\prime/q}$ which by Theorem  \ref{hol} is holomorphic in the interior of the limbs.

For a given $c\in M_{p/q}$, the second part of Theorem A claims the existence of a homeomorphism $\bfihat_c$ between a neighborhood of $K_c$ and a neighborhood of $K_{c'}$ where $c'=\finho(c)$. 

Let $\la=\ho(c)=\phi_{p'/q}(c)$. Let $H_c:X_c \longrightarrow W_\la$ and $H_{c'}:X_{c'} \longrightarrow W_\la$ be  quasi-conformal homeomorphisms as in Theorem G. 
Recall that these maps conjugate the first return maps to the polynomial $P_\la$ outside of the sectors and send the sector $S_c(\thetatil)$ (resp~$S_{c'}(\thetatil)$) around the identified ray ($R_c(\thetatil)$, resp.~$R_{c'}(\thetatil)$) to the sector $S_\la(0)$ around $R_\la(0)$. Define,
\[ 
	\bfihat_c^T:X_c\setminus S_c(\thetatil) @>{H_c}>>  W_\la\setminus S_\la(0)
@>{H_{c^\prime}^{-1}}>> X_{c'}\setminus S_{c'}(\thetatil).
\]
Finally,  extend $\bfihat_c^T$ to $W_c$ - the filled level set of the chosen potential -  by defining 
\[ 
	\bfihat_c(z)= 
		\begin{cases}
			\bfihat_c^T(z) & \text{if defined} \\
			-\bfihat_c^T(z) & \text{otherwise}
		\end{cases}
\]
\begin{lemma}
The map $\bfihat_c$ is well defined.
\end{lemma}
\begin{proof}{}
We must check that along the boundaries of $S_c(\thetatil^{0})$ and $S_c(\thetatil^{q-1})$ in $V_c^0$ we have
\[	\bfihat_c^T(-z)=-\bfihat_c^T(z) .\]
We will actually prove it for all points in $V_c^0$ not in the interior of the sectors. We know that in this region $\bfihat_c^T$ conjugates $Q_c^q$ to $Q_{c'}^q$. Hence
\[
	\bfihat_c^T(Q_c^q(z))=Q_{c'}^q(\bfihat_c^T(z)) 
\]
and the same holds for $-z$ since $-z \in V_c^0$. But $Q_c^q(z)=Q_c^q(-z)$ so 
\[
	Q_{c'}^q(\bfihat_c^T(z))=Q_{c'}^q(\bfihat_c^T(-z)).
\]
Since $Q_c^q$ is two-to-one on $V_c^0$ the conclusion follows.
\end{proof} 

This concludes the proof of Theorem A.

\begin{remark} The map $\bfihat_c$ is a quasi-conformal homeomorphism satisfying $\overline{\partial} \bfihat_c=0$ on $K_c$. This map does not conjugate $Q_c$ and $Q_{c'}$ but  it conjugates the first return map on $K_c^T$ to the first return map on $K_{c^\prime}^T$.
\end{remark}
%*****************************************************************************
%*****************************************************************************

\subsection{Compatibility with Tuning} \label{tuning}

The Mandelbrot set contains copies of itself. For each center $c_0$ of a hyperbolic component there is associated a copy of $M$ obtained by tuning (see \cite{douady}):
\[
	\begin{array}{ccl}
		M & \longrightarrow & M_{c_0}=c_0 \perp M \\
		c & \longrightarrow & c_0 \perp c
	\end{array}
\]
The tuning map is a homeomorphism mapping centers of hyperbolic components to centers of hyperbolic components and Misiurewicz points to Misiurewicz points;  it is analytic in the interior of $M$ and $\partial M_{c_0} \subset \partial M$. Moreover, Misiurewicz points are dense in $\partial M$ (see \cite{bodil}), hence dense in $\partial M_{c_0}$.

\begin{proposition} 
Let $c_0$ be a  center of a hyperbolic component in $M_{p/q}$. The homeomorphism $\finho$ is compatible with tuning, i.e.
\begin{equation} \label{tunning}
	\finho (c_0 \perp c) = \finho (c_0) \perp c
\end{equation}
for all $c \in M$.
\end{proposition}

\begin{figure}[htbp]
%\vspace{2cm}
\centerline{\epsfysize=6cm\epsffile{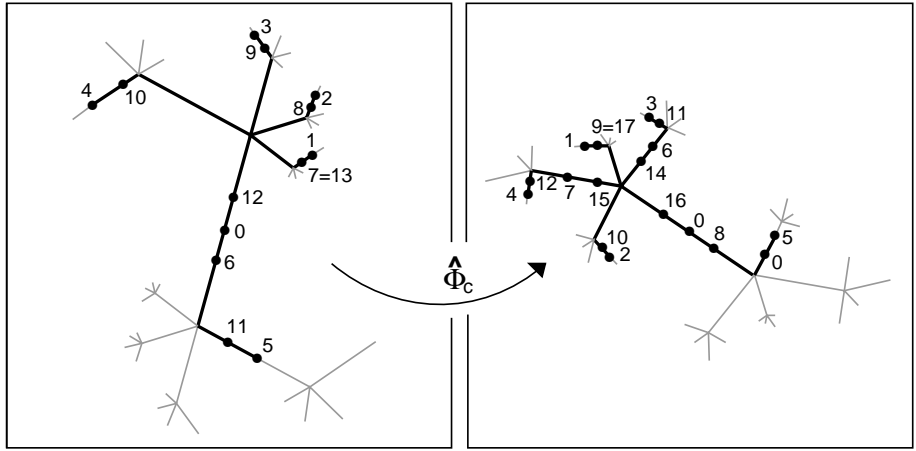}}
\caption{\small Left: the Hubbard tree for $c_0 \perp (-2)$, where $c_0$ is the center of the period 6 hyperbolic component in $M_{1/5}$ (as in Fig.~\ref{juliafifths} left). Right: the Hubbard tree for $\phi_{12}^5(c_0) \perp (-2)$, where $\phi_{12}^5(c_0)$ is the center of the period 8 hyperbolic component (as in Fig.~\ref{juliafifths} right). Compare also with Fig.~\ref{treesfifths}. }
\label{treesfifths2}
\end{figure}

\begin{proof}{}
The parameter $c'_0=\finho(c_0)$ is the center of a hyperbolic component in $M_{p'/q}$. Following the surgery constructions it is very easy to see that eq.~(\ref{tunning}) is satisfied for any Misiurewicz point $c \in M$ (compare the Hubbard trees, see Figs.~\ref{treesfifths} and \ref{treesfifths2} for an example), hence satisfied for all $c \in \partial M$, so $\finho (\partial M_{c_0}) = \partial M_{c'_0}$. The proof is completed through the following topological lemma:
\begin{lemma}
Let $X_0 \subset X \subset \bc$ and $X'_0 \subset X' \subset \bc$. Suppose $X_0$ and $X'_0$ are compact, connected and full in $\bc$ and satisfy $\partial X_0 \subset \partial X$ and $\partial X'_0 \subset \partial X'$. If $\phi:X \longrightarrow X'$ is a homeomorphism satisfying $\phi(\partial X_0)=\partial X'_0$ then $\phi (X_0)=X'_0$.
\end{lemma} \ 
\end{proof}

%*****************************************************************************
%*****************************************************************************
%*****************************************************************************

\section{Involutions: Proof of Theorem B} \label{involutions} 

If we look at the pictures of the connectedness loci $L_q$ (see Figs. \ref{pol1} to \ref{pol19}), we observe a symmetry with respect to complex conjugation. Indeed, the function
\[ \begin{array}{rccc}  C: &L_q  & \longrightarrow & L_q \\
        & \la & \longmapsto	& \overline{\la}
\end{array}
\] is a homeomorphism, and $C^2=Id$. The map $C$ maps the $0$-limb to itself and the polynomials $P_{q,\la}$ and $P_{q,\overline{\la}}$ are conjugate through the anti-holomorphic map $C:z \mapsto \overline{z}$.
 
It is then clear that the function
\[ 
	\ical_{p/q}:M_{p/q} @>{\ho}>> L_{q,0} @>{C}>> L_{q,0} @>{\ho^{-1}}>> M_{p/q} 
\]
is a non-trivial homeomorphism from the $p/q$-limb of $M$ to itself, such that
$\ical_{p/q}^2=Id$.

Since $\ho$ is holomorphic in the interior of $M\p$ it follows that $\ical\p$ is anti-holomorphic in the interior of $M\p$.

From the construction, we see that $\ical_{p/q}$ is like a reflection
with respect to an arc of symmetry.  This arc is the preimage of
the real axis (up to the ``last'' Misiurewicz point in $L_{q,0}$) under
$\ho$, and hence it is always a topological arc (see Figs. \ref{inv1} and \ref{inv1-19}). 

The involution $\ical\p$ can also be obtained by using the symmetry of the Mandelbrot set and the homeomorphisms between limbs given by Theorem A. The map $C$ maps the limb $M\p$ onto  the limb $M_{(q-p)/q}$ and the polynomials $Q_c$ and $Q_{\overline{c}}$ are conjugate through the anti-holomorphic map $C:z \mapsto \overline{z}$. We have the following lemma:

\begin{lemma} \label{compat}
The homeomorphisms  $\ho:M\p \longrightarrow L_{q,0}$ are compatible with complex conjugation, i.e.~the diagram 
\[
	\begin{CD}
		M\p @>{\ho}>> L_{q,0} \\
		@VCVV  @VVCV			\\
		M_{(q-p)/q} @>{\phi_{(q-p)/q}}>> L_{q,0}
	\end{CD}
\]
commutes.
\end{lemma}
\begin{proof}{}
For this proof we use notation from section \ref{polike}. Fix the potential $\eta>0$ and the slope $s$. The first return map satisfies $f_{\overline{c}}\one (z)=\overline{f_c\one(\overline{z})}$. We can then modify the maps in a similar fashion, i.e.~ set $f_{\overline{c}}\two (z)=\overline{f_c\two(\overline{z})}$. The map $\varphi_{\overline{c}}(z)=\overline{\varphi_c(\overline{z})}$ is an integrating map and it defines the polynomial-like map $f_{\overline{c}}\three$ to be
\[	 f_{\overline{c}}\three (z)=\overline{f_c\three(\overline{z})}	\]
Finally, the map $\chi_{\overline{c}}(z)=\overline{\chi_c(\overline{z})}$ is a hybrid equivalence between $f_{\overline{c}}\three$ and $P_{q,\overline{\la}}$.
\end{proof}

From the lemma above it follows that
\begin{align*}
	\ical\p = \ho^{-1} \circ C \circ \ho &= \ho^{-1} \circ \phi_{(q-p)/q} \circ C  = \Phi_{(q-p)p}^q \circ C\\
			&=C \circ \phi_{(q-p)/p}^{-1} \circ \ho = C \circ \Phi^q_{p(q-p)}.
\end{align*}

Note that without reference to the polynomials $\pol$ it is not clear that the fixed points under $\ical\p$ form a topological arc.

\begin{remark} The map $\ical_{p/q}$ maps hyperbolic components to hyperbolic components and the corresponding multiplier maps are
complex conjugates of each other.
\end{remark}

Very much in particular, and as an example, the map $\ical_{1/3}$ gives a homeomorphism between the two main antennas in the $1/3$-limb, sending $c=i$ to the landing point of $R_M(1/4)$, and the center of the hyperbolic component of period 4 to the center of the hyperbolic component of period 5. 

The dynamical statement of Theorem B follows similarly. Set $\ical=\ical\p$. On one hand, for $c \in M\p$  and $\la=\ho(c)$ we can define 

\[  
	\begin{array}{rccccccc} 
  		\icalhat_c^T:& N(K_c^T) &\stackrel{H_c}{\longrightarrow} & N(K_\la) &\stackrel{\mathcal C}{\longrightarrow}	& N(K_{\overline{\la}}) & \stackrel{H_{\ical(c)}^{-1}}{\longrightarrow}  & N(K_{\ical(c)}^T) \\
 	   & z	& \longmapsto 	& H_c(z) & \longmapsto & \overline{H_c(z)}& \longmapsto 	& H_{\ical(c)}^{-1}(\overline{H_c(z)}) 
	\end{array}
\]
where   $H_c$ and $H_{\ical(c)}$ are as in Theorem G.  The map $\icalhat_c^T$ is an orientation reversing homeomorphism between $N(K_c^T)$ and $N(K_{\ical(c)}^T)$ (which in fact conjugates the two first return maps $f^{(1)}_c$ and $f^{(1)}_{\ical(c)}$ on $K_c^T$ and $K_{\ical(c)}^T$ respectively). The map $\icalhat_c^T$ can be extended to a homeomorphism $\icalhat_c$ (similarly to the extension of $\bfihat_c^T$ to $\bfihat_c$) between neighborhoods $N(K_c)$ and $N(K_{\ical(c)})$ respectively, as stated in Theorem B.

On the other hand, it follows from the lemma above that the maps in dynamical plane are compatible with complex conjugation. That is, for any $c \in M\p$ the following diagram commutes:
\[
	\begin{CD}
		N(K_c) @>{\bfihat_c}>>  N(K_{\Phi_{p(q-p)}^q}) \\
		@VCVV  @VVCV  \\
		N(K_{\overline{c}}) @>{\bfihat_{\overline{c}}}>> N(K_{\ical\p(c)})
	\end{CD}
\]
so $\icalhat_c:N(K_c) \longrightarrow N(K_{\ical\p(c)})$ can alternatively be defined as
\[
	\icalhat_c=\bfihat_{\overline{c}} \circ C = C \circ \bfihat_c.
\]

%*****************************************************************************
%*****************************************************************************
%*****************************************************************************
\section{Combinatorial Surgery: \\ Proof of Theorems H, C, D and E} \label{combsurg} 

\subsection{Proof of Theorems H and C} \label{theoremc}
The main goal in this section is to prove Theorem C. The strategy is the same as for proving Theorem A, namely via the higher degree polynomials $\pol$ with $\la \in L\qo$.

We first introduce the necessary notation and formulate Theorem H, the analogue of Theorem C building the bridge to the higher degree polynomials. 

In Sect.~ \ref{dynplane} we prove the part of Theorem H which deals with the dynamical planes and deduce the analogous part of Theorem C.

In Sect.~\ref{param} we prove the part of Theorem H which deals with the parameter spaces and deduce the analogous part of Theorem C.

Let $c \in M\p$ and $\la =\ho(c) \in L\qo$. Recall from Theorem G that there is a quasi-conformal homeomorphism 
\[
	H_c : X_c \longrightarrow W_\la
\]
conjugating $f_c\two$ and $P_\la$. In particular, it is uniquely defined on $K_c^T$ and conjugates the first return map $f_c\one$ and $P_\la$ on the filled Julia sets.

In this section we extend this conjugacy to the exterior of the filled Julia sets in a ``combinatorial way'' which is explained below.

Let $m_d: \bt \longrightarrow \bt$ denote the map induced on the circle $\bt=\br / \bz$ by $\mathcal{M}_d:\br \longrightarrow \br$, where $\mathcal{M}_d(x)=dx$ is multiplication by $d$.

Set $\bi\p:=[\thetatil^{q-1},\thetatil^0]$ as a subset of $\bt$ (and hence with the metric of $\bt$). Define the truncated circle to be
\[	\bt^T=\bt\p^T := \bi\p / \sim   \]
where $\sim$ denotes the equivalence relation identifying the endpoints of $\bi\p$. Denote as before the identified argument by $\thetatil$.

Let $m_2\one=m_{p/q,2}\one : \bt\p^T \longrightarrow \bt\p^T$ denote the first return map of the doubling map $m_2$, i.e.
\[  
	m_2\one(\theta)= 
	\begin{cases} 
  		2^q \theta    &  \text{if $\theta^{q-1} \leq\theta \leq
\theta^{0}$} 	\\
  		2^{q-1}\theta	&  \text{if $\theta^{p-1} \leq\theta \leq
\theta^{p}$}	\\
  		2^{q-2}\theta	&  \text{if $\theta^{2p-1} \leq\theta \leq
\theta^{2p}$}\\
  		\ldots			& 						\\
  		2 \theta &  \text{if $\theta^{q-p-1} \leq\theta \leq \theta^{q-p}$}
 	\end{cases}
\]
The first return map $m_2\one$ is an orientation preserving, expanding covering map of degree $q+1$, fixing $\thetatil$. Moreover, consider
\[	m_{q+1} : \bt \longrightarrow \bt	. \]
The map induced from multiplication by $q+1$ is also an orientation preserving, expanding covering map of degree $q+1$, fixing 0.

The lemma below is the tool we need in order to define the combinatorial extension of $H_c$. 
\begin{lemma} \label{expanding}
Let $\bt_j$, $j=1,2$ be a metric space homeomorphic to $\bt$ and with basepoint $t_j \in \bt_j$. Suppose $E_j :\bt_j \longrightarrow \bt_j$, $j=1,2$ is an orientation  preserving, expanding covering map of degree $d>1$, fixing $t_j$. Then, there exists a unique orientation preserving homeomorphism $T:\bt_1 \longrightarrow \bt_2$, conjugating $E_1$ and $E_2$ and mapping $t_1$ to $t_2$.
\end{lemma}

Hence, it follows that there exists a unique orientation preserving homeomorphism
\begin{equation}  \label{star}
	\bigthat = \bigthat\p := \bt\p^T \longrightarrow \bt
\end{equation}
conjugating $m_2\one$ and $m_{q+1}$ and mapping $\thetatil$ to 0. Note that $\bigthat$ is defined independently of $c \in M\p$ and $\la \in L\qo$.

\begin{proof}{of Lemma \ref{expanding}}
We give a sketch of the proof. Choose any orientation preserving homeomorphism $T_0 : \bt_1 \longrightarrow \bt_2$ which maps $t_1$ to $t_2$. Define by induction $T_n : \bt_1 \longrightarrow \bt_2$ to be the homeomorphism which makes the diagram 
\[	\begin{CD}
		\bt_1 @>{E_1}>> \bt_1 \\
		@V{T_n}VV   @VV{T_{n-1}}V \\
		\bt_2 @>{E_2}>>  \bt_2
	\end{CD}
\]
commutative and which maps $t_1$ to $t_2$. One can prove, using that $E_i$ is expanding, that $T_n$ converges to a homeomorphism $T$ with properties as claimed.
\end{proof}

For a given $p/q$ recall that $\theta\p^-=\thetatil^{p-1}$ and $\theta\p^+ = \thetatil^p$ are the arguments of the external rays of $M$ landing at the root point of the limb $M\p$. Moreover, recall that we have defined 0 and 1 to be the arguments relative to the 0-wake $WL\qo$ of the external rays of $L_q$ landing at the root point $\la=1$ of the 0-limb $L\qo$. Note that the interval $[\theta\p^-, \theta\p^+]=[\thetatil^{p-1}, \thetatil^p]$  is mapped homeomorphically onto the interval $\bi\p=[\thetatil^{q-1}, \thetatil^0]$ by the $(q-1)$-st iterate of the doubling map
\[	m_2^{q-1} :  [\theta\p^-, \theta\p^+] \longrightarrow \bi\p . \]
Define the homeomorphism 
\[  \Theta\p: [ \theta\p^-, \theta\p^+ ] \longrightarrow [0,1] \]
as
\begin{equation} \label{star2}
	\Theta\p = \bigthat\p \circ m_2^{q-1}
\end{equation}
on the open interval $(\theta\p^-, \theta\p^+)$.

\begin{theorem*}{H}
Given $p/q$, let $\bigthat\p$ and $\Theta\p$ be as defined in eqs.~(\ref{star}) and (\ref{star2}) respectively.
\begin{enumerate}
	\item Suppose $\theta=\frac{r}{s} \in [\theta\p^-, \theta\p^+]$ with $\gcd(r,s)=1$ and $s$ even. Then the ray $R_M(\theta)$ lands at a point $c \in M\p$ if and only if $R_{L\qo}(\Theta\p(\theta))$ lands at $\ho(c) \in L\qo$.

	\item Suppose $c \in M\p$ and $\la = \ho(c) \in L\qo$. Recall that a hybrid equivalence $H_c$ between $f_c\three$ and $P_\la$ is uniquely determined on $K_c^T$. For any $\theta \in \bt\p^T \setminus \{ \thetatil \}$ the ray $R_c(\theta)$ lands at a point $z \in K_c^T$ if and only if the ray $R_\la(\bigthat\p (\theta))$ lands at $H_c(z) \in K_\la$.
\end{enumerate}  
\end{theorem*} 

%***************************************************
\subsubsection{Dynamical Plane} \label{dynplane}

Let $c \in M\p$ and $\la = \ho(c) \in L\qo$. Moreover let $H_c :X_c \longrightarrow W_\la$ be given as in Theorem G.

\begin{proposition} \label{Bigtheta}
The map $H_c$ induces an orientation preserving homeomorphism
\[ 
	\Theta_c : \bt\p^T \longrightarrow \bt
\]
satisfying:
\begin{enumerate}
\item $\Theta_c(\thetatil)=0$.
\item For any $\theta \in \bt\p^T \setminus \{ \thetatil \}$, the external ray $R_c(\theta)$ lands at a point $z \in K_c^T$ if and only if the external ray $R_\la (\Theta_c (\theta))$ lands at $H_c(z) \in K_\la$.
\item $\Theta_c$ conjugates $m_2\one$ and $m_{q+1}$, hence $\Theta_c = \bigthat\p$.
\end{enumerate}
\end{proposition}

From part 3 follows this corollary:
\begin{corollary} \label{rationals}
$\Theta_c$ maps periodic (resp.~preperiodic) arguments under $m_2\one$ to periodic (resp.~preperiodic) arguments under $m_{q+1}$.
\end{corollary}

\begin{proof}{of \ref{Bigtheta}}
This would follow immediately from a stronger result called The Prime Ends Theorem (\cite{pomerenke}). However, we give here a proof for this weaker version. In order to define the map $\Theta_c$ we need some preliminaries.

Given $K$ a compact, connected and full subset of $\bc$, a point $z\in K$ is called {\em accessible} if there exists a curve $\gamma:[0,\varepsilon] \longrightarrow \bc$ with $\gamma(0)=z$ and $\gamma((0,\varepsilon]) \subset \bc \setminus K$. The homotopy class of such curves is call an {\em access} to $z$. If $K$ is locally connected, then all points in $K$ are accessible.

We shall need the following result:

Suppose $z$ is accessible. Then, to each access $[\gamma]$ there is associated a unique external ray  $R(\theta)$ that lands at $z$ and such that an arc of the ray belongs to the homotopy class $[\gamma]$ (see \cite{carsten,milnor}).

We define $\Theta_c$ first in the set of arguments that correspond to accesses. Let
\begin{eqnarray*}
	\cala(K_c^T)  &=& \{\theta \in \bt\p^T \mid \theta \mbox{\ corresponds to an access of } K_c^T \},\\
    \cala(K_\la)  &=& \{\theta \in \bt \mid \theta \mbox{\ corresponds to an access of } K_\la\}
\end{eqnarray*}
Note that if $K_c^T$ (and hence $K_\la$) are locally connected, then $\cala(K_c^T)= \bt\p^T$ and $\cala(K_\la)=\bt$.

\begin{definition}{of $\Theta_c$ on $\cala(K_c^T)$} \\ Given $\theta \in \cala(K_c^T)$, let $\gamma(\theta)$ be the curve which is the connected component of $\overline{R_c(\theta)} \cap N(K_c^T)$ containing $z=R_c^\ast(\theta)$ (the landing point of $R_c(\theta)$). Then, $H_c(\gamma(\theta))$ is a curve landing at $H_c(z)$. Define $\Theta_c(\theta)$ as the external argument of the external ray  associated with $[H_c(\gamma(\theta))]$, $R_\la(\Theta_c(\theta))$.
Rays landing at $H_c(z)$ have arguments in $\cala(K_\la)$, so $\Theta_c(\theta) \in \cala(K_\la)$. 
\end{definition}

Note that  $\Theta_c(\thetatil)=0$.

\begin{lemma}\label{primeends}
The map $\Theta_c$ is an orientation preserving bijection between $\cala(K_c^T)$ and $\cala(K_\la)$.
\end{lemma}

\begin{proof}{}
To show injectivity, take two elements in $\cala(K_c^T)$, $\theta \neq \theta^\prime$ and first suppose that $R_c(\theta)$ and $R_c(\theta^\prime)$ land at different points. Then, since $H_c$ is a homeomorphism,  $R_\la(\Theta_c(\theta))$ and $R_\la(\Theta_c(\theta^\prime))$ must land at different points. Hence $\Theta_c(\theta) \neq \Theta_c(\theta^\prime)$. Now suppose that $R_c(\theta)$ and $R_c(\theta^\prime)$ land at the same point. Then  $\gamma(\theta)$ and $\gamma(\theta^\prime)$ are not homotopic which implies that $H_c(\gamma(\theta))$ and $H_c(\gamma(\theta'))$ are not homotopic. Hence $\Theta_c(\theta) \neq \Theta_c(\theta^\prime)$.

To show surjectivity, let $R_\la(\theta^\prime)$ be a ray landing at $z^\prime\in K_\la$ so that $\theta^\prime \in \cala(K_\la)$, and let $\gamma(\theta^\prime)$ be the curve as in the definition above. Then, $H_c^{-1}(\gamma(\theta^\prime))$ is a curve landing at $H_c^{-1}(z^\prime)$. Let $\theta \in \cala(K_c)$ be the external argument of the ray associated with $[H_c^{-1}(\gamma(\theta^\prime))]$. It is clear that $\Theta_c(\theta)=\theta^\prime$, since the homotopy classes are preserved by the homeomorphism.

The fact that $\Theta_c$ is orientation preserving follows from $H_c$ being orientation preserving.  
\end{proof} 

\begin{corollary}{}
For every accessible point $z\in \cala(K_c^T)$, the number of rays landing at $z$ equals the number of rays landing at $H_c(z) \in K_\la $.
\end{corollary}

The map $\Theta_c$ is now defined from a dense set in $\bt\p^T$ to a dense set in $\bt$. Note that we can also consider it as a map from a dense set in $[\thetatil^q,\thetatil^0]$ to a dense set in $[0,1]$ by defining
\[ 
	\Theta_c(\thetatil^q)=0 \mbox{ \ \ and \ \ } \Theta_c(\thetatil^0)=1.
\]

The next step is to extend it to the whole interval. 

If $K_c^T$ is locally connected, then $\Theta_c$ is already defined in the whole interval. Moreover, since it is bijective and order preserving between two compact sets of $\br$, it follows that $\Theta_c$ is a homeomorphism and we are done.

If $K_c^T$ is not locally connected we will use the following topological lemma to extend $\Theta_c$. 

\begin{lemma} \label{extension} Let  $S$ be a dense set in an interval $[a,b]$, and $\tau: S \longrightarrow [c,d] $ an order preserving injective function. Assume also that $\tau (S)$ is dense in $[c,d]$. Then, $\tau$ has a unique  extension to a homeomorphism $\tau:[a,b] \longrightarrow [c,d]$.
\end{lemma}

The map $\Theta_c: \cala(K_c^T) \longrightarrow \cala(K_\la)$ satisfies all the hypothesis of Lemma \ref{extension}. Abusing notation, we obtain the unique extension
\[ 
	\Theta_c: [\thetatil^q,\thetatil^0] \longrightarrow [0,1]  \quad \text{or} \quad  \Theta_c: \bt \longrightarrow \bt, 
\] 
which is an orientation preserving homeomorphism. By definition, this map satisfies the second statement of Prop.~\ref{Bigtheta}.

To prove the last statement first suppose that $\theta$ belongs to $\cala(K_c^T)$,  $z=R^\ast_c(\theta)\in K_c^T$. Then $\Theta_c(m_2^{(1)}(\theta))$ is the argument of a ray landing at $H_c(f_c\one(z))\in K_\la$. On the other hand $m_{q+1} (\Theta_c(\theta))$ is the argument of a ray landing at $P_\la(H_c(z))$. Since $H_c$ conjugates $f\one$ to $P_\la$, these two points are the same. It follows that $\Theta_c$ conjugates $m_2^{(1)}$ to $m_{q+1}$ at least up to arguments of rays landing at the same point. Therefore, we must show that equivalence classes of rays landing at the same point cannot be permuted, but this follows from the fact that $\Theta_c$, $m_2^{(1)}$ and $m_{q+1}$ are order preserving for arguments that belong to the same $V_c^i$ (in particular, for rays in the same equivalence class) and from $\Theta_c(\thetatil)=0$.

Therefore, the  map $\Theta_c$  conjugates $m_2^{(1)}$ and $m_{q+1}$ on a dense set and hence everywhere. By Lemma \ref{expanding} such a homeomorphism is unique and therefore $\Theta_c=\bigthatp$. This ends the proof of Theorem \ref{Bigtheta}   
\end{proof}

We proceed now to combine the maps above to construct  
\[ 
	(\bigthatpp)^T: \bt\p^T@>{\bigthatp}>> \bt @>{\bigthat_{p^\prime/q}^{-1}}>>  \bt_{p^\prime/q}^T.
\]

The next proposition follows from Prop.~\ref{Bigtheta}:

\begin{proposition}
Let $c \in M_{p/q}$ and $c^\prime=\finho(c) \in M_{p^\prime/q}$. Then,
\begin{enumerate}

\item  A ray $R_c(\theta)$ lands at $z \in K_c^T$ if and only if the ray $R_{c^\prime}((\bigthatpp)^T(\theta))$ lands at $\widehat{\Phi}_c^T(z)$.

\item The map $(\bigthatpp)^T$ conjugates the first return maps on the arguments in $\bt\p^T$ and $\bt_{p^\prime/q}^T$ respectively, that is the following diagram commutes:
\[
	\begin{CD}
		\bt\p^T @>{m_{p/q,2}\one}>> \bt\p^T \\
		@V{(\bigthatpp)^T}VV  @VV{(\bigthatpp)^T}V    	\\
        \bt_{p'/q}^T @>{m_{p'/q,2}\one}>> \bt_{p^\prime/q}^T	 	\end{CD}	
\]   
\end{enumerate}
\end{proposition}

Finally, we extend the map $(\bigthatpp)^T$ to the rest of the arguments of the filled Julia set to obtain the map claimed in part two of Theorem C. Note that this is a somewhat ``artificial'' extension, since we will not use the dynamics of the polynomial. However, it is a nice complement to the map $\widehat{\Phi}_c$ on the filled Julia sets.

As before let $c \in M_{p/q}$ and $c^\prime=\finho(c) \in M_{p^\prime/q}$.

\begin{definition}{}
We define 
\[
	\bigthatpp: \bt \longrightarrow \bt
\]
by
\[
  \bigthatpp(\theta)= 
		\begin{cases}
    		(\bigthatpp)^T(\theta)		& \text{if $\theta \in \overset{\circ}{\bi}\p$}\\
  		(\bigthatpp)^T(\theta+\tfrac{1}{2})+\tfrac{1}{2} & \text{otherwise}
 		\end{cases}
\]
\end{definition}
Note that this map is well defined since the arguments of external rays in $\smash[t]{\overset{\circ}{V}}_c^0$ satisfy
\[ (\bigthatpp)^T(\theta) = (\bigthatpp)^T(\theta + \frac{1}{2})+\frac{1}{2}
\]
The map $\bigthatpp$ is an orientation preserving homeomorphism satisfying part two of Theorem C, that is:
\begin{proposition}
A ray $R_c(\theta)$ lands at $z\in K_c$ if and only if the ray $R_{c^\prime}(\bigthatpp(\theta))$ lands at $\widehat{\Phi}_c(z) \in K_{c^\prime}$.
\end{proposition}

\begin{remark}
We could have deduced the map on arguments directly from the homeomorphism $\widehat{\phi}_c: N(K_c) \longrightarrow N(K_{c'})$. Instead, we have once again applied the bridge to the higher degree polynomials in order to show that the map $\bigthatpp$ only depends on $p$, $p'$ and $q$.
\end{remark}

\begin{proposition}
Suppose $K_c$ (and therefore $K_{c^\prime}$) is locally connected. Then, there exists a homeomorphism $\widehat{\widehat{\phi}}_c$ which coincides with $\fihat_c$ on $K_c$ and which maps the external ray $R_c(\theta)$ to the external ray $R_{c'}(\bigthatpp(\theta))$ equipotentially.
\end{proposition}

\begin{remark}
If $K_c$ is not locally connected, we can still define the homeomorphism from $\bc \setminus K_c$ to $\bc \setminus K_{c'}$ as above, but this map may not match continuously with $\fihat_c$ restricted to $K_c$.
\end{remark}

%*****************************************************************************

\subsubsection{Parameter Plane} \label{param}

We are now ready to prove the combinatorial part of the parameter space theorems.

We start by proving the first part of Theorem H. Then we obtain the first part of Theorem C for rational arguments with even denominators. Finally, we extend it to all rational arguments.

\begin{proof}{of the first part of Theorem H}
Let $\Theta\p=\bigthat\p \circ m_2^{q-1} : (\theta\p^-, \theta\p^+) \longrightarrow (0,1)$ as defined earlier.

Suppose $c$ is a Misiurewicz point in $M\p$. Recall that
\[
	R_c(\theta)  \text{\ lands at } c \in K_c \Longleftrightarrow  R_M(\theta) \text{\  lands at } c \in M\p.
\]
 Moreover 
\[
R_c(\theta)  \text{\ lands at } c \in K_c \Longleftrightarrow  R_c(2^{q-1} \theta) \text{\  lands at } Q_c^{q-1}(c)=Q_c^q(0) \in K_c.
\]
Note that $Q_c^q(0)$ is the critical value of $f_c\three$, hence mapped by $H_c$ to the critical value $v_\la=\la (-1)^q(\frac{q}{q+1})^{q+1}$ of $P_\la$ where $\la=\ho(c)$.

It follows from the dynamical part of Theorem H that
\[
R_c(2^{q-1} \theta)  \text{\ lands at } Q_c^q(0) \in K_c \Longleftrightarrow  R_\la(\bigthat\p ((2^{q-1} \theta)) \text{\  lands at } v_\la \in K_\la.
\]
Finally, apply Prop.~\ref{misus} to obtain that
\[
R_\la(\bigthat\p ((2^{q-1} \theta)) \text{\  lands at } v_\la \in K_\la
\Longleftrightarrow  R_{L\qo}(\bigthat\p ((2^{q-1} \theta)) \text{\  lands at } \la \in L\qo.
\] \ 
\end{proof}

Define 
\[
	\finte: [\theta\p^-, \theta\p^+] @>{\Theta_{p/q}}>> [0,1] @>{\Theta_{p'/q}^{-1}}>> [\theta_{p'/q}^-, \theta_{p'/q}^+].
\]
The first part of Theorem C concerning rational arguments with even denominators, that is arguments of rays landing at Misiurewicz points,  follows immediately from Theorem H.

The missing part is to prove part one of Theorem C for rational arguments with odd denominators, that is arguments of rays landing at roots of hyperbolic components.

\begin{remark}
From Cor.~\ref{rationals} it follows that $\finte$ maps rational angles with even (resp.~odd) denominator to rational angles with even (resp.~odd) denominators.
\end{remark}

\begin{lemma}
Let $c_0 \in M\p$ be the root of a hyperbolic component and $\theta_0 \in \bt$ the argument of one of the rays landing at $c_0$. Then, there exists a sequence of Misiurewicz points $c_n \rightarrow c_0$ and a sequence of arguments $\theta_n \rightarrow \theta_0$ such that $R_M^\ast(\theta_n)=c_n$.
\end{lemma}
For a proof, see \cite{hubbard,dan}.

Fix $c_0$, $\theta_0$, $c_n$, and $\theta_n$ as in the lemma. Let $c'_0=\finho(c_0)$, $c'_n=\finho(c_n)$, $\theta'_0=\finte(\theta_0)$ and $\theta'_n=\finte(\theta_n)$. By the remark above, $\theta'_0$ is rational with odd denominator. By continuity of $\finte$ and $\finho$ we have that $\theta'_n \rightarrow \theta'_0$ and $c'_n \rightarrow c'_0$. By remarks \ref{conju} and \ref{roots}, $c'_n$ are Misiurewicz points and $c'_0$ is the root of a hyperbolic component. It follows immediately from above that $R^\ast_M(\theta'_n)=c'_n$.

To conclude the proof of Theorem C we must show
\begin{equation} \label{conclude}
	R_M^\ast(\theta'_0)=c'_0.
\end{equation}

We recall the definition of the {\em impression} of an external ray of a compact, connected, full subset $K \subset \bc$ (see \cite{dan}).
\begin{definition}{}
For any $\epsilon >0$ and any $\theta \in \bt$ let $I_{K,\epsilon} (\theta)$ denote the closed subset of $\partial K$ given by
\[	
	I_{K,\epsilon} (\theta)= \overline{R_K([\theta -\epsilon, \theta+\epsilon ])} \setminus R_K([\theta -\epsilon, \theta+\epsilon ]).
\]
The {\em impression} $I_K(\theta)$ is defined as 
\[	I_K(\theta)=\bigcap_{\epsilon>0} I_{K,\epsilon} (\theta).
\]
\end{definition}

\begin{lemma}
$c'_0 \in I_M(\theta'_0)$.
\end{lemma}

\begin{proof}{}
For each $\epsilon > 0$, there exists $n_\epsilon \in \bn$ such that $\theta'_n \in [\theta'_0 -\epsilon, \theta'_0 +\epsilon ]$ for all $n \geq n_\epsilon$. It follows that $c'_n \in I_{M,\epsilon} (\theta'_0)$ for all $n \geq n_\epsilon$, therefore that $c'_0 \in I_{M,\epsilon} (\theta'_0)$ for all $\epsilon >0$. Hence $c'_0 \in I_M(\theta'_0)$.
\end{proof}

The conclusion will follow from the following lemma.

\begin{lemma}
Let $c \in M$ be the root of a hyperbolic component $\Omega$ and $\theta$ a rational angle with odd denominator. If $c \in I_M(\theta)$ then $I_M(\theta)=\{ c \}$.
\end{lemma}
\begin{proof}{}
There are two cases, depending on $\Omega$ being a satellite component or a primitive component. The proof in the satellite case and ``half'' of the primitive case can be found in \cite{dan}. The missing part is the following, using a technique from \cite{dh2}, part II.

Suppose $\Omega$ (different from $\Omega_0$)  is primitive and $\theta^\pm$ are the rational arguments of the rays landing at $c$. Then, there exists a sequence of Misiurewicz points on the combinatorial vein in $M$ between $0$ and $c$, and sequences of arguments $\theta'_n$ and $\theta''_n$, where $R_M^\ast(\theta'_n)=R_M^\ast(\theta''_n)=c_n$ and $\theta'_n \nearrow \theta^-$ and $\theta''_n \searrow \theta^+$. The existence is proved by modifying the Hubbard tree for the center of $\Omega$ to become the Hubbard tree of Misiurewicz points as described, transferring arguments from the dynamical plane to the parameter plane as usual. (See also \cite{lavaurs,alfredo,tanlei}).
\end{proof}

With this lemma, (\ref{conclude}) follows immediately.

\begin{remark} \label{stronger}
The result above can be carried over by the same technique to the set of irrational arguments that correspond to rays that land at c-values on boundaries of hyperbolic components and also to c-values described by  Yoccoz para-puzzles.
\end{remark}

%*****************************************************************************

\subsection{Proof of Theorem D}\label{theoremD} 

In this section we combine the symmetries of the 0-wakes $WL\qo$ with respect to complex conjugation and Theorem C to obtain a non trivial homeomorphism from the set of arguments of rays in a $p/q$-wake $WM\p$ to itself. This map is the combinatorial analogue of  the involutions $\ical_{p/q}$ in Theorem B.

Recall from section \ref{involutions} that the map
\[ 
	\begin{array}{rccc}  
		C : & WL_{q,0}  	& \longrightarrow 	& WL_{q,0} 		\\
        		& \la 		& \longmapsto		& \overline{\la}
	\end{array}
\] 
is a homeomorphism mapping the ray $R_L(\theta)$ to the ray $R_L(1-\theta)$ for $\theta \in (0,1)$. We define the orientation reversing homeomorphism from $[\theta\p^-, \theta\p^+]$ onto itself by
\[
	\begin{array}{rccccccc}
		\linte\p: & [\theta^-\p,\theta^+\p] & @>{\Theta\p}>> &[0,1]& @>>> & [0,1] & @>{\Theta\p^{-1}}>> & [\theta^-\p,\theta^+\p] \\
	          		&   \theta    & \mapsto	&     \Theta\p(\theta) & \mapsto   &   1-\Theta\p(\theta) & \mapsto & \Theta\p^{-1}(1-\Theta\p(\theta))
	\end{array}
\]

Note that this is an orientation reversing homeomorphism from the arguments of the rays in $WM_{p/q}$ to themselves. It also satisfies $\overline{\Theta}_{p/q}^2=Id$. In particular,
the arguments that bound the wake are mapped to each other.

The fixed point of $\overline{\Theta}_{p/q}$ is $\theta^s:=\Theta_{p/q}^{-1}(1/2)$,  hence rational. 

This proves properties 1 and 2. It is clear that by the properties of $\Theta\p$, property 3 is satisfied, i.e. 
\[
	R_M^\ast(\theta)=c \Longleftrightarrow R_M^\ast(\linte\p(\theta))=\icalp(c).
\]

For property 4 in dynamical plane we define
\[ \begin{array}{rccccccc} 
\smash[t]{\widehat{\linte}}\p^T:&I_{p/q} &\stackrel{\bigthatp}{\longrightarrow} &[0,1] &\longrightarrow &[0,1] &\stackrel{\bigthatp^{-1}}{\longrightarrow} & I_{p/q}\\
 & \theta & \longmapsto & \bigthatp(\theta) & \longmapsto &1- \bigthatp(\theta) & \longmapsto & \bigthatp^{-1}(1- \bigthatp(\theta))
\end{array}
\]

Then, as we did in section \ref{dynplane} to define $\finte$, we use the symmetries of $K_c$ with respect to the origin to extend $\smash[t]{\widehat{\linte}}\p^T$ to an orientation reversing homeomorphism $\widehat{\linte}\p:\bt \longrightarrow \bt$ satisfying the required properties.
  
This concludes the proof of Theorem D.

%*****************************************************************************
%*****************************************************************************
\subsection{Proof of Theorem E} \label{theoremE}

It is clear that we can use the map $\finte$ in Theorem C to define a homeomorphism from $WM\p \setminus M\p$ to $WM_{p'/q} \setminus M_{p'/q}$ by mapping a ray of argument $\theta$ to the ray of argument $\finte(\theta)$ equipotentially. By the properties of $\finte$, this homeomorphism would match with $\finho$ at all the landing points of rays with rational arguments (together with the landing points of rays with irrational arguments covered by remark 7.16).

Now suppose the Mandelbrot set were locally connected. It would follow that all external rays of $M$ land continuously and therefore, the radial extension would match continuously with $\finho$ at all points in the boundary of $M\p$.

The extension of $\ical\p$ to the $p/q$-wake would follow immediately by composition of the extensions of $\Phi_{p(q-p)}^q$ and complex conjugation.

This concludes the proof of Theorem E.

%*****************************************************************************
%*****************************************************************************
%*****************************************************************************

%\bibliographystyle{apanamed}
%\bibliography{/usr/users/guest/fagella/Paper/math}

\begin{thebibliography}{AB9}

\bibitem[A]{atela} P.~Atela, Bifurcations of Dynamic Rays in Complex Polynomials of Degree Two, {\em Erg.~Th.~\& Dyn.~Sys.} {\bf 12} (1991) 401-423.

\bibitem[Bl]{Blanchard} P.Blanchard, Complex Analytic Dynamics on the Riemann Sphere, {\em Bull.~Amer.~Math.~Soc.} {\bf 11} (1984), 85-141.

\bibitem[Br]{bodil} B.~Branner, The Mandelbrot Set, {\em Proc.~Symp.~Applied Math.}, (1989), 75-105. 

\bibitem[BD]{Branner} B. Branner \& A. Douady, Surgery on Complex
Polynomials,   Proceedings of the Symposium on Dynamical Systems, Mexico,
1986, {\em Lecture Notes in Math.}, {\bf 1345} 11-72, Springer.

\bibitem[BH]{bh} B. Branner \& J. H. Hubbard, The Iteration of Cubic
Polynomials. Part I: The Global Topology of Parameter Space,   {\em Acta
Mathematica}, {\bf 160}, (1988), 143-206.

\bibitem[D]{douady} A.~Douady, Chirugie sur les Applications Holomorphes, {\em Proc.~Int.~Congr.~Math.}, Berkeley, 1986, 724-738.

\bibitem[DGH]{dgh} R. Devaney, L.Goldberg \& J.Hubbard, A dynamical
approximation to the Exponential Map by polynomials, {\em Preprint}.

\bibitem[DH1]{dh} A.Douady \& J.Hubbard, Iterations des Polynomes
Quadratiques Complexes, {\em C.R. Acad. Sci., Paris, Ser.I} {\bf 29}
(1982), 123-126.

\bibitem[DH2]{dh2} A.Douady \& J.Hubbard, Etude dynamique des polynomes
complexes, I,II, {\em Publ. Math. Orsay} (1984, 1985).

\bibitem[DH3]{dh3} A. Douady \& J.Hubbard,  On the Dynamics of
Polynomial-like Mappings, {\em Ann. Scient., Ec. Norm. Sup. $4^e$
series}, {\bf 18} (1985) 287-343.

\bibitem[F]{me} N.~Fagella, Limiting Dynamics of the Complex Standard Family, {\em Int.~J.~of Bif.~and Chaos} (1995), To appear.

\bibitem[GM]{gm} L.~Goldberg \& J.~Milnor, Fixed Points of Polynomial Maps II: Fixed Point Portraits, {\em Ann.~Sci.~Ec.~Norm.~Super.} {\bf 26} (1993) 51-98.

\bibitem[H]{hubbard} J.~H.~Hubbard, Local Connectivity of Julia Sets and Bif.~Loci: Three Theorems of J.~C.~Yoccoz, {\em Topological Methods in Modern Mathematics}, Publish or Perrish, Inc.~(1993) 467-511.

\bibitem[L]{lavaurs} P.~Lavaurs, Syst\'emes Dynamiques Holomorphes: Explosion de Points P\'eriodiques Paraboliques, {\em These}, Univ.~Paris Sud Orsay (1989).

\bibitem[LS]{dierk} E.~Lau \& D.~Schleicher, Internal Adresses in the Mandelbrot Set and Irreducibility of Polynomials, {\em Prepint} n.~1994/19, SUNY Stony Brook.

\bibitem[Mc]{mcmullen} C.~T.~McMullen, Complex Dynamics and Renormalization, {\em Annals of Mathematics Studies}, 135, Princeton University Press, 1994.

\bibitem[Mi]{milnor} J.~Milnor, Dynamics on One Complex Variable: Introductory Lectures, IMS at Stony Brook Preprint n.~1990/5 (1992).

\bibitem[MSS]{sullivan} R. Ma\~{n}e, P. Sad \& D. Sullivan,  On the Dynamics of Rational Maps, {\em Ann. Scie.  \'{E}cole Norm. Sup. } (4) {\bf 16} (1983), 193-217.

\bibitem[Pe]{carsten} C.~L.~Petersen, On the Pommerenke-Levin-Yoccoz Inequality, {\em Erg.~Th~\& Dyn.~Sys.} {\bf 13} (1993) 785-806. 

\bibitem[Poi]{alfredo} A.~Poirier, On Postcritically Finite Polynomials. Part One: Critical Portraits, {\em Prepint} n.~1993/5, SUNY Stony Brook. 

\bibitem[Pom]{pomerenke} Ch.~Pommerenke, Boundary Behavior of Conformal Maps, GMV 299, Springer 1992.

\bibitem[S]{dan} D.~S{\o}rensen, Complex Dynamical Systems: Rays and Non-Local Connectivity, {\em Ph.D.~Thesis}, 1994, Tech.~Univ.~of Denmark. 

\bibitem[T]{tanlei} Tan Lei, Voisinages Connexes des Points de Misiurewicz, {\em Ann.~Inst.~Fourier} {\bf 42, 4} (1992) 707-735
\end{thebibliography}
%\end{document}

\end{document}